\documentclass [12pt] {article}
\title { Knot polynomials via one parameter knot theory}
\author{Thomas Fiedler}
\usepackage{amssymb,amsfonts,epsfig}
\begin{document}
\newtheorem{proposition}{Proposition}
\newtheorem{theorem}{Theorem}
\newtheorem{lemma}{Lemma}
\newtheorem{corollary}{Corollary}
\newtheorem{example}{Example}
\newtheorem{remark}{Remark}
\newtheorem{definition}{Definition}
\newtheorem{question}{Question}
\maketitle
\begin{abstract}
We construct new knot polynomials.

Let $V$ be the standard solid torus in 3-space and let $pr$ be its standard projection onto an annulus.
Let $M$ be the space of all smooth oriented knots in $V$ such that the restriction of $pr$ is an immersion
(e.g.  regular diagrams of a classical knot in the complement of its meridian).

There is a canonical one dimensional homology class for each connected component of $M$.

We construct homomorphisms from the first homology group of $M$  into rings 
of Laurent polynomials. Each such homomorphism 
applied to the canonical homology class gives a knot invariant.

Let $\gamma$ be a generic smooth oriented loop in $M$ (i.e. a one parameter family of knot diagrams in the annulus).
For finitely many points in $\gamma$ the corresponding knot diagram has in the projection $pr$ an ordinary triple point
or an ordinary auto-tangency. To each such diagram we associate some Laurent polynomial by using extensions of
the Kauffman 
bracket or of the Kauffman state model  for the Alexander polynomial. We take then an 
algebraic sum of these polynomials over all triple points and all autotangencies in $\gamma$. The resulting polynomial 
depends only on the homology class of $\gamma$ if and only if it verifies two sorts of equations: the tetrahedron equations
and the cube equations. We have found five different non trivial solutions.
\footnote{
Knot polynomials, tetrahedron equation, cube equation.}
\end{abstract}
\tableofcontents
\section{Introduction}
We work in the smooth category and we use the standard orientation conventions. We suppose that the reader is familiar
with the state models of Kauffman for the Jones and for the Alexander polynomial.

Knots in 3-space can be given by diagrams in the plane. A knot polynomial is usually a polynomial which is defined
by using a generic diagram (i.e. the projection in the plane has as singularities only ordinary double points). In order
to be a knot invariant the polynomial has to be invariant under the Reidemeister moves of type III and type II . 
Reidemeister moves of type I multiply the polynomial usually by some well defined factor. Quantum knot invariants 
are knot polynomials which are defined by using solutions of the Yang-Baxter equation. The Yang-Baxter equation is 
an equation on the level of operators which is associated to Reidemeister III moves (see \cite{Tu}, \cite{Jo} and e.g. \cite{K3}). 
If one considers oriented diagrams 
then there are eight different (local) types of Reidemeister III moves and four different types of Reidemeister II moves
(see e.g. \cite{AF}).
Let us call the {\em positive Reidemeister III move} those in which all three involved crossings are positive.
Fortunately, it turns out that in order to check that a polynomial is a knot invariant it suffices to check invariance under 
only the positive Reidemeister III move and certain two of the types of Reidemeister II moves (see e.g. Section 1 in \cite{F1}).
Let $M$ be the space of all regular knot diagrams as introduced in the abstract (see also Subsection 2.1.). Quantum 
knot invariants can be seen as certain homomorphisms from $H_0(M)$ to some rings of Laurent polynomials.

In this paper we introduce a different approach in order to define new knot polynomials. First of all we replace the 3-sphere 
by the standard solid torus $V$ in $S^3$. Indeed, two oriented knots in the 3-sphere are isotopic if and only if the two oriented links 
consisting of the knot and its positive meridian are isotopic. In fact, we can identify the meridians of the two knots 
and keep them fixed in the whole isotopy. The meridian is a trivial knot and its complement in the 3-sphere is a standard 
solid torus. Consequently, the knots in the 3-sphere are isotopic if and only if the corresponding knots in the standard solid torus
are isotopic. Moreover, this remains true if we replace the two framed knots in the solid torus by the same satellite.
Consequently, the question about isotopy of classical knots can be reformulated for knots in the solid torus which
represent arbitrary homology classes in $H_1(V)$.

Why are we doing this? As well known, the identity component of the diffeomorphism group of the 3-sphere 
retracts by deformation onto $SO(3)$ (see \cite{H1}). The identity component of the diffeomorphism group of $V$ retracts 
by deformation onto $SO(2)\times SO(2)$. So, the dimension becomes smaller but the first homology becomes much 
bigger. The latter is in fact generated by the rotations (from the angle 0 to $2\pi$) of the solid torus along its core and 
by those around its core. These rotations act on the knot diagrams and hence define loops in $M$. We want to define
homomorphisms from $H_1(M)$ to (torsion free) polynomial rings. Hence, we need that our knot space has non trivial
first Betti number.

Moreover, a famous theorem of Waldhausen says that knots in 3-space are classified by their peripheral system (\cite{W}). This
suggests that for the construction of knot invariants in 3-space it should be very usefull to consider framed knots in the complement of 
their meridian. This is exactly what we are doing.

  The rotations along the core 
are not interesting: they correspond just to rotating the diagrams together with the annulus. However, the rotations around 
the core change the knot diagrams in a non trivial way. We call this the {\em canonical loop}. Besides the case of closed braids
the canonical loop will always contain Reidemeister I moves. Hence, it does not yet define a loop in $M$. But we can 
approximate the canonical loop by a loop in $M$ using the Whitney trick. It turns out that the homology class of the
approximating loop, called $[rot]$, is unique up to adding those homology classes which correspond to sliding a small curl (created 
by a Reidemeister I move) once along the whole knot in the solid torus. We call the latter homology classes in $H_1(M)$
the {\em sliding classes}. Sliding classes are always non trivial if the knot is homological non trivial in $V$. Let the 
generic diagram $D$ be a point in $M$.
Let $w(D)$ be the writhe of $D$ and let $n(D)$ be the {\em Whitney index of D in the annulus}. ( We smooth all double points
of $pr(D)$ with respect to the orientation of $D$. The Whitney index of D in the annulus $n(D)$ is then defined as the number positive {\em contractible}
 Seifert circles
minus the number of negative {\em contractible} Seifert circles. For example, if $D$ is the diagram of a closed braid
then $n(D) = 0$.) The component
of $M$ which contains $D$ is completely characterized by the knot type in $V$ represented by $D$ together with the two
integers $w(D)$ and $n(D)$.

As a summary so far: for each connected component $M_c$ of $M$ we have a homology class 
$[rot_c] \in H_1(M_c)$ which is well defined modulo sliding classes. We call $[rot_c]$ modulo sliding classes 
the {\em canonical class} of $M_c$. The component 
$M_c$ (and hence the canonical class) is determined by the knot type together with two easily calculable integers.

The topology and in particular the homology of knot spaces is studied in \cite{H2},  \cite{BC} and  \cite{Bu}. In this 
paper we are not interested in the homology of the (disconnected space) $M$ but in its cohomology!

In order to define new knot invariants it is of crucial importance that each component of $M$ contains a distinguished
one dimensional integer homology class (called the canonical class). We evaluate our cohomology classes (with 
values in polynomial rings) on this distinguished homology class.

Instead of homomorphisms from $H_0(M)$ into polynomial rings we will construct homomorphisms  from $H_1(M)$ into
polynomial rings.  Each such homomorphism $\rho$ gives rise to a serie of knot polynomials, indexed by the
two integers $w$ and $n$, provided that $\rho$ is trivial on all sliding classes. The question is now: how to obtain such
homomorphisms?

To that goal we have to study loops in $M$. This is called {\em one parameter knot theory}. The basic notions and facts
of one parameter knot theory are worked out in our joint paper with Vitaliy Kurlin \cite{FK}. It contains in particular
a {\em higher order Reidemeister theorem}. The solid torus $V$ is naturally fibered over the circle. The fibers are 
the planar discs which are orthogonal to the core of $V$.
$M$ contains the discriminant $\Sigma$ of all non generic diagrams.
$\Sigma$ has a natural stratification. The strata of codimension 1 correspond to diagrams which have exactly
one ordinary triple point in $pr$ (sometimes called {\em triple crossing}) or one ordinary autotangency in $pr$.
The higher order Reidemeister theorem says the following: two generic loops 
in $M$ are homologic if and only if the two families of diagrams can be transformed into each other by isotopies of 
families of diagrams and finitely many {\em moves} of certain types. These moves correspond to passing generically 
through strata of codimension 2 of $\Sigma$, to touching generically strata of codimension 1 of $\Sigma$
 and to Morse modifications of the loops in $M$. It turns out that the following four types of strata of codimension
2 are important:

(1) diagrams with an ordinary quadruple point in $pr$. We denote the union of these strata by $\Sigma^{(2)}_q$.

(2) diagrams with an ordinary autotangency through which passes transversally another branch in $pr$. We denote 
the union of these strata by $\Sigma^{(2)}_{a-t}$.

(3) transverse intersections of two strata of codimension 1, i.e. diagrams which have exactly two ordinary
triple points in $pr$ and so on.

(4) diagrams with an autotangency in an ordinary flex. We denote the union of these strata by $\Sigma^{(2)}_{f}$.

Our strategy is now the following: let $\gamma$ be an oriented generic loop in $M$. Let $\Sigma^{(i)}$ denote
the union of all strata of codimension i. $\gamma$ intersects $\Sigma^{(1)}$ transversally in a finite number of points.
We consider only those intersection points which correspond to a diagram which has exactly one ordinary triple
 point or one
ordinary autotangency in $pr$. The union of the strata which correspond to diagrams with an ordinary triple point is 
denoted by $\Sigma^{(1)}_t$ and those which correspond to diagrams with an ordinary autotangency by 
$\Sigma^{(1)}_a$. To each of these diagrams we want to associate some Laurent polynomial which should 
be a knot invariant relative to the singularity. Lets concentrate for the moment only on triple points.
{\em Relative invariance with respect to the triple point} means that the polynomial is invariant under 
all regular isotopies of the knot which preserve the ordinary triple point and such that in the isotopy no other branch of 
the diagram  moves through 
the triple point in $pr$. We want to define such polynomials with {\em state models}. The strata of codimension 2 
from (3) 
force us to use state models which are invariant under Reidemeister III and II moves outside of the triple point.
In this paper we will use the {\em Kauffman bracket} \cite{K2} and the {\em Kauffman state model for the 
Alexander polynomial} \cite{K1}.
These state models are defined for diagrams which have only ordinary double points as singularities in $pr$.
However, we have a triple point too and we have to extend the definitions of the state models at the triple point.
As already mentioned , there are eight different types of triple points. Let us consider the positive triple point and the
case of the Kauffman bracket. Let us consider a small disc in the annulus centered at the triple point. Its boundary
intersects the diagram in six points. We replace now the triple point by all possible chord diagrams in the disc 
which connect the six points on the boundary. We call this the {\em simplification} of the triple point.
( This is the analogue of the {\em A-smoothing and $A^{-1}$-smoothing} of
double points
in the definition of the Kauffman bracket.) However, for a triple point we have fifteen different chord diagrams. This
allows us to introduce fifteen new variables. In the case of the state model for the Alexander polynomial we have to put
two dots in exactly two of the six regions of the complement of the diagram in the disc. We call this the {\em marking}
of the triple point. There are again fifteen 
possibilities to do this and hence fifteen new variables. For each of the new variables we define the polynomial now 
almost as usual.

 In the case of the state model for the Alexander polynomial we have to choose two {\em *-regions},
i.e. two regions of the complement of the diagram which will never contain a dot. The choice has to be invariant under 
Reidemeister III and II moves outside of the triple point. Luckily, we can do this by choosing the two regions which
are adjacent to the boundary components of the annulus. This is a correct choice if the diagram of the knot is not 
contained in a disc. But this was achieved from the very beginning when we replaced knots in 3-space by knots in 
the solid torus. Notice, that the two *-regions need not to be adjacent, in difference with the original state model 
constructed by Kauffman.

In the case of the Kauffman bracket there are no longer all circles embedded in the plane. There can be one, two or 
three double points. We take an arbitrary resolution of the double points (i.e. we separate abstractly the two branches)
and we count the number of circles now as usual in the Kauffman bracket (i.e. each contractible circle in the annulus
counts as a factor $-A^2 -A^{-2}$ and each non contractible circle counts as a factor $h$, where $h$ is an independent 
variable, see \cite{HP}).

As a summary so far: for an oriented knot diagram in the annulus with a triple point $p$ we have constructed 
(in each of the two 
cases, Kauffman bracket or Kauffman state model for the Alexander polynomial) a Laurent polynomial of the form 
$L_p = C_1P_1 + C_2P_2 +... +C_{15}P_{15}$. Here $C_1$ up to $C_{15}$ are independent variables, which
depend only on the type of the triple point $p$ (remember that there are exactly eight types) and of its 
simplification or its marking and each $P_i$ is a 
Laurent polynomial of either one variable $A$ ({\em Alexanders case}) or of two variables $A,h$ ({\em Jones case}).

We want to associate to the loop $\gamma$ now the sum of the polynomials $L_p$ over all triple points in the loop.
In a homotopy the loop $\gamma$ can become tangential to a stratum of $\Sigma^{(1)}$ corresponding to an ordinary
triple point. Consequently,  triple points can appear or disappear in pairs in the loop. Therefore we have to define an 
intersection index of $\Sigma^{(1)}$ with $\gamma$ at each triple point. For this purpose we have to define a 
coorientation for each stratum of $\Sigma^{(1)}$ which corresponds to a triple point. It turns out, that this can be done 
consistently, i.e. the closure in $M$ of the union of all strata of $\Sigma^{(1)}$ which correspond to triple points
becomes an integer cycle of codimension 1 in $M$. Let $sign(p)$ denote the sign of the intersection of the oriented loop
$\gamma$ with the cooriented $\Sigma^{(1)}$ at the diagram with the triple point $p$.

We define now

$L(\gamma) = \sum_p{sign(p) L_p}$  

where the sum is over all diagrams with triple points $p$ in the loop $\gamma$.

It follows immediately from its construction that $L(\gamma)$ is invariant under a homotopy of $\gamma$  
which passes through the transverse intersection of two strata from $\Sigma^{(1)}$ (see (3) above), or
which touches a stratum of $\Sigma^{(1)}$. On easily sees that $L(\gamma)$ is invariant under Morse modifications
of $\gamma$, because we can assume that the singularity of the Morse modification does not correspond to
a diagram which has a triple point or an autotangency in $pr$. Hence $L(\gamma)$ is unchanged when we take
the sum over all components of $\gamma$.

The main body of our work is to achieve invariance of $L(\gamma)$ under homotopies which passe through 
strata of (1) or (2). 

At a stratum of a quadruple point there are four strata of triple points which intersect mutually transverse.
Let $s$ be the boundary of a small normal disc in $M$ for a point in $\Sigma^{(2)}_q$. The loop $s$ intersects
$\Sigma^{(1)}_t$ in exactly eight points and it does not intersect other strata from $\Sigma^{(1)}$. We number 
these intersection points. The situation is illustrated in Fig.1 which shows the normal disc and its intersection with
\begin{figure}
\centering 
\psfig{file=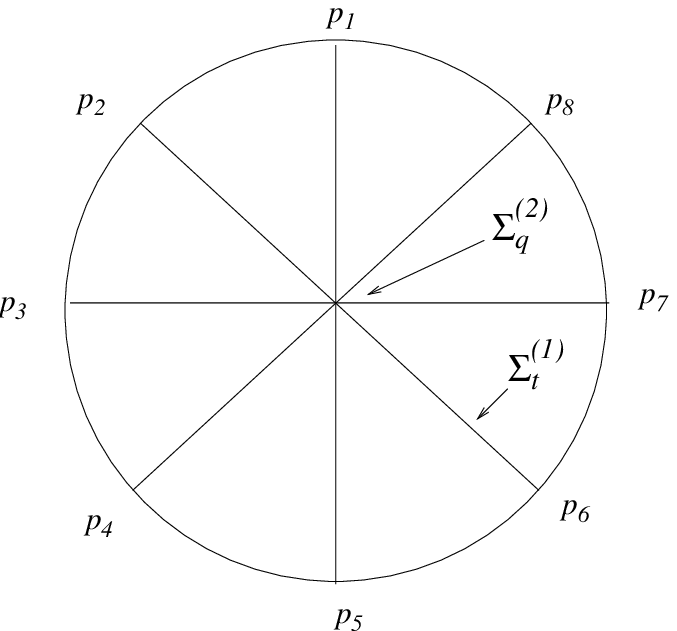}
\caption{}
\end{figure}
$\Sigma$.
 $L(\gamma)$ is invariant under a generic homotopy of $\gamma$ which passes through $\Sigma^{(2)}_q$
if and only if 

$L(s) = \sum_i sign(p_i)L_{p_i} = 0$.  (a)

We call this condition a {\em tetrahedron equation}. Indeed, let us consider
 diagrams which are colored pure closed braids.
In analogy with the quantum group approach we associate to 
each triple point in a diagram of a colored pure closed braid
some matrix. This matrix should depend only on the three colors of the strands in the triple point and on the 
cyclical order of their images in the annulus. Let $s$ be a meridian of $\Sigma^{(2)}_q$.  Let $A_i$ be the matrix associated to the
 triple crossing $p_i$, $i \in \{1,...,8\}$ in $s$. In this case the tetrahedron equation would be the following equation:

$\Pi_i A^{sign(i)}_i = Id$.  (b)

We expect that this equation should be closely related to the {\em Knizhnik-Zamolodchikov} 
equations  \cite{EFK}. If one associates to a one parameter
family of diagrams its {\em trace graph} then one can really see a tetrahedron associated to $\Sigma^{(2)}_q$ (see 
\cite{FK}).

However, in this paper we do not solve the operator equation (b). Instead of a matrix we associate to a diagram 
with a triple 
crossing directly some Laurent polynomial, as described above. This approach has the advantage that it works
for {\em all} knots, not just for colored pure closed braids. 

The tetrahedron equation becomes then the
 equation $L(s) = 0$ (a).

In the case of $\Sigma^{(2)}_{a-t}$ the situation is a bit simpler. A normal disc for $\Sigma^{(2)}_{a-t}$ and its
intersection with $\Sigma$ is shown in Fig.2. The loop $s$ contains exactly two diagrams with a triple crossing
\begin{figure}
\centering 
\psfig{file=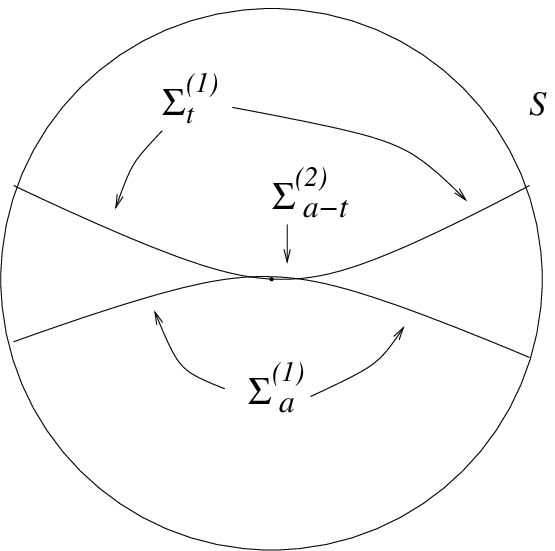}
\caption{}
\end{figure}
and two diagrams with an autotangency. Moreover, the two triple crossings are always from different types.
We have to solve the much simpler equation $L(s) = 0$.

The problem now is the following: there are exactly 48 types of strata of $\Sigma^{(2)}_{q}$ and 24 types of
strata of $\Sigma^{(2)}_{a-t}$. Let us call the {\em positive quadruple point} those for which all six involved
crossings are positive (and hence, all four involved triple crossings are positive too). Fortunately, an analogous 
phenomen as in the quantum group approach occurs: if $L(s) = 0$ for the positive quadruple point and for certain 
12 types
of autotangencies with a transverse branch then $L(s) =0$ for {\em all}  strata of $\Sigma^{(2)}_{q}$ and
$\Sigma^{(2)}_{a-t}$.

As a summary so far: in order that $L(\gamma)$ depends only on the homology class of $\gamma$ we have to 
show that $L(s) =0$ for the meridian $s$ of the positive quadruple point and for the meridians $s$ of 12 types of 
autotangecies with a transverse branch.
Here we have to work hard.

Let us consider the positive quadruple point in Jones case. For each of the eight triple crossings in $s$ we have 
to consider each of the fifteen simplifications. Let $T_4$ be the free $\mathbb{Z}[A,A^{-1}]$-modul generated by all 
chord diagrams 
with exactly four chords in the disc with distinguished end points and which have no more than three
 double points (and such that this number is minimal).
To each of the simplifications $C_i(p)$ for each of the triple crossings $p$ we associate an element $t_i(p)$ in $T_4$.

The tetrahedron equation $L(s) =0$ reduces then to the following equation:

$\sum_p\sum_i  t_i(p)C_i  =0$  (c).

It turns out that the left hand side of (c) contains  960 summands in $T_4$.

Alexanders case is even more complex.

For the convenience of the reader who wants to verify our results, we give a  part of the calculations of (c) 
in Subsections 4.2 respectively 4.3.

Let us consider now the boundary $s$ of a normal disc to a stratum of $\Sigma^{(2)}_{a-t}$. We replace the modul $T_4$
by its analoge $T_3$ with three chords. There are only two possible simplifications of an autotangency. We call
them $C_0$ and $C_x$ and show them in Fig.3. (The evident third possibility leads only to a summand which
\begin{figure}
\centering 
\psfig{file=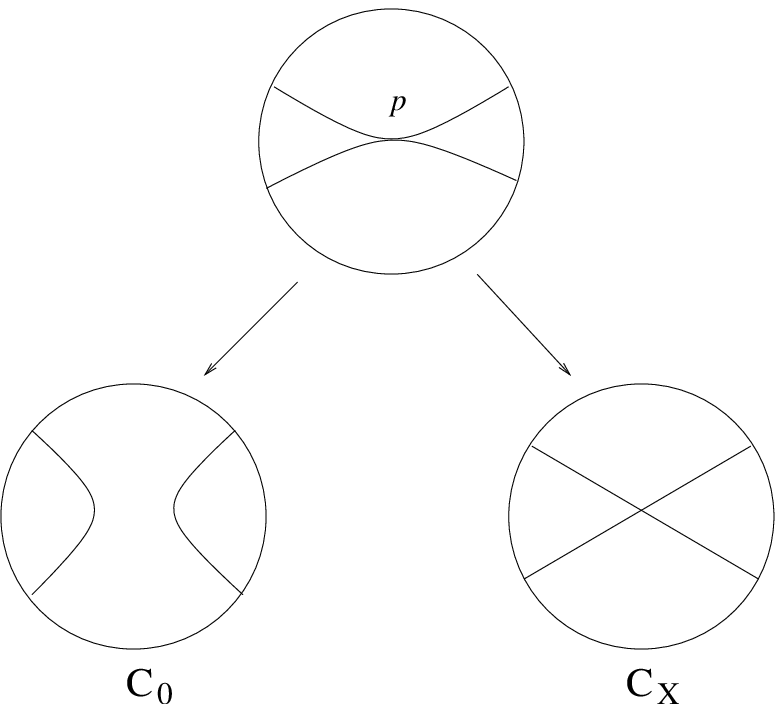}
\caption{}
\end{figure}
is a multiple of the Jones polynomial.) There are four types of autotangencies and consequently, there are four  
 independent variables. Autotangencies in an ordinary flex give two relations for these four variables. Let $p$ and $p'$ be the triple crossings  and let $r$ and $r'$ be the 
autotangencies in $s$. To each of the simplifications of triple crossings and of autotangencies  we associate 
an element in $T_3$, also denoted by $t_i(p)$.
The equation $L(s) =0$ leads now to an equation of the following form:

$\sum_i t_i(p) C_i + \sum_i t_i(p') C'_i + t_1(r) C_0 + t_2(r) C_x + t_1(r') C_0 + t_2(r') C_x =0$  (d).

Here, all $C_i$ and $C'_i$ from triple points  are independent variables, because the triple points are of different types.
This equation has 34 summands on the left hand side and we have to consider 12 such equations.

It turns out that the variables $C_i$ and $C'_i$ from triple points determine each other and that the variables 
$C_0$ and $C_x$ from autotangencies
become combinations of the variables $C_i$.

We are not yet done, because the equations (d) for different types of strata in $\Sigma^{(2)}_{a-t}$ are not independent!

Let $G$  be the graph which is constructed in the following way: the eight vertices correspond to the eight types
of triple crossings. Two vertices are connected by an edge if and only if the two triple crossings can come together in
a stratum of $\Sigma^{(2)}_{a-t}$. Hence, the edges correspond exactly to the 24 types of strata in $\Sigma^{(2)}_{a-t}$.
It turns out that $G$ is a cube, where each edge has to be doubled. Each of the 2-faces of the cube and each
doubled edge gives now
a relation for the set of variables $C_i$. We call these relations the {\em cube equations}. It turns out that there are exactly
nine independent cube equations. Notice, that autotangencies enter only in the cube equations and not in the
tetrahedron equation.

Let us summarize now: each solution of the positive tetrahedron equation produces automatically a solution of the
12 equations coming from meridional loops to $\Sigma^{(2)}_{a-t}$. Each of these solutions which verifies in addition 
the cube equations gives rise to a homomorphism from $H_1(M)$ to some polynomial ring.

We have found four such solutions in Jones case and one in Alexanders case. In Jones case there is a homomorphism

$S : H_1(M) \rightarrow  \mathbb{Z}[A,A^{-1},h,t,t^{-1}]$.

The new variable $t$ in this polynomial is rather surprising (see Subsection 3.1).

There are two homomorphisms

$S^+$ and $S^- : H_1(M) \rightarrow \mathbb{Z}[A,A^{-1},h,t,t^{-1}]$.

They have the property that $S^-([rot]) =0\quad  (mod 2)$ for each positive closed braid and that $S^+([rot]) =0\quad  (mod 2)$ for
each negative closed braid.

There is a homomorphism

$X : H_1(M) \rightarrow \mathbb{Z}/2\mathbb{Z}[A,A^{-1},h,r,s]$.

Again, the new variables $r$ and $s$ are a surprise (see Subsection 3.3).

The homomorphism $X$ can be lifted to a homomorphism with integer coefficients in the case of closed braids.

In Alexanders case we could construct only one homomorphism. 
 There is a homomorphism

$\Phi : H_1(M) \rightarrow \mathbb{Z}[A,A^{-1}]$.

$\Phi$ can detect chirality of a knot and it might be that it can detect non invertibility as well. Notice, that the Alexander
polynomial can not detect chirality and that no quantum knot invariant can detect non invertibility.

All five homomorphisms are trivial on the sliding classes.

Let K be a knot type with given framing $w$ and with Whitney index $n$ in the annulus. Hence, $(K,w,n)$
corresponds to a component of $M$. Let $[rot]$ be the corresponding canonical class of this component. For any
homomorphism $\rho$ from $H_1(M)$  into some polynomial
ring and such that $\rho$ is trivial on sliding classes we denote by $\rho_{K,w,n}$ the Laurent polynomial 
$\rho [rot]$. As we have explained, the $\rho_{K,w,n}$ are invariants of framed knots. But of course, the whole homomorphism 
$\rho$ applied to the corresponding first homology group is an invariant of framed knots too.

We finish this long introduction with some remarks.
\begin{remark}
The main question in knot theory can be formulated in the following way: are two given (framed) knot diagrams in the same
 component of $M$? Quantum knot 
invariants allow sometimes to answer this question (in the negative).
Our five homomorphisms $\rho$ are certainly new invariants because they allow sometimes to answer (in the negative)
a new question: are two given loops in $M$ homologic? (Of course, if they are 
homologic then they are in the same component of $M$. The inverse is true for canonical loops modulo sliding loops.) Notice, that the first Betti number of single 
components of 
$M$ can be arbitrary 
large. (Each boundary incompressible torus in the complement of $K$ in $V$ bounds a solid torus $V'$ in $V$. 
Consequently, we can perform our two rotations for $V'$ and we obtain in general new loops in $M$.)

In the case of classical knots all our invariants are polynomials of just the variable $A$. However, it is known that 
 cabling leads to new knot invariants already for the Jones polynomial . The same is true for our invariants, which become
polynomials with unexpected variables.
\end{remark}

\begin{remark}
The main shortcoming of our work is the lack of a computer program in order to calculate the invariants $\rho_{K,w,n}$
in examples. The calculation is much harder than the calculation of quantum knot invariants. We have calculated by
hand  very view (and very simple) examples, just enough in order to show that each of our invariants is not always trivial.
\end{remark}

\begin{remark}
We do not know what are the relations of our invariants with the known knot polynomials. But notice, that we use at four different places the fact that the diagrams $D$ represent knots and not links
 (compare Subsection 3.5.). It seems to be impossible to generalize
our results to links without passing to proper quotients of the polynomial rings. Notice also, that our invariants are specifique 
invariants for knots in the standard solid torus in 3-space. Indeed, the complement of the meridian of a knot in a 3-manifold other than $S^3$ is not
a solid torus. Moreover,  a  knot in a closed 3-manifold, which is not null homotopic and which is not a satellite, does not allow any non trivial isotopies to itself. 
\end{remark}

\begin{remark}
It is essential in our approach to exclude Reidemeister I moves. Otherwise $\Sigma^{(2)}$ would contain strata 
which correspond to a cusp with a transverse branch in $pr$. A homotopy through such a stratum would create or eleminate 
a single triple crossing and our method would break down completely. Moreover, it seems that the invariants
$\rho_{K,w,n}$ depend indeed rather non trivial on the framing $w$, in difference to quantum knot invariants.
But examples have to be calculated (see Remark 2).

It is interesting to notice analogies of our approach with the classical approach. In the latter one needs invariance under
passing $\Sigma^{(1)}_t$ and $\Sigma^{(1)}_a$. Passing a diagram with a cusp is not really important. In our approach 
$\Sigma^{(1)}_t$ is replaced by $\Sigma^{(2)}_{q}$ and $\Sigma^{(1)}_a$ is replaced by $\Sigma^{(2)}_{a-t}$ . Cusps are replaced by sliding loops, which are 
again not really important.
\end{remark}

\begin{remark}
The space of knots in the 3-sphere, in particular its homology, was studied in \cite{H2}, \cite{BC}, \cite{Bu}. Cocycles with integer or 
$ \mathbb{Z}/ 2 \mathbb{Z}$ coefficients for the space of knots in $ \mathbb{R}^3$, in particular for long knots, were constructed in
\cite{T} and \cite{V}. Let $F$ be a planar surface and let $M_F$ be the space of all regular knots in $F \times I$. Integer 
valued one and two dimensional coboundaries in $M_F$ were  
constructed in \cite{F1}.
\end{remark}

\begin{remark}
Our first attempt was to use one parameter knot theory in order to refine finite type knot invariants. This did not work
for classical knots (because of the {\em extrem pair move} for the {\em trace graph}, see \cite{FK} and \cite{F2}). Fortunately, this move is 
irrelevant in our construction of knot polynomials (because we do not use the trace graph at all). Notice, that nevertheless our first
attempt has worked perfectly for closed braids and 
almost closed braids. It has led to new invariants which are no longer of finite type but which are still calculable in 
polynomial time (with respect to the number 
of  crossings)
and which can detect non invertibility of closed braids (see \cite{F2}).
\end{remark}

\section{Basic notions of one parameter knot theory and extensions of state models}

In the first two subsections we recall briefly the basic notions of one parameter knot theory. All details with
complete proofs can be found in \cite{FK} (see also \cite{F2}). In the remaining subsections we study triple points,
 autotangencies and we extend Kauffman's state models.

\subsection{The space of regular knots and the higher order Reidemeister theorem}

We fix once for all a coordinate system in $\mathbb{R}^3$ : $(\phi , \rho , z )$.
Here, $(\phi ,\rho ) \in S^1\times \mathbb{R}^+$ are polar coordinates of the plane  $\mathbb{R}^2 = \{ z = 0 \}$.
Let $pr : \mathbb{R}^3 \setminus z-axes \to \mathbb{R}^2 \setminus 0$ be the canonical projection 
$(\phi ,\rho , z) \to (\phi ,\rho )$.
A regular knot diagram $D$
 is an oriented knot in the solid torus $V = \mathbb{R}^3 \setminus  z-axes$, such that the restriction
$pr : D \to \mathbb{R}^2 \setminus 0$ is an immersion. In all our considerations we will not distinguish regular
knot diagrams which differ only by an isotopy which changes only the $z-coordinate$. However, the formal definition 
of the space of diagrams as a quotient would be rather complicated and it is easier to work in the space of knots.
Let us mention two special cases: {\em closed braids} are links $L$ in $\mathbb{R}^3 \setminus  z-axes$ such that
the restriction of $\phi$ to $L$ is non singular. {\em Almost closed braids} are knots $K$ in $\mathbb{R}^3 \setminus  z-axes$
such that the restriction of $\phi$ to $K$ has exactly two critical points (compare \cite{F2}).

 Let $M$ be the infinit dimensional space of all regular knot diagrams
in V. It has natural finite dimensional approximations by spaces of polynomial knots of given degree.

 A generic regular knot $D$ has only ordinary double points 
as singularities of $pr(D)$. Let $\Sigma$ be the discriminant in $M$ which consists of all non-generic 
diagrams of regular knots. 

The discriminant $\Sigma$ has a natural stratification:  $\Sigma = \Sigma^{(1)} \cup \Sigma^{(2)} \cup ...$,
where $\Sigma^{(i)}$ are the union of all strata of codimension i in $M$.

\begin{proposition}
(Reidemeisters theorem for regular knots)

$\Sigma^{(1)} = \Sigma_t \cup \Sigma_a$,

where $\Sigma_t$ is the union of all strata which correspond to diagrams with exactly one ordinary triple point 
(besides ordinary double points)  and $\Sigma_a$ is the union of all strata which correspond to diagrams with exactly one 
ordinary autotangency.
\end{proposition}

In the sequel we need also the description of $\Sigma^{(2)}$.

\begin{proposition}
$\Sigma^{(2)} = \Sigma_q \cup \Sigma_{a-t} \cup \Sigma_f \cup \Sigma_{int}$ 

where $\Sigma_q$ is the union of all strata which correspond to diagrams with exactly one  ordinary quadruple point,
$\Sigma_{a-t}$ is the union of all strata which correspond to diagrams with exactly one ordinary autotangency 
through which passes another branch transversally,
$\Sigma_f$ corresponds to the union of all strata of diagrams with an autotangency in an ordinary flex,
$\Sigma_{int}$ is the union of all transverse intersections of strata from $\Sigma^{(1)}$.
\end{proposition}

Let $\gamma$ and $\gamma'$ be two generic collections of oriented loops in $M$ and which are homologic. Then there is an 
oriented embedded surface $F$ in $M$ and a generic Morse function $f : F \to [0, 1]$ such that $\gamma$ and 
$\gamma'$ are the oriented boundary of $F$ and $f^{-1}(0) = \gamma$ and $f^{-1}(1) = \gamma'$. Let $f^{-1}(t), t \in [0, 1]$
denote the fibres of $f$ in $F$.
\begin{theorem}
(higher order Reidemeister theorem)

 $\gamma$ and $\gamma'$ are homologic in $M$ if and only if there is a surface $F$ and a Morse function $f$ as described above
and such that for generic values of $t$ the fiber $f^{-1}(t)$ intersects $\Sigma$ transversally (hence, only in $\Sigma^{(1)}$) 
and for finitely many values of $t$ exactly one of the following {\em moves} happens:

(i) the non singular fiber $f^{-1}(t)$ has an ordinary tangency with $\Sigma^{(1)}$

(ii) the non singular fiber $f^{-1}(t)$ intersects $\Sigma^{(2)}$ transversally (as a fiber in a one parameter family of non 
singular fibers)

(iii) the fiber $f^{-1}(t)$ is singular and it intersects $\Sigma$ transversally (hence, the collection of loops $f^{-1}(t)$
changes just by a Morse modification).

\end{theorem}

The theorem with homology replaced by homotopy follows immediately as a special case from Theorem 1.10.
in \cite{FK}. The statement about homology follows then immediately from the usual general position arguments in the (infinite dimensional)
space $M$.

\subsection{The canonical loop and the sliding loops}

We identify $\mathbb{R}^3 \setminus z-axes$ with the standard solid torus $V = S^1 \times D^2 \hookrightarrow \mathbb{R}^3 \setminus z-axes$.
We identify the core of V with the unit circle in $\mathbb{R}^2$.

Let $rot(V)$ denote the $S^1$-parameter family of diffeomorphismes of V which is defined in the following way:
we rotate the solid torus monotoneously and with constant speed around its core by the angle t , $t \in [0 ,2\pi]$, i.e.
all discs $( \phi = const ) \times D^2$ stay invariant and are rotated simoultaneously around their centre.

Let us introduce as preliminary object in this subsection the space $M'$ of {\em all} (not necessarily regular) knots in $V$.
Let $D$ be a generic knot in $V$.

\begin{definition}
The {\em pre-canonical loop \/}  $rot(D) \in M'$ is the oriented loop induced by $rot(V)$.
\end{definition}

Notice that the whole loop $rot(D)$ is completely determined by an arbitrary point in it.

The following lemma is an immediat corollary of the definition of the pre-canonical loop .
\begin{lemma}
Let $D_s , s\in [0 ,1 ]$ , be an isotopy of knots in the solid torus. Then $rot(D_s ) , s \in [0 ,1 ]$,
is a homotopy of loops in $M'$.
\end{lemma}

Let $D$ and $D'$ be two knot diagrams in the annulus. Let $w(D)$ be the {\em writhe} (see e.g. \cite{BZ}) and let
$n(D)$ be the Whitney index in the annulus (see the Introduction). The following lemma is an immediate consequence
of Lemma 1.1. in \cite{F1}:

\begin{lemma}
$D$ and $D'$ are in the same component of $M$ if and only if they represent the same knot type in $V$ and $w(D) = w(D')$
and $n(D) = n(D')$.
\end{lemma}

If in a knot isotopy from $D$ to $D'$ we would have to performe a Reidemeister I move which eleminates a crossing 
(together with a
small {\em curl}) then we keep simply the crossing. The four types of curls are shown in Fig. 4. If we have to create a
\begin{figure}
\centering 
\psfig{file=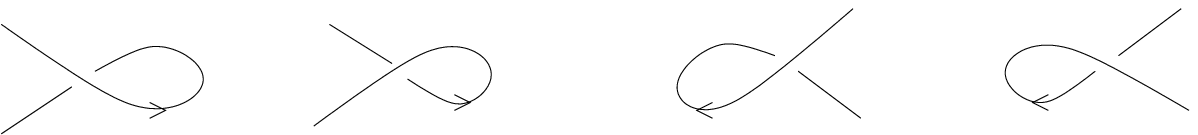}
\caption{}
\end{figure}
new crossing then we do this with the 
{\em Whitney trick} instead of the Reidemeister I move. This is illustrated in Fig. 5. In this way we replace the isotopy
\begin{figure}
\centering 
\psfig{file=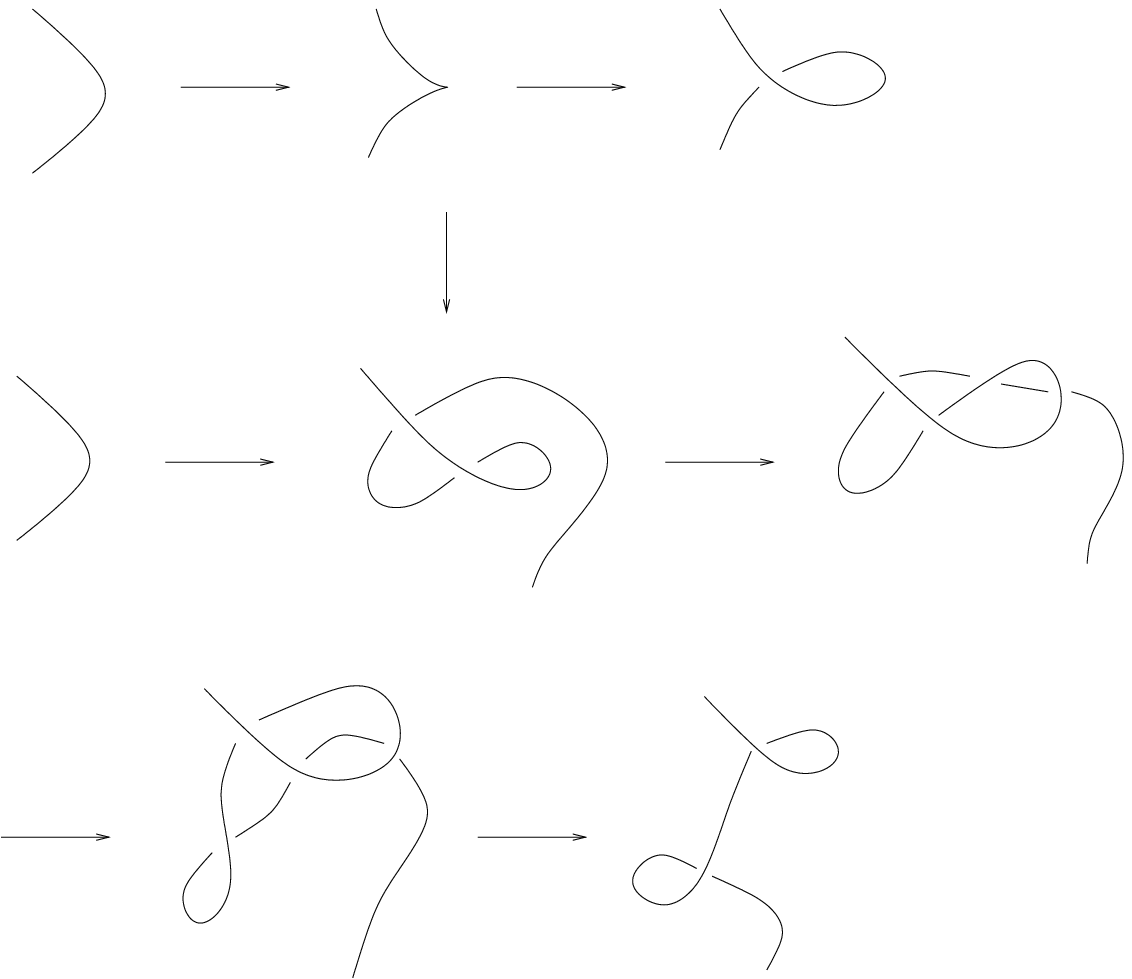}
\caption{}
\end{figure}
from $D$
to $D'$ by a regular isotopy to some $D''$. The latter differs from $D'$ by small curls at different places. We have to 
push these curls 
along $D''$ to collect them all in a small disc and then to eleminate them pairwise by using the Whitney trick in the 
opposite direction. In this way we approximate the given isotopy from $D$ to $D'$. The only ambiguity in this approximation is the
direction in which we push the curls along $D''$.

\begin{definition}
A {\em sliding loop} is a loop in $M$ which can be represented by pushing a small curl once along the whole diagram $D$.
\end{definition}

Let $\gamma$ be a loop in $M'$ and let $D$ be a generic point in it.
The loop $\gamma$ can be seen as an isotopy from $D$ to itself. We approximate $\gamma$ with the above construction 
by a loop in the component of $M$ which contains $D$. This approximation is unique for a given $D$ up to adding sliding 
loops. 

\begin{definition}
The above constructed approximation in $M$ of the pre-canonical loop $rot(D)$ in $M'$ is called the {\em canonical loop} and is
denoted by $rot(D)$ too. Its homology class, which is well defined  up to adding homology classes of sliding loops, is
called the {\em canonical class} of the corresponding component of $M$.
\end{definition}

This definition is correct because it follows from Lemma 1 that if we choose another $D$ in the same component of $M$ 
then the canonical class $[rot(D)]$ modulo sliding classes does not change.

\begin{remark}
The trace graph of an isotopy lives in a thickened torus $T^2$ (see \cite{FK}). One easily sees that the trace graph of
a sliding loop is homological non trivial in $H_1(T^2)$ provided that $D$ is homological non trivial in $V$. Consequently,
sliding loops represent in general non trivial elements in $H_1(M)$. We will not carry this out, because we do not need it
in the sequel. In the contrary, we will show that our 1-cohomolgy classes are trivial on all homology classes which can be
represented by sliding loops.
\end{remark}

In practice, a canonical loop is given by a sequence of pictures of diagrams corresponding to Reidemeister II or Reidemeister III
moves. It is rather tedious to draw these pictures. However, in the special case of closed braids in $V$ the canonical loop has a simpler 
combinatorial description. A closed n-braid $\hat \beta$ is a knot in the solid torus $V = \mathbb{R}^3 \setminus  z-axes$, such that 
$\phi : \hat \beta \to  S^1$ is non-singular and $[\hat \beta] = n \in H_1(V)$.

We define the space of closed braids, denoted by $M^{braid}$ exatly as $M$ but with the additional condition
that it consists only of closed braids.

Let $\Delta \in B_n$ be Garside's element, i.e. $\Delta^2$ is a generator of the centre of $B_n$ (see \cite{B}).
Geometrically, $\Delta^2$ is the full twist of the n strings.
\begin{definition}
Let $\gamma \in B_n$ be a braid with closure isotopic to $\hat \beta$. Then the {\em combinatorial canonical loop \/}
$rot(\gamma)$ is defined by the following sequence of braids:

$\gamma \to \Delta\Delta^{-1}\gamma \to \Delta^{-1}\gamma\Delta \to \dots \to \Delta^{-1}\Delta\gamma' \to \gamma' 
\to \Delta\Delta^{-1}\gamma' \to \Delta^{-1}\gamma'\Delta \to \dots \to \Delta^{-1}\Delta\gamma \to 
 \gamma$. 

Here, the first arrow consists only of Reidemeister II moves, the second arrow is a cyclic permutation of the braid word 
(which corresponds to an isotopy of the braid diagram in the solid torus) and the following arrows consist of "pushing $\Delta$
monotoneously from the right to the left through the braid $\gamma$". We obtain a braid  $\gamma'$ and we start again.
\end{definition}We give below a precise definition in the case $n = 3$. The general case is a straightforward generalization which is left to the reader.
$\Delta = \sigma_1\sigma_2\sigma_1$ for $n = 3$. We have just to consider the following four cases:$\sigma_1\Delta = \sigma_1(\sigma_1\sigma_2\sigma_1) \to \sigma_1(\sigma_2\sigma_1\sigma_2) = \Delta\sigma_2$$\sigma_2\Delta = (\sigma_2\sigma_1\sigma_2)\sigma_1 \to (\sigma_1\sigma_2\sigma_1)\sigma_1 = \Delta\sigma_1$

$\sigma_1^{-1}\Delta = \sigma_1^{-1}(\sigma_1\sigma_2\sigma_1) \to \sigma_2\sigma_1 \to (\sigma_1\sigma_1^{-1})\sigma_2\sigma_1 \to \sigma_1(\sigma_2\sigma_1\sigma_2^{-1}) = \Delta\sigma_2^{-1}$

$\sigma_2^{-1}\Delta = \sigma_2^{-1}(\sigma_1\sigma_2\sigma_1) \to (\sigma_1\sigma_2\sigma_1^{-1})\sigma_1 \to \sigma_1\sigma_2 \to \sigma_1\sigma_2\sigma_1\sigma_1^{-1} = \Delta\sigma_1^{-1}$.

Notice, that the sequence is canonical in the case of a generator and almost canonical in the case of an inverse generator. Indeed,
we could replace the above sequence   $\sigma_1^{-1}\Delta \to \Delta\sigma_2^{-1}$ by

$\sigma_1^{-1}(\sigma_1\sigma_2\sigma_1) \to \sigma_2\sigma_1 \to \sigma_2\sigma_1\sigma_2\sigma_2^{-1} \to (\sigma_1\sigma_2\sigma_1)\sigma_2^{-1}$.

But it turns out that the corresponding canonical loops in $M(\hat \beta)$ differ just by a homotopy which passes once transversally through 
a stratum of $\Sigma_{a-t}^{(2)}$.

Let c be the word lenght of $\gamma$. Then we use exactly $2c(n-2)$ braid relations (or Reidemeister III moves) in the combinatorial canonical loop.
This means that the corresponding loop in $M$ cuts $\Sigma^{(1)}_t$ transversally in exactly $2c(n-2)$ points. Notice, that
we could also push directly $\Delta^2$ once through the braid.

One easily sees that the combinatorial canonical loop $rot(\gamma)$ from Definition 4 is homologic in $M$ 
to the canonical loop $rot(\hat \beta)$ from Definition 3 modulo loops which are generated by rotations of the solid torus $V$
along its core. But as already mentioned, all our homomorphisms are trivial on such loops (because they consist only
of isotopies of {\em diagrams}). Moreover, there are no sliding loops here because we can assume that the 
isotopy which connects two closed braids is already regular, see e.g. \cite{M}, and in the rotation of a closed braid there appear 
no cusps in the projection.
\subsection{The classification of triple points and of autotangencies. The distinguished crossing of a triple point. Homological markings}

We number the eight types of triple points as shown in Fig. 6. Notice, that there are exactly six types of
\begin{figure}
\centering
\psfig{file=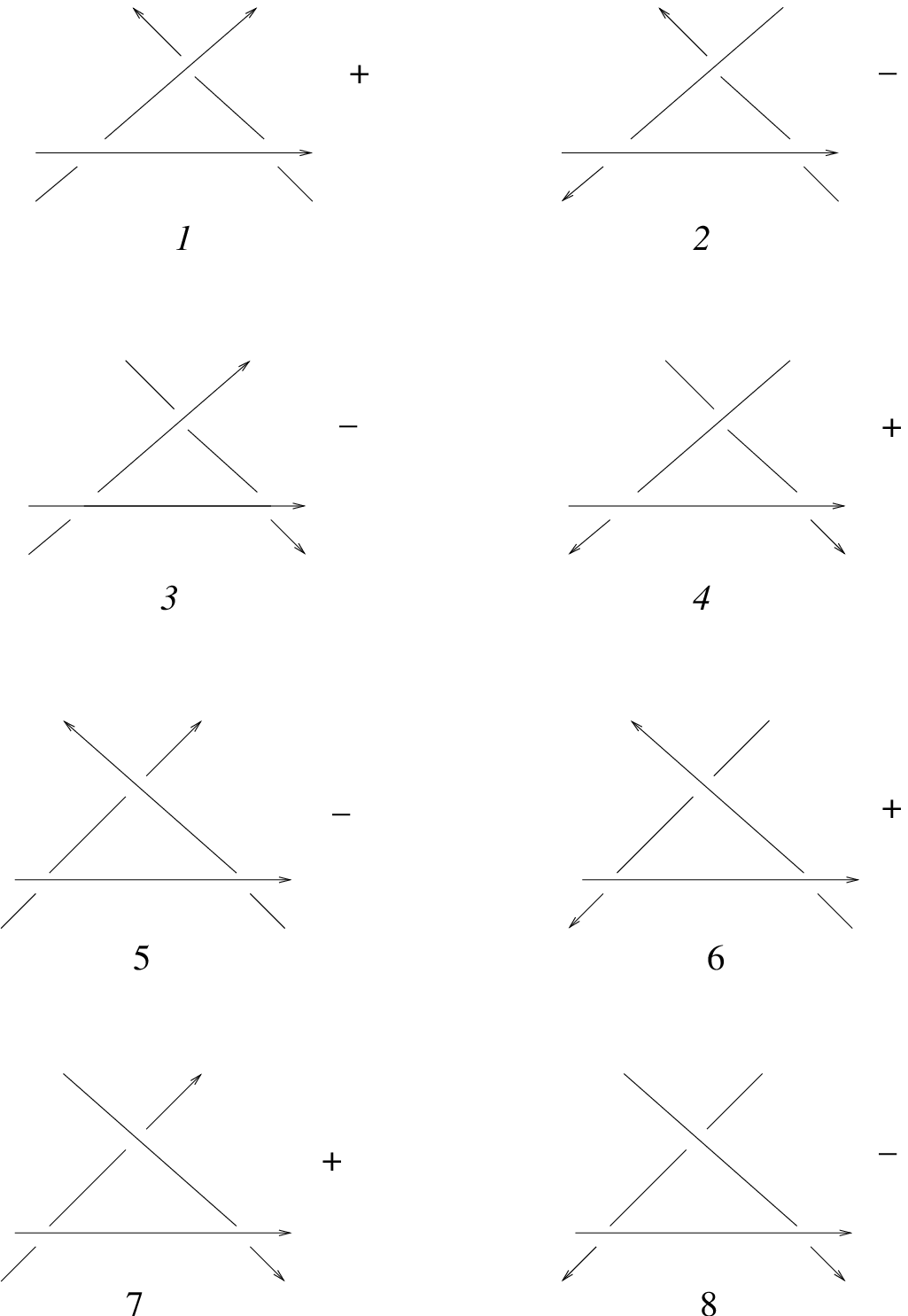}
\caption{}
\end{figure}
{\em braid-like} triple points and two types of {\em star-like} triple points, namely $2$ and $6$ (see \cite{AF}).
The pictures show the diagram near the triple crossing at one (local) side of the corresponding stratum of $\Sigma^{(1)}_t$.
The signs in Fig. 6 define the given sides as positive or negative. The positive coorientation is then the direction 
from the negative to the positive side. We will show in the next section that the closure in $M$ of $\Sigma^{(1)}_t$ with
this coorientation becomes an integer cycle of codimension one. 

Let $M_0$ be the union of all components of $M$ with vanishing Whitney index $n(D)$ in the annulus.
Let $rot_\pi$ be the rotation in the canonical loop by the angle $\pi$. $rot_\pi$ acts as an {\em involution} on 
$M_0$ and it interchanges the remaining components of $M$. Let $M_0/rot_\pi$ be the quotient.

Let $p$ be a triple point of a diagram $D$. Let $rot_{\pi}(p)$ be the corresponding triple point in the diagram $rot_{\pi}(D)$.

\begin{lemma}
The involution $rot_{\pi}$ changes the sign of a triple point, i.e. $sign(p) = - sign(rot_{\pi}(p))$ and it changes the types as follows:
1, 2, 6, 8 stay invariant, 5 is interchanged with 7 and 3 is interchanged with 4.
\end{lemma}

The proof is an easy case by case verification and it is left to the reader.

Let $p$ be a triple crossing. The three branches of the diagram at the triple crossing are ordered by the hight of the
$z-coordinate$.

\begin{definition}
The crossing of the highest branch with the lowest branch is called the {\em distinguished crossing} of the triple crossing
$p$ and it is denoted by $d(p)$ (compare Fig. 7).
\end{definition}
\begin{figure}
\centering 
\psfig{file=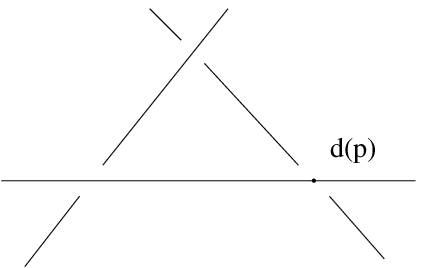}
\caption{}
\end{figure}
We identify $H_1(V)$ with $\mathbb{Z}$ by sending the core of V to the generator $+1$. If $D$ is a knot then we can attache to each crossing $q$ a 
{\em homological marking \/} $[q] \in H_1(V)$ in the following way: 
we smooth $q$ with respect to the orientation of $D$. The result is an oriented 2-component link. The component of this link which contains 
the undercross  which goes to the overcross at $q$ is called $q^+$. We associate now to $q$ the homology class $[q] = [q^+] \in H_1(V)$ (compare also \cite{F1}).
Notice, that the two crossings involved in a Reidemeister II move have the same homological marking and that a Reidemeister III move 
does not change the homological marking of any of the three involved crossings.

\begin{definition}
To a triple point $p$ we associate its {\em type} $j(p) \in \{1, ...,8\}$ (see Fig. 6).

To a distinguished crossing $d(p)$ of a triple point $p$ we associate its {\em homological marking} $[d(p)] = [d^+(p)] \in \mathbb{Z}$ 
as well as 
its {\em Whitney index in the annulus} $n(d(p)) = n(d^+(p)) \in \mathbb{Z}$.
\end{definition}

The following lemma is evident.
\begin{lemma}
Let $[D] = m \in H_1(V) = \mathbb{Z}$.
Then $d(rot_{\pi}(p)) = rot_{\pi}(d(p))$ and $[d(rot_{\pi}(p))] = m - [d(p)]$.
\end{lemma}

The four types of autotangencies together with their coorientation are shown in Fig. 8. We call the autotangencies $1$
\begin{figure}
\centering 
\psfig{file=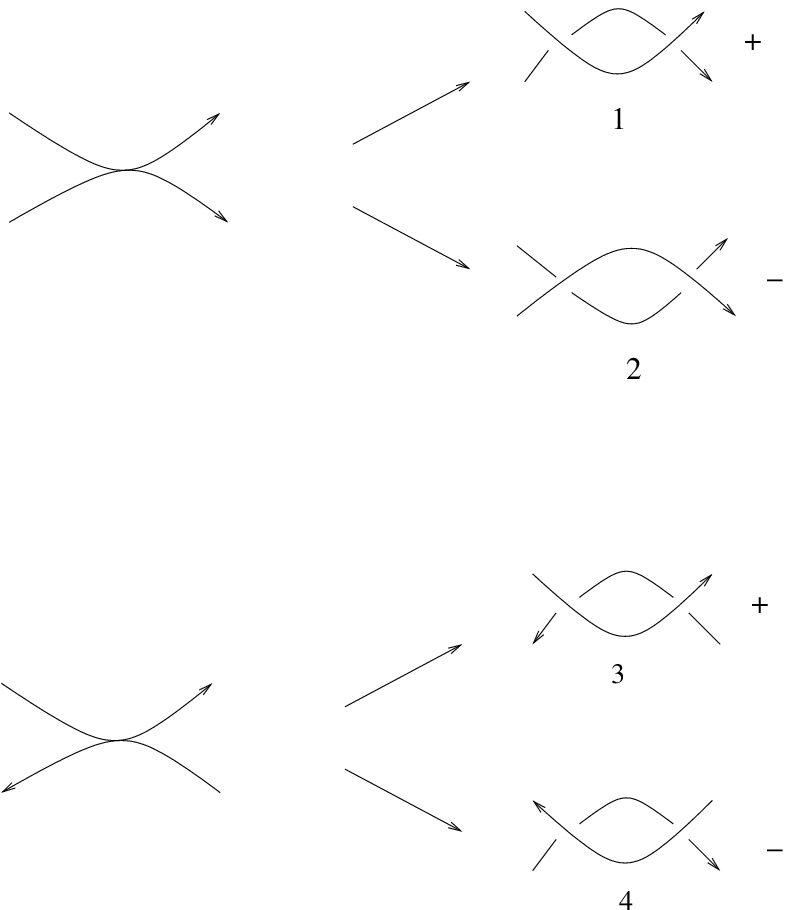}
\caption{}
\end{figure}
 and $2$ , respectively $3$ and $4$, the {\em dual} autotangencies. We will see in the next
subsection that the closure of $\Sigma^{(1)}_a$  with this coorientation is {\em not} an integer cycle of codimension 
one in $M$. However, it will turn out that this is the right coorientation in order to define our homomorphisms $S^+,S^-,X$
with integer coefficients instead of coefficients in $ \mathbb{Z}/ 2 \mathbb{Z}$.

\begin{definition}
To an autotangency $p$ we associate its {\em type} $j(p) \in \{1, ...,4\}$ as well as its {\em homological marking}
$[d(p)] = [d^+(p)] \in \mathbb{Z}$. Here $d(p)$ is any of the two crossings (they have the same homological marking). 

\end{definition}

The verification of the following easy lemma is left to the reader .

\begin{lemma}
The involution $rot_{\pi}$ preserves the types and the signs of the autotangencies 1 and 2. It interchanges the types
3 and 4 and changes their signs. 
In all four cases we have $[d(rot_{\pi}(p))] = m - [d(p)]$.

\end{lemma}
\subsection{The cube of triple points and the coorientations of triple points and of autotangencies}

There are exactly 24 types of autotangencies with a transverse branch. We show 12 of them in Fig. 9. The remaining 
12 types are obtained by those from Fig. 9 by either taking the {\em mirror image} (i.e. switching all crossings) or by
\begin{figure}
\centering
\psfig{file=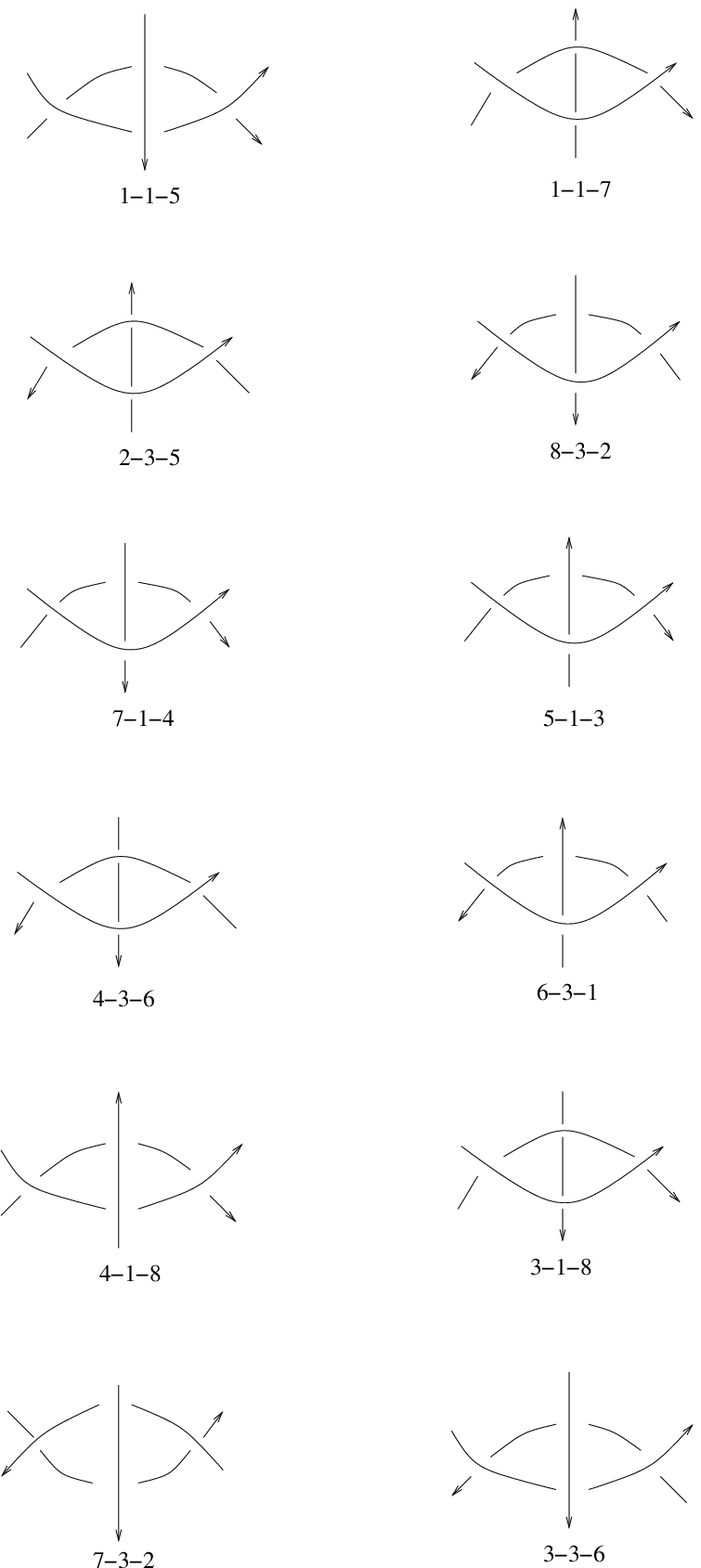}
\caption{}
\end{figure}
inversing the orientation of all three branches simultaneously.
Our notation convention is the following : e.g. 2-3-5 means that the triple crossing on the left is of type 2, the autotangency 
is of type 3, and the triple crossing on the right is of type 5.
 An example of the intersection with $\Sigma$ of a meridional disc of
$\Sigma_{a-t}^{(2)}$ is shown in Fig. 10 (compare  \cite{FK}). We encode the type of the stratum of $\Sigma_{a-t}^{(2)}$
\begin{figure}
\centering 
\psfig{file=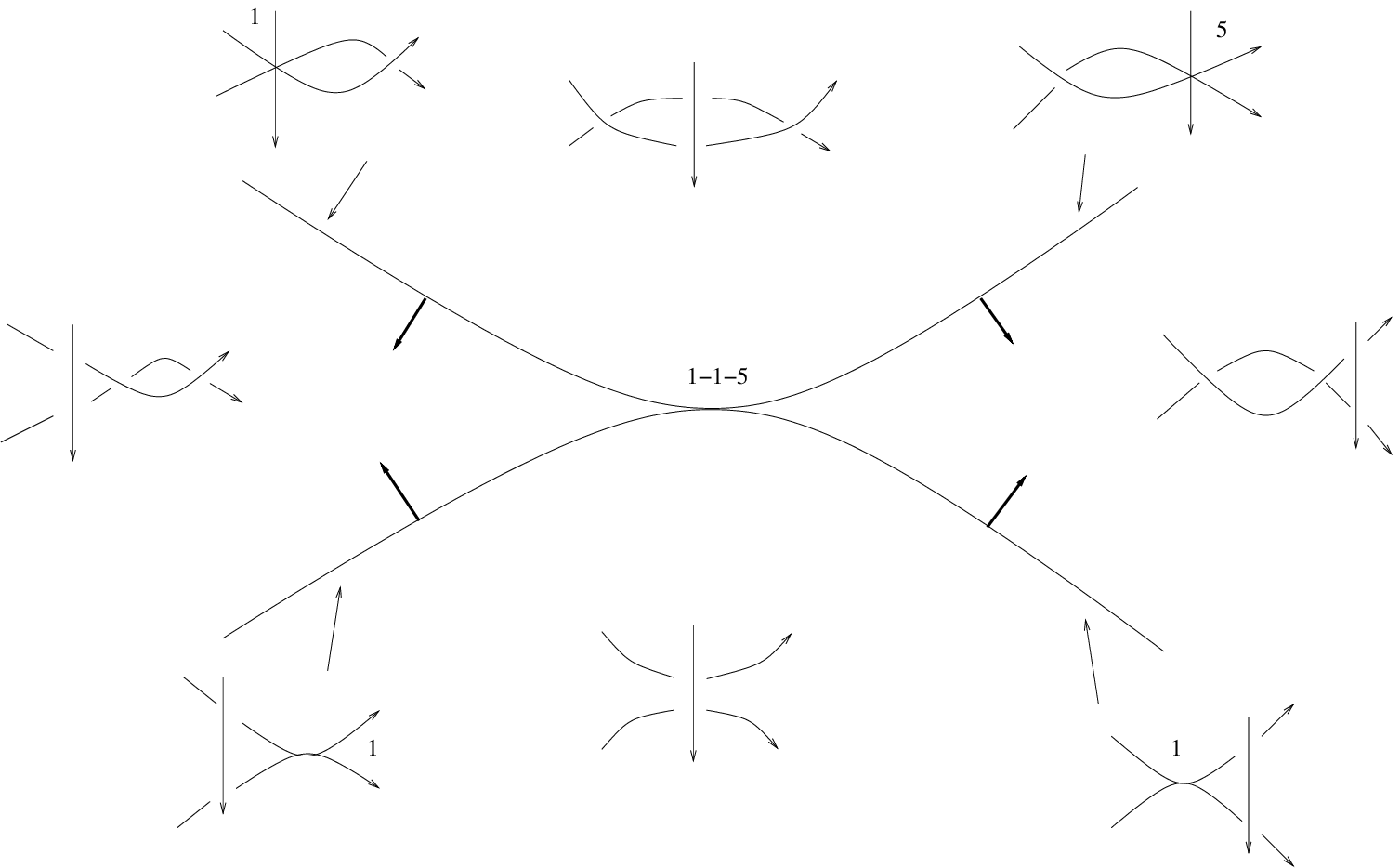}
\caption{}
\end{figure}
by the numbers of the two involved triple points and of the autotangency as shown in Fig. 9. The coorientations
of the strata are indicated by the positive normal vector. We construct now the graph $G$
as described in the Introduction.
\begin{lemma}
The graph $G$ is obtained from the 1-skeleton of a 3-dimensional cube by doubling each edge. 
\end{lemma}
{\em Proof:\/} The projection of a triple crossing to a small disc consists of six arcs which all meet at the triple point.
The six arcs divide the disc into six regions. Take such a region. The two arcs in its boundary have a writhe at the 
triple crossing. We close this region by adding a new intersection point of the two arcs , but with opposite writhe
(see Fig. 11). The region corresponds now to an autotangency by shrinking it to a point. Antipodal regions correspond
\begin{figure}
\centering
\psfig{file=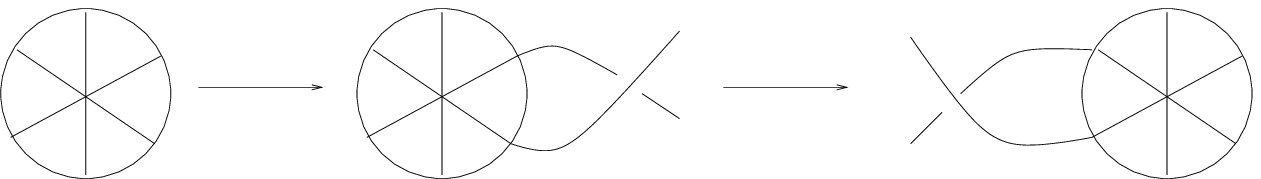}
\caption{}
\end{figure}
to dual autotangencies. The starting triple crossing together with the choosen region determines now another triple
crossing as shown in Fig. 11 too. The two triple crossings together with the autotangency corresponds to a stratum 
of $\Sigma_{a-t}^{(2)}$. Hence, each vertex of $G$ is adjacent to exactly six edges of $G$. A case by case consideration 
(which we left to the reader) shows that antipodal regions connect always the same couple of triple crossings.
Hence, each vertex of $G$ is connected with exactly three different vertices. Each of the edges can be replaced by another 
edge corresponding to the dual autotangency. We show $G$ in Fig. 12 (where all edges have to be considered as
\begin{figure}
\centering
\psfig{file=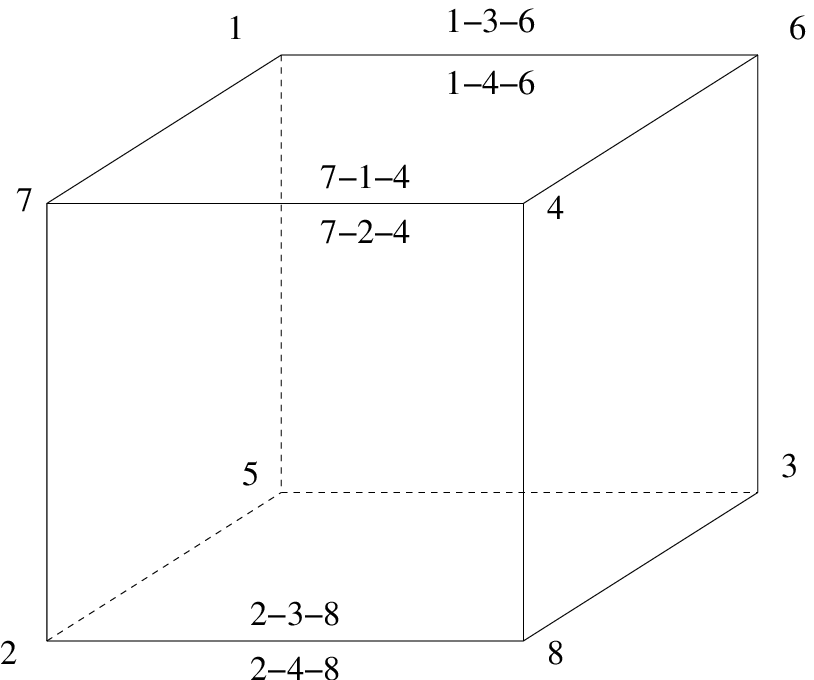}
\caption{}
\end{figure}
double edges). $\Box$

The closure of $\Sigma^{(1)}_t$ with the choosen coorientations becomes an integer cycle of codimension one in $M$
if (and only if) for each edge of $G$ the coorientations of the two vertices fit at $\Sigma_{a-t}^{(2)}$ as shown in Fig. 10. (The strata of $\Sigma^{(1)}_t$
intersect mutually transvers at  $\Sigma_{q}^{(2)}$ and hence, the coorientations do not change under passing a quadruple point.)
This is indeed the case. We left the verification to the reader.

There are evidently exactly four types of strata of $\Sigma_{f}^{(2)}$. At each type a couple of dual autotangencies 
come together as shown for example in Fig. 13 (compare \cite{FK}). We see immediately from Fig. 13 that the closure
\begin{figure}
\centering 
\psfig{file=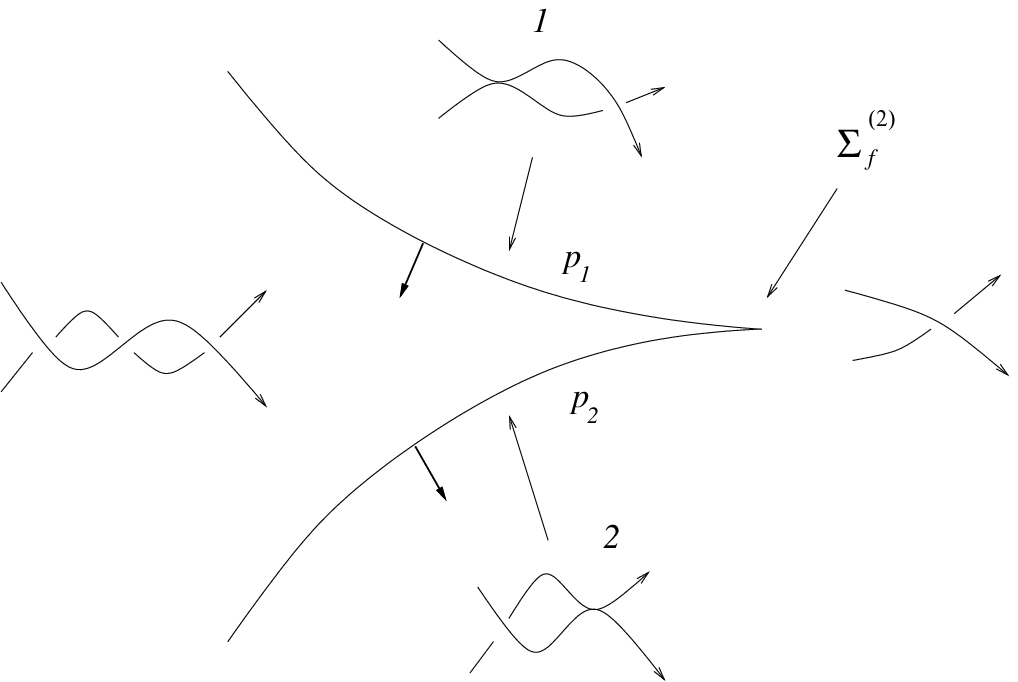}
\caption{}
\end{figure}
of $\Sigma^{(1)}_a$  in $M$ provided with our choosen coorientation is not an integer cycle.

\begin{definition}
Let $\gamma$ be an oriented arc which cuts $\Sigma^{(1)}$ transversally at a point $p$. Then the {\em intersection index}
$sign(p)$ is +1 if the orientation of $\gamma$ at $p$ coincides with the coorientation of $\Sigma^{(1)}$ at $p$ and
it is -1 otherwise.
\end{definition}

\begin{remark}
Our coorientation for strata of triple crossings does not coincide with the coorientation used previously in \cite{FK}. The
latter was defined globally but using only the underlying planar curve. The definition of the present coorientation is local,
but it uses the knot instead of the underlying planar curve. This fact is of crucial importance in order to find solutions of the 
tetrahedron equation (2) (compare Subsection 4.2.). With the previous coorientation there aren't any non trivial solutions.
\end{remark}
\subsection{The tetrahedrons of triple points}
One easily sees that there are exactly 48 different types of quadruple points (and consequently 48 different types of strata
of $\Sigma_{q}^{(2)}$).
We consider the (unique) positive quadruple point (see the Introduction). The intersection of $\Sigma$ with a meridional
disc for $\Sigma_{q}^{(2)}$ is shown in Fig. 1 (see also \cite{FK}). In Fig. 14 we show the eight diagrams with a triple
\begin{figure}
\centering
\psfig{file=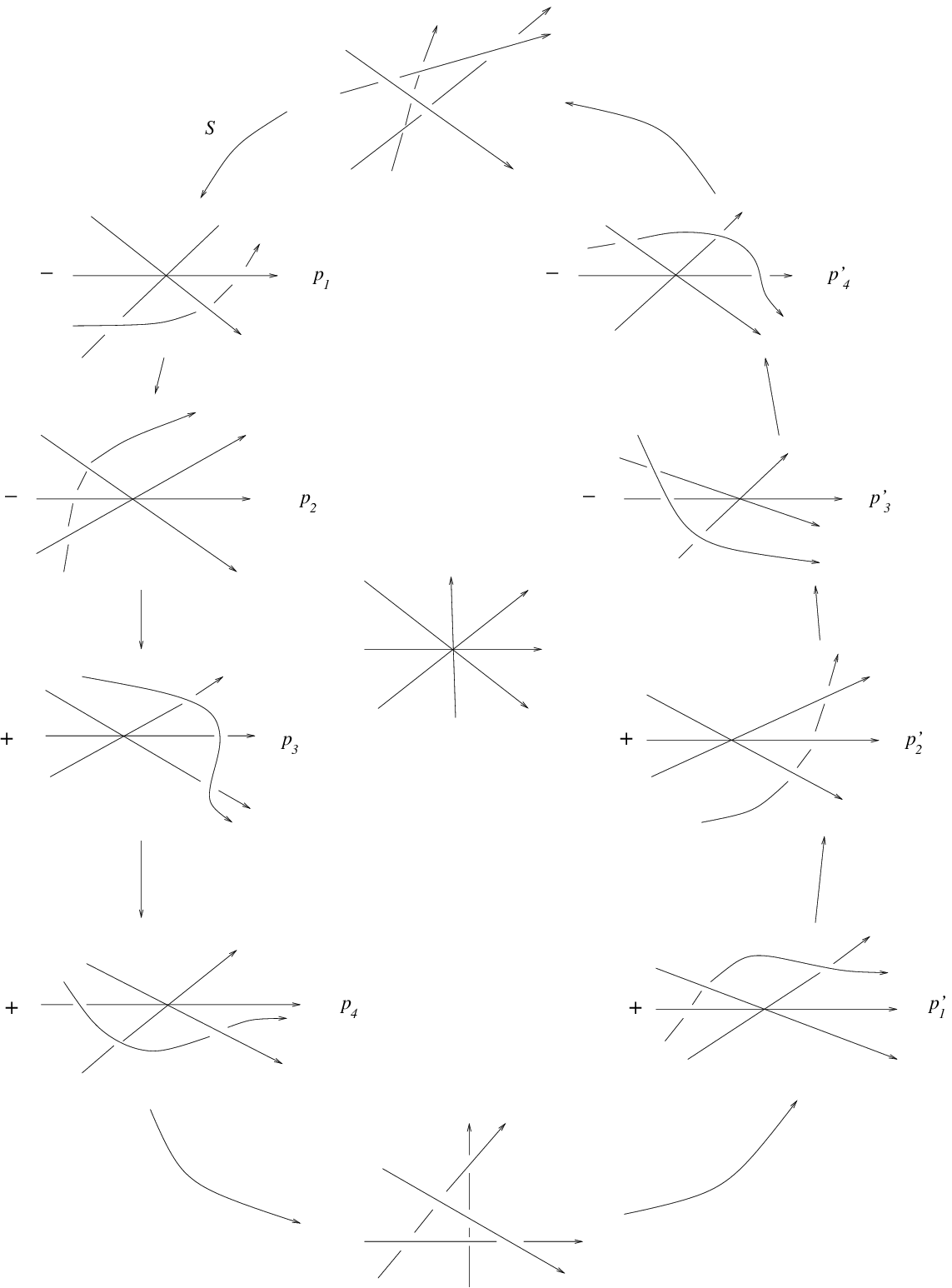}
\caption{}
\end{figure}
crossing and with their (easily established) signs in the boundary of the meridional disc. The triple crossings are all positive and we draw them just as
triple points. The eight diagrams are all identical outside of the drawn parts. Notice, that the meridional loop $s$ cuts twice each
of the four strata of triple points. The corresponding diagrams with triple crossings are denoted by $p_i$ and $p'_i$. Fig. 14
will be of great importance in the Subsections 4.2. and 4.3.

We indicate how to  associate tetrahedrons to $s$  (also we do not need this in 
the sequel): the four vertices correspond to four consecutive intersections of $s$ with $\Sigma^{(1)}_t$. Two vertices are connected
by an edge if and only if the two corresponding triple crossings have a commun crossing. The 2-faces of the tetrahedron correspond 
to the "vanishing triangles" in the nearby generic diagrams.

\begin{remark}
Exactly the triple crossings $p_1$ and $p_4$ (and hence also $p'_1$ and $p'_4$) share the same distinguished crossing.
This fact will lead to the new variable $t$ in the homomorphisms $S, S^+, S^-$. 
\end{remark}
\subsection{The simplifications of triple points and of autotangencies and the extension of the Kauffman bracket}
Let $p$ be a braid-like triple crossing. We define the {\em simplifications} of $p$ as shown in Fig. 15. Notice, that they
\begin{figure}
\centering 
\psfig{file=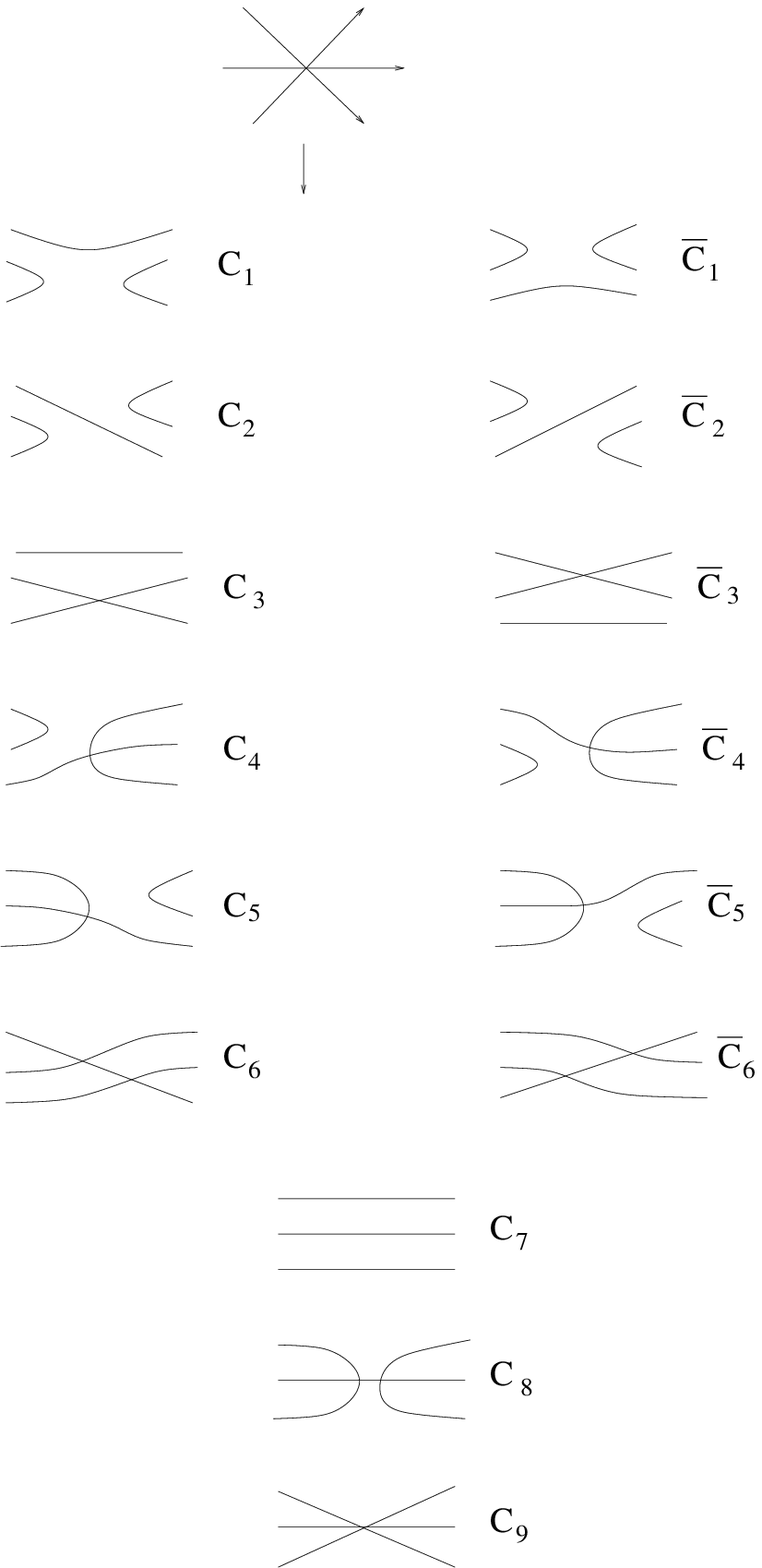}
\caption{}
\end{figure}
are determined by the planar triple point alone. 

\begin{remark}
The involution $rot_{\pi}$ acts on the simplifications. We denote by $C_i$ and $\bar C_i$ the couples of simplifications
which are interchanged by $rot_{\pi}$. Some simplifications are invariant under $rot_{\pi}$. But, according to Lemma 3 
the triple crossings $p$ and $rot_{\pi}(p)$ have different signs. It turns out, that their contributions cancel out in $L(rot)$
 ( we will see in 
Section 4  that triple crossings of  type 3 can be identified with those of type 4 and respectively, type 5 can be 
identified with type 7). Therefore, in the sequel we will not consider configurations which are invariant under the
involution $rot_{\pi}$.
\end{remark}

Our notation convention is the following: to the simplification of number $i$ of a triple point $p$ of type $j$ we 
associate the (independent) variable $C_i(j)$. In the case of positive triple points we write instead of $C_i(1)$ just
$C_i$.

Let $p$ be a star-like triple point (i.e. of type $2$ or $6$). Then the simplifications are given in Fig. 16.
\begin{figure}
\centering 
\psfig{file=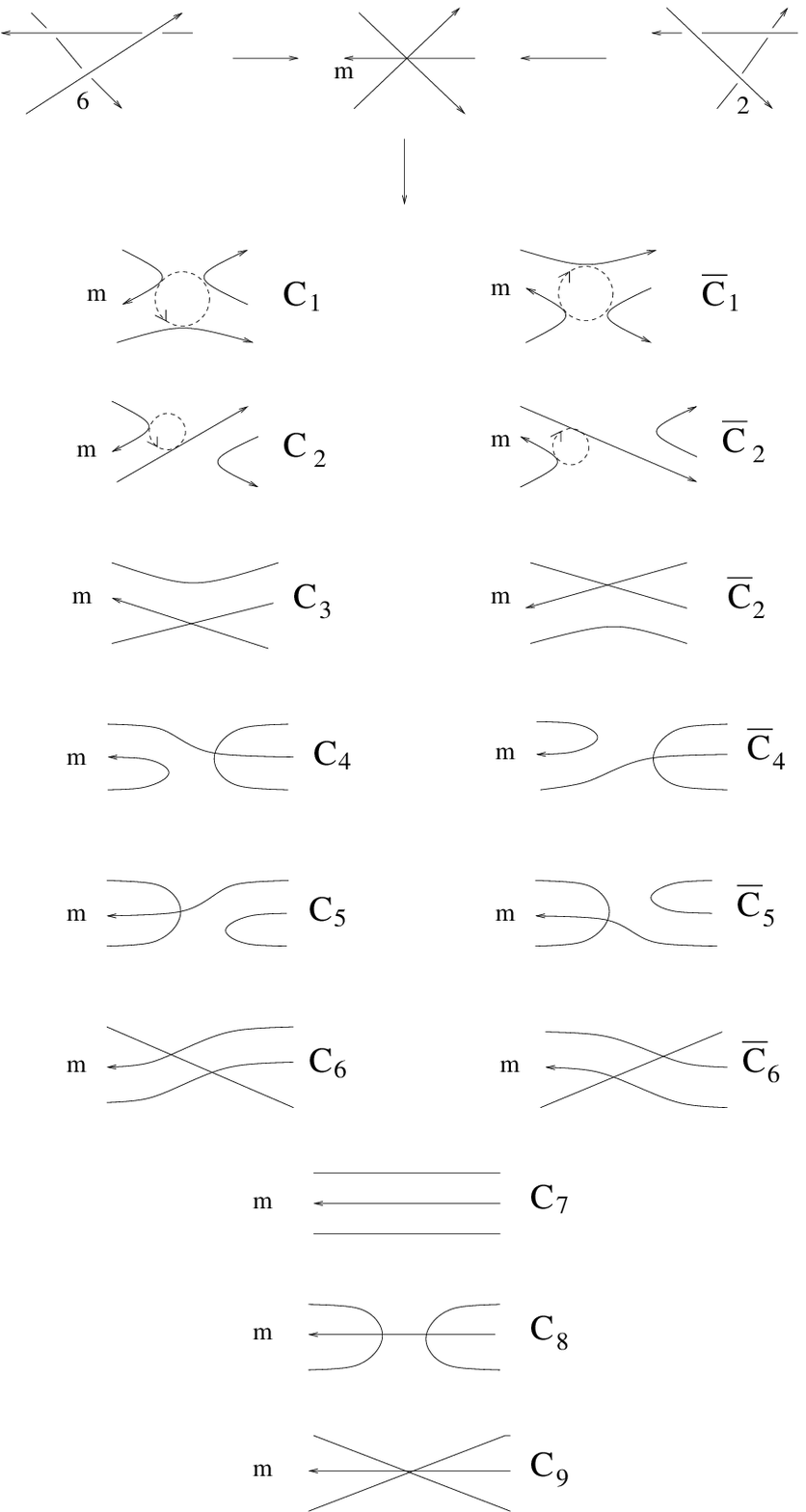}
\caption{}
\end{figure}
Notice, that for the definition of the simplifications we need  the diagrams (which allow us to distinguish the three branches)
instead of  the planar curves alone. We denote the branch with the z-coordinate in the middle by $m$. We define the
 simplifications
$C_i(2)$ and $C_i(6)$ by using the middle branch (equiped with a $m$) and by using the orientation of the small
dotted circle as shown in Fig.16. Again, the types of the corresponding simplifications of $rot_{\pi}(p)$ are denoted by $\bar C_i(2)$
and $\bar C_i(6)$.

There are four {\em coarse} types of simplifications: with  0, 1, 2, or 3 double points. We call the simplifications with 0
double points also the {\em smoothings}. As we have seen so far: {\em to each type of a triple point we can associate
six couples of independent variables.}

We define now the extension of Kauffman's bracket. Let $D_p(C_i)$ (or $D_p(\bar C_i)$be a diagram with triple crossing
 $p$ which
is simplified by $C_i$ (respectively $\bar C_i$). The result is an unoriented link diagram in the annulus.  If $D_p(C_i)$ has no (double point)
singularities, then the Kauffman bracket $<D_p(C_i)> \in \mathbb{Z}[A,A^{-1},h]$ is defined as usual (see  \cite{HP}).
Otherwise, each {\em Kauffman state} consists of embedded circles in the annulus besides one, two or three double 
points. We take now an {\em abstract resolution} of the double points (i.e. we separate the two branches without distinguishing 
the overcross from the undercross in the solid torus). It will turn out that the case of two or three double points can be ignored,
because these cases will not occure in non trivial solutions of the tetrahedron equation (compare Subsection 4. 2.).
There are exactly three possibilities in the case of one double point. They are shown in Fig. 17. The circle in a) counts
\begin{figure}
\centering 
\psfig{file=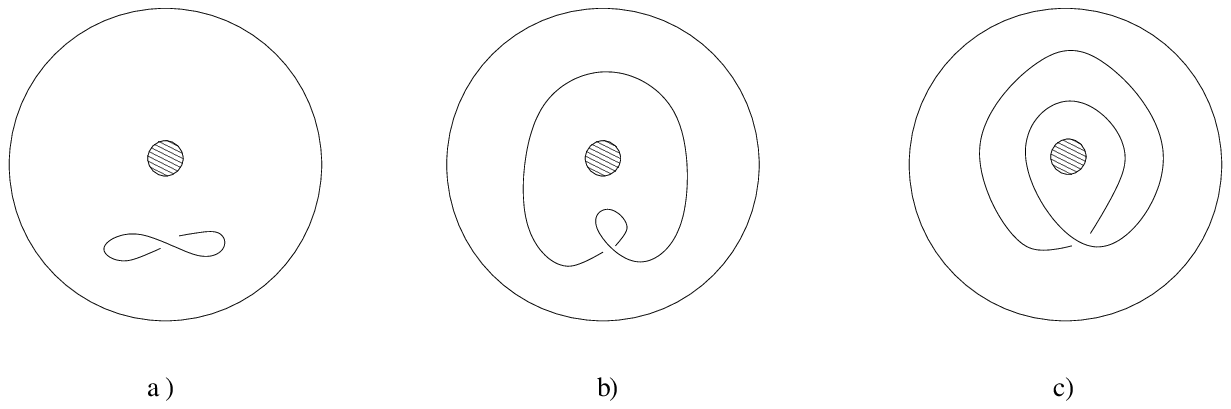}
\caption{}
\end{figure}
as a $d-circle$ and contributes consequently a factor  $-A^2 -A^{-2}$ in the monomial associated to the Kauffman state.
The circle in b) counts for $h$ as usual. Let us consider a Kauffman state which contains a  circle as shown in c), 
which we call a {\em double circle}. All other circles  in this Kauffman state are of type $d$ or $h$. The double circle 
separates the two boundary components of the annulus. Therefore we can consider the 
configuration of the $h-circles$ with respect to the double circle. Assume that there are exactly j $h-circles$ in the region 
adjacent to the inner boundary of the annulus and that there are exactly k  $h-circles$ in the region adjacent to the outer boundary
of the annulus. We illustrate the situation in Fig.18. Let \#(.) denote the number of (.).
\begin{figure}
\centering
\psfig{file=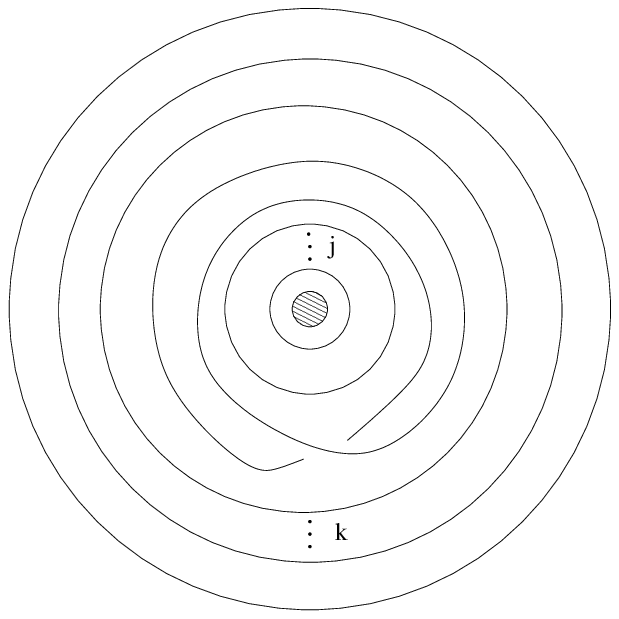}
\caption{}
\end{figure}
\begin{definition}
To the Kauffman state shown in Fig.18  we associate the monomial 

$A^{\#(A-smoothings) - \#(A^{-1}-smoothings)} (-A^2 -A^{-2})^{\#(d-circles)}
 h^{\#(h-circles)} r^j s^k$ 
 
Here, $r$ and $s$ are new independent variables.

We define now the Kauffman bracket $<D_p(C_i)> \in \mathbb{Z}[A,A^{-1},h,r,s]$ as the sum of the monomials over all Kauffmann states as usual.
\end{definition}

Let us consider the simplifications of autotangencies. There are exactly four types of autotangencies $p$ denoted by
$j(p) \in \{1, 2, 3, 4\}$. There are only two types of simplifications: $C_0$ and $C_x$. They are shown in Fig. 3. (The 
remaining simplification
corresponds just to a Reidemeister II move (which does not change the knot type) and is consequently  not interesting.)
As in the case of triple points, we denote the simplification of type $C_0$ or $C_x$ of an autotangency of type $j$ by
$C_0(j)$ respectively $C_x(j)$.

We define now the extension of the Kauffman bracket exactly as we have done for the simplifications of triple points.

Dual autotangencies can come together in strata of $\Sigma_{f}^{(2)}$. The meridional loops to $\Sigma_{f}^{(2)}$ give
relations for the variables associated to the simplifications. Let $s$ be a meridian for $\Sigma_{f}^{(2)}$. Let $p_1$ and
$p_2$ be the two autotangencies in $s$ (compare Fig.13).

\begin{lemma}
$L(s) = sign(p_1) \sum_{i \in \{0, x\}} <D_{p_1}(C_i)> + sign(p_2) \sum_{i \in \{0, x\}} <D_{p_2}(C_i)> = 0$

if and only if   $C_0(1) = - C_0(2)$, $C_x(1) = - C_x(2)$, $C_0(3) = - C_0(4)$, $C_x(3) = - C_x(4)$.
\end{lemma}
(In other words: the variables for the same simplification of an autotangency and its dual differ by the sign.)

{\em Proof:\/} The extension of the Kauffman bracket for autotangencies does not depend on the orientation of the diagram 
(in difference to its extension for triple crossings). It is therefore sufficient to prove that $L(s) = 0$ in the case shown in Fig. 13
and in its mirror image. The two autotangencies have the same sign. Let $T_2$ be the free $\mathbb{Z}[A,A^{-1}]$-modul generated by all 
chord diagrams 
with exactly two chords in the disc and which have no more than one double point. Each of the two autotangencies 
determines
 an element in $T_2$ as shown in Fig. 19. Evidently, $L(s) = 0$ independently of the rest of the diagram outside the
\begin{figure}
\centering 
\psfig{file=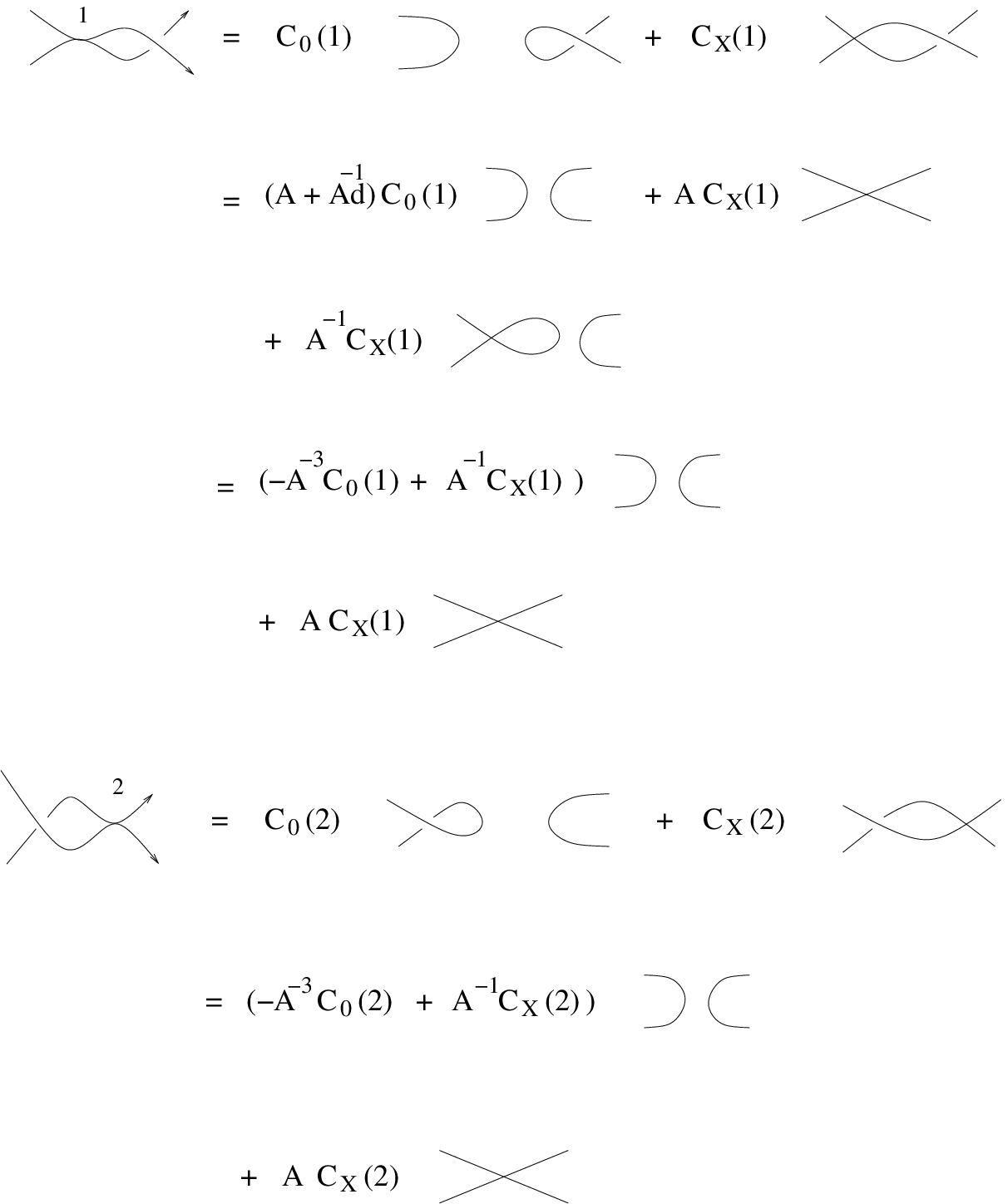}
\caption{}
\end{figure}
chord diagram 
if and only if for each generator of $T_2$ the signed sum of the coefficients is zero.  This gives us the equations:

$- A^{-3} C_0(1) - A^{-3} C_0(2) + A^{-1} C_x(1) + A^{-1} C_x(2) = 0$

$A C_x(1) + A C_x(2) = 0$.

The only solution is $C_x(1) = - C_x(2)$ and $C_0(1) = - C_0(2)$. The case of the mirror image is completely analogous and leads to
the same solution. $\Box$

\subsection{Kauffman's state model for the Alexander polynomial and markers for triple points and autotangencies}
The beautiful state model for the Alexander polynomial was introduced by Kauffman in \cite{K1}.
We modify it by putting the stars into the adjacent regions of the boundary of the annulus (compare the Introduction).
At a triple crossing a triangle (which is not a *-region) has shrinked to a point and there have to be markings in exactly 
two of the six regions of a small disc around  the triple point. Let $p$ be a braid-like triple point. We define the 
markings of $p$ in Fig. 20. Again, the markings come in couples, interchanged by the involution $rot_{\pi}$.
\begin{figure}
\centering 
\psfig{file=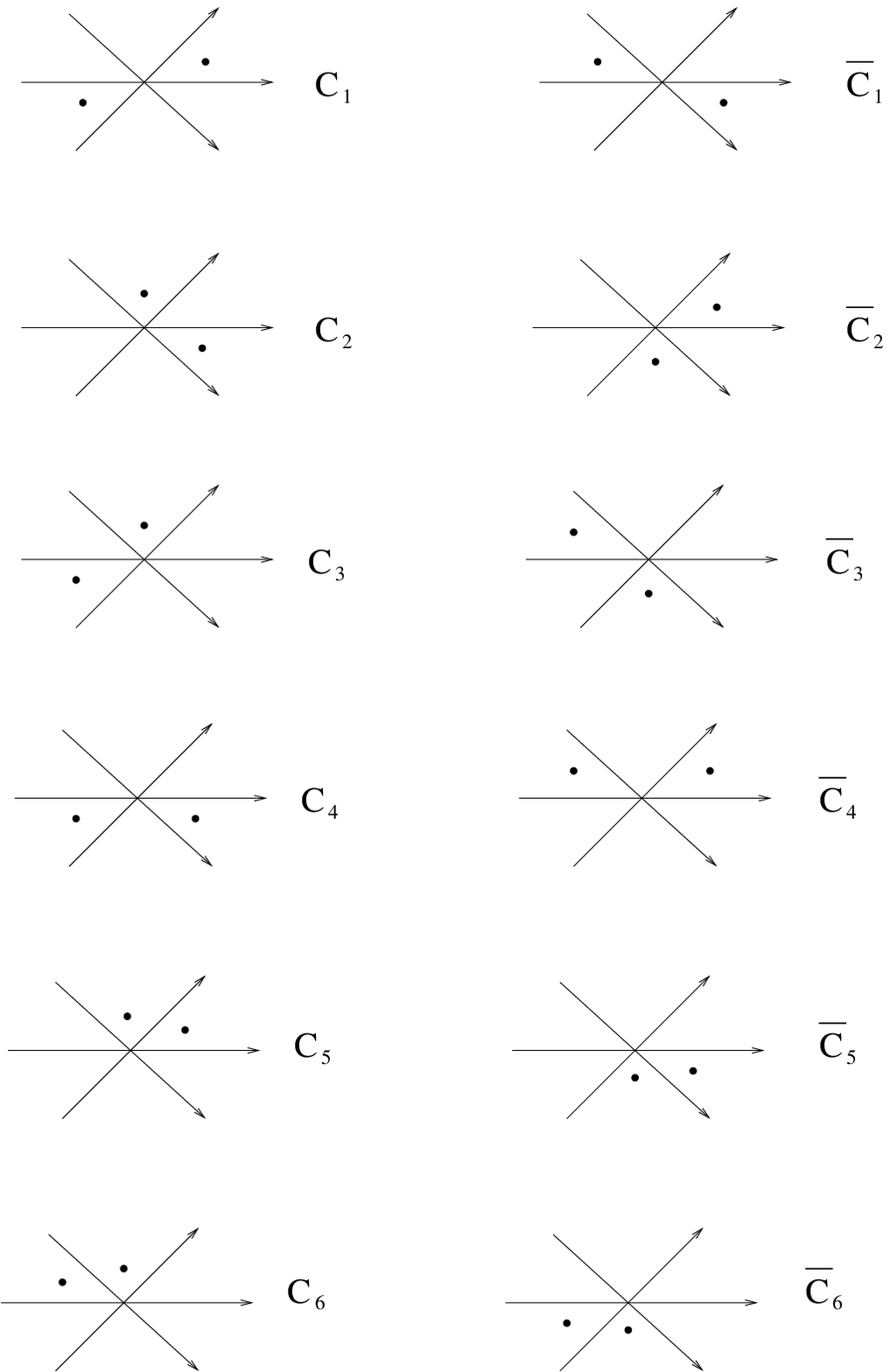}
\caption{}
\end{figure}

As in Jones case there are six couples of independent variables, and which we denote by the same symbols.
(Notice however, that there seems to be no natural duality between the simplifications in Jones case and the markings
in Alexanders case.)

Let $p$ be a star-like triple point. We define the markings of $p$ in Fig. 21. Again, we need to consider diagrams in
\begin{figure}
\centering 
\psfig{file=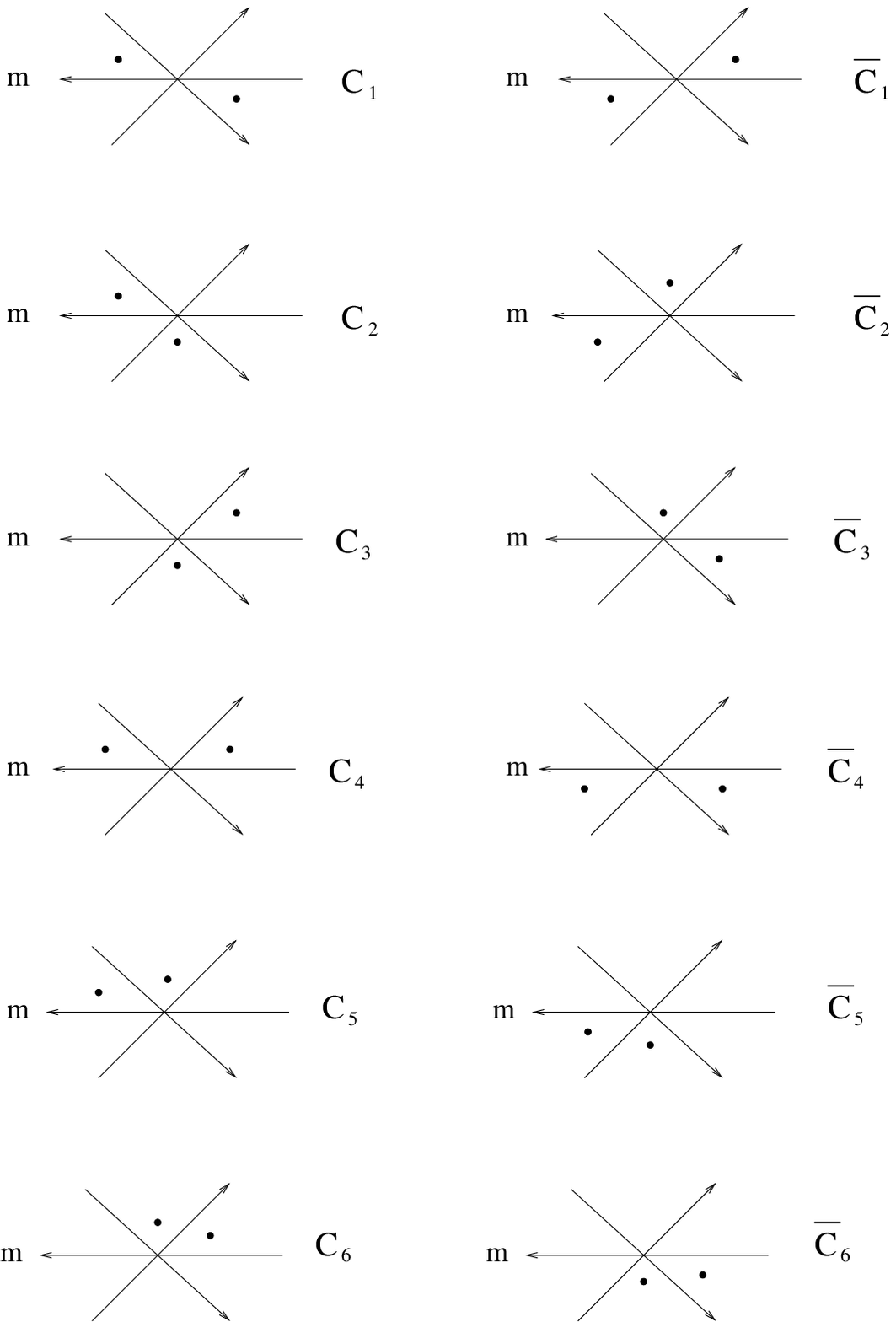}
\caption{}
\end{figure}
this definition instead of just planar curves. The middle branch is denoted by $m$ as previously. We end up with six 
couples of independent variables exactly as in the braid-like case.

For a marking $C_i$ (or $\bar C_i$) at a triple crossing $p$ we define the extended
Kauffman state sum, also called
$<D_p(C_i)>$ (respectively $<D_p(\bar C_i)>$, as usual
but without taking into account the two markings at the triple point.

\begin{remark}
Let $p$ be a positive triple point. We consider the marking $C_1$ (or likewise $\bar C_1$). Let us consider a nearby
generic diagram $D$ (no matter on which local side of the discriminant) and let $<D(C_1)>$ be the corresponding 
Kauffman state sum (i.e. at the three crossings coming from $p$ we consider only those markings which become $C_1$
when we shrink the triangle to the triple point $p$). One easily calculates that  $<D(C_1)> = 0$. Consequently, these 
Kauffman states do not contribute to the Alexander polynomial. However, they will contribute to our invariant $\Phi_K$.
\end{remark}

The markings at autotangencies are shown in Fig. 22. For a marking $B_i$ at an autotangency $p$ we define the
\begin{figure}
\centering 
\psfig{file=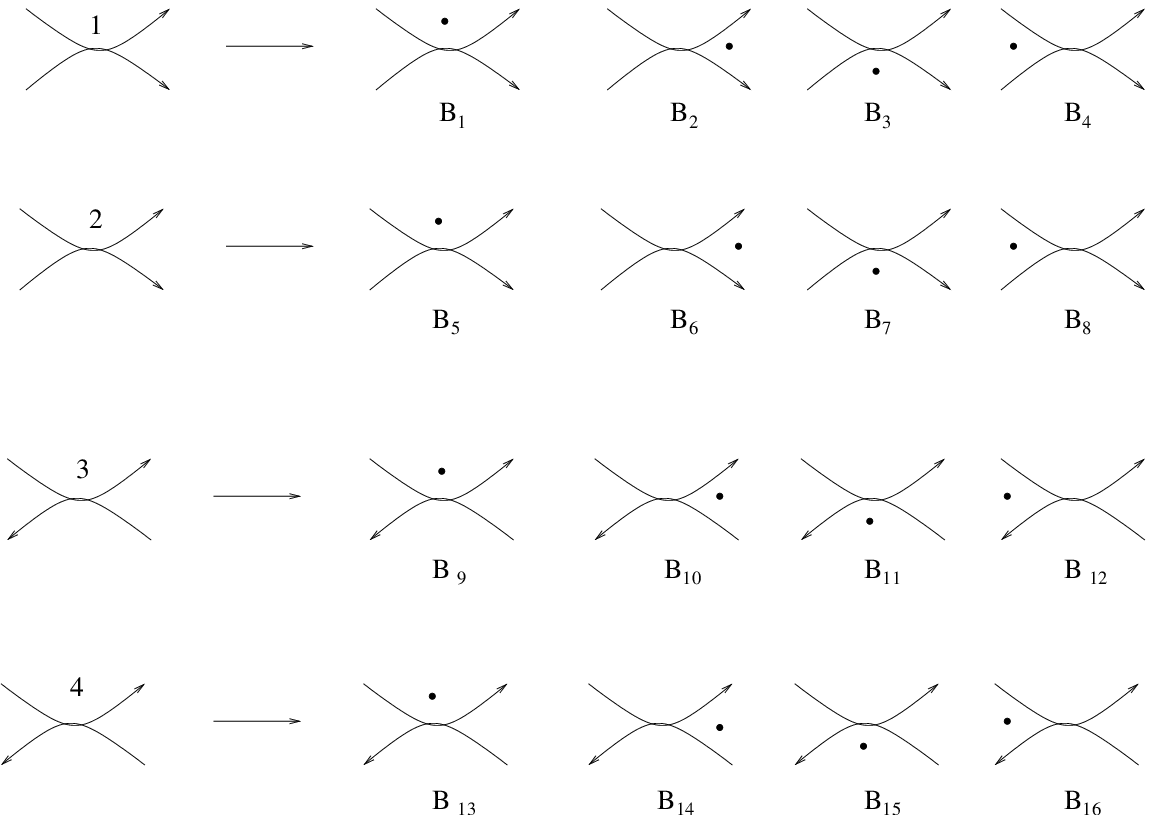}
\caption{}
\end{figure} 
extended
Kauffman state sum
$<D_p(B_i)>$ as usual,
but without taking into account the marking at the autotangency.

As in Jones case, the strata of $\Sigma_{f}^{(2)}$ give relations for the 
corresponding variables.  
Let $s$ be a meridian of $\Sigma_{f}^{(2)}$ and let $L(s)$ be the sum defined as in Lemma 4
but with the extended Kauffman brackets  replaced by the above Kauffman state sums in Alexanders case.
Let $X_2$ be the chord diagram of exactly two once intersecting chords in the disc.  
We distinguish all four endpoints of the chords on the boundary of the disc. Let $T_2$ be the free 
$\mathbb{Z}[A,A^{-1}]$-modul generated by the four {\em markings} of the chord diagram $X_2$ shown in Fig. 23.
\begin{figure}
\centering 
\psfig{file=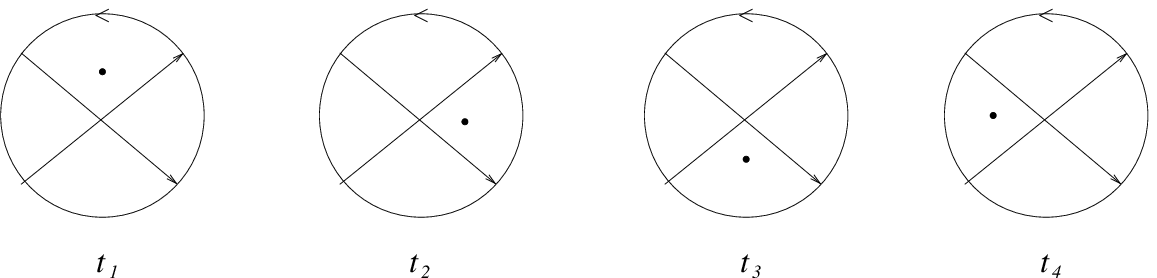}
\caption{}
\end{figure}
\begin{lemma}
$L(s) = 0$ if and only if 
$B_6 = - B_2$, $B_8 = - B_4$, $B_5 = A^{-2}B_1 - A^{-1}B_2 + A^{-1}B_4$, $B_7 = A^{-2}B_3 - A^{-1}B_2 + A^{-1}B_4$,
$B_{16} = - B_{12}$, $B_{14} = - B_{10}$, $B_{13} = - AB_{12} + AB_{10} - B_9$, $B_{15} = - A^{-1}B_{10} + A^{-1}B_{12} - B_{11}$.
\end{lemma}
{\em Proof:\/} We have to consider one stratum of $\Sigma_{f}^{(2)}$ where the tangent directions in the autotangency
coincide and another one where they are opposite. Let us consider the first case.
The two autotangencies in $s$ determine elements in $T_2$ as 
shown in Fig. 24. Evidently, $L(s) = 0$ independently of the rest of the diagram outside the marked diagrams in the
\begin{figure}
\centering 
\psfig{file=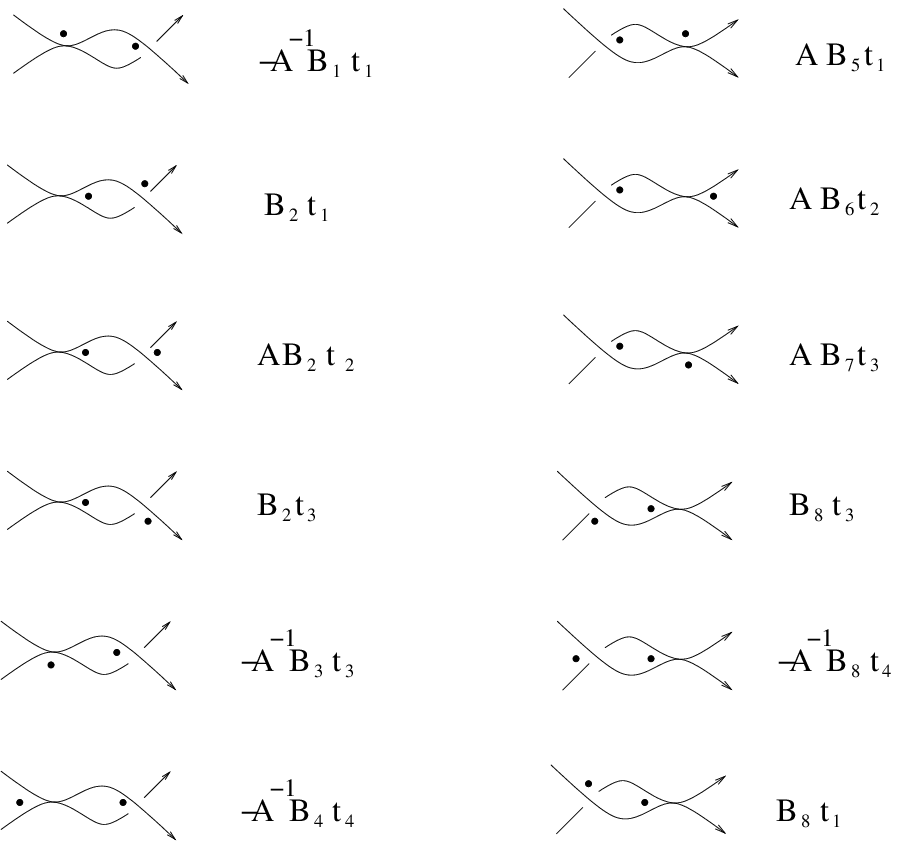}
\caption{}
\end{figure}
small disc  
if and only if for each generator of $T_2$ the signed sum of the coefficients is zero.  This gives us exactly the first
four required identities. The second case is  analogous and gives the remaining four identities.
$\Box$

\section{Main results}
In this section we define our five homomorphisms and give first examples and applications. There are four 
homomorphisms in Jones case ($S$, $S^+$, $S^-$, $X$) and one in Alexanders case ($\Phi$). They have a
very simple characterisation: $S$ is the only homomorphism which uses only triple points and only non singular
simplifications. $S^+$ and $S^-$ are the only homomorphisms which use only non singular simplifications of both, triple
points and of autotangencies. $X$ is the only homomorphism which uses only singular simplifications of both, triple
points and of autotangencies. All this will be proved in Section 4, where we show that the above homomorphisms are 
the {\em only}  solutions of the tetrahedron and the cube equations in Jones case.

$\Phi$ is the only homomorphism which uses only triple points in Alexanders case. However, we have not checked 
wether there are
other homomorphisms which use in addition also autotangencies.

\subsection{The homomorphism $S$}

It turns out that the tetrahedron equation has exactly two solutions in Jones case: one which uses only smoothings
and one which uses only simplifications with exactly one double point. The homomorphism $S$ will be
defined by using only smoothings (the letter $S$ stands for "smoothing").

We have already almost all ingridients (compare Definition 6) in order to give the definition of $S$.

\begin{definition}
The {\em weight function} $f_S$ is defined as follows: let $j \in \{1, ...,8\}$ be the type of the triple point.
Then $f_S(1) = f_S(2) = 1$, $f_S(3) = f_S(4) = A^4$, $f_S(5) = f_S(7) = -A^2$, $f_S(6) = f_S(8) = -A^6$.
\end{definition}

\begin{definition}
Let $\gamma$ be a generic oriented loop in $M$. Then $S(\gamma) \in  \mathbb{Z}[A,A^{-1},h,t,t^{-1}]$ is defined
as follows:

$S(\gamma) = \sum_{p \in \gamma}sign(p)[<D_p(C_1)> - <D_p(\bar C_1)>]f_S(j(p))A^{-8n(d(p))}t^{[d(p)]}$.

Here the sum is over all intersections $p$ of $\gamma$ with $\Sigma^{(1)}_t$, i.e. diagrams with an ordinary triple 
crossing, and $t$ is an independent new variable. For the definition of the brackets compare Subsection 2.6.
\end{definition}

\begin{remark}
Surprisingly, the Whitney index in the annulus $n(d(p))$ of certain subdiagrams enter in the definition of $S$.
\end{remark} 

\begin{theorem}
$S(\gamma)$ depends only on the homology class of $[\gamma] \in H_1(M)$ and it is zero on all sliding classes.
Consequently, $S: H_1(M) \rightarrow  \mathbb{Z}[A,A^{-1},h,t,t^{-1}]$ is a homomorphism and $S_K = S([rot])$ is 
well defined.
\end{theorem}
\subsection{The homomorphisms $S^+$ and $S^-$}
The homomorphisms $S, S^+, S^-$ come all from the same solution of the tetrahedron equation (using only smoothings)
, but from different
solutions of the cube equations. In the definition of $S$ we have used only triple points. In the definition of $S^+$ and
$S^-$ we use also autotangencies. The surprising result is that $n(d(p))$ enters no longer in the definition and that 
the types of triple points splitt into two families. Each of these families gives rise to an invariant.

\begin{definition}
The {\em relatives of the positive triple point $1$} are the types $2, 5, 7$. The {\em relatives of the negative triple
point $8$} are the types $3, 4, 6$.
\end{definition}

\begin{definition}
Let $j \in \{1, ...,4\}$ be the type of the autotangency.

The {\em weight function} $f_+$ is defined by $f_+(1) = A^3$, $f_+(2) = - A^3$, $f_+(3) = A^{-1}$, $f_+(4) = - A^{-1}$.
The {\em weight function} $f_-$ is defined by  $f_-(j) = 1/f_+(j)$.
\end{definition}

\begin{definition}
Let $\gamma$ be a generic oriented loop in $M$. Then $S^+(\gamma) \in  \mathbb{Z}[A,A^{-1},h,t,t^{-1}]$   is defined as follows:

$S^+(\gamma) = (A^{-2} -  A^2) \sum_{p \in \gamma}sign(p)[<D_p(C_1)> - <D_p(\bar C_1)>]f_S(j(p))t^{[d(p)]} 
+ \sum_{q \in \gamma}sign(q)<D_q(C_0)>f_+(j(q))t^{[d(q)]}$

Here, the first sum is only over all positive triple points and their relatives $p$ and the second sum is over all 
autotangencies $q$.

$S^-(\gamma) = (A^{-2} -  A^2) \sum_{p \in \gamma}sign(p)[<D_p(C_1)> - <D_p(\bar C_1)>]f_S(j(p))t^{[d(p)]}
+ \sum_{q \in \gamma}sign(q)<D_q(C_0)>f_-(j(q))t^{[d(q)]}$

Here, the first sum is only over all negative triple points and their relatives $p$ and the second sum is over all 
autotangencies $q$.

\end{definition}

\begin{theorem}
$S^{+(-)}(\gamma)$ depend only on the homology class of $[\gamma] \in H_1(M)$ and they are zero on all sliding classes.
Consequently, $S^{+(-)}: H_1(M) \rightarrow  \mathbb{Z}[A,A^{-1},h,t,t^{-1}]$ are homomorphisms and $S^+_K$ and
$S^-_K$ are well defined.

\end{theorem}

\subsection{The homomorphism $X$}
The homomorphism $X$ (the letter $X$ stands for "a double point") corresponds to the other solution of the tetrahedron
 equation in Jones case. In order to find 
a solution which satisfies also the cube equations we have to consider triple points and autotangencies.

\begin{definition}
The {\em weight function} $f_i(j)$ is defined as follows: let $j \in \{1, ...,8\}$ be the type of the tripel point and let
$C_i$ with $i \in \{3, 4, 5\}$ be the type of the simplification. Then

$f_3(1) = - A^4$, $f_3(3) =f_3(4) = - 1$, $f_3(5) = f_3(7) = A^6$, $f_3(8) = A^2$.

$f_4(1) = 1$, $f_4(3) = f_4(4) = A^4$, $f_4(5) = f_4(7) = - A^2$, $f_4(8) = - A^6$.

$f_5(j)$ is identical to $f_4(j)$.

$f_3(2) = f_4(2) = f_5(2) = A^4$.

$f_3(6) = f_4(6) = f_5(6) = A^2$.

\end{definition}

Notice, that triple points of type $3$ and $4$ (respectively type $5$ and $7$) have always the same value 
of the weight function. So we can treat the couples as triple points of the same type.

\begin{definition}
The {\em weight function} $f_X(j)$ is defined as follows: let $j \in \{1, ...,4\}$ be the type of the autotangency. 
We consider only the simplification $C_2$ with a double point. Then

$f_X(1) = A^5 + A$, $f_X(2) = - A^5 - A$, $f_X(3) = f_X(4) = 0$.
\end{definition}

\begin{definition}
Let $\gamma$ be a generic oriented loop in $M$.  
Then $X(\gamma) \in  \mathbb{Z}[A,A^{-1},h,r,s]$ is defined
as follows:

$X(\gamma) = \sum_{p \in \gamma}sign(p) \sum_{i \in \{3, 4, 5\}}[<D_p(C_i)> - <D_p(\bar C_i)>]f_i(j(p))
+ \sum_{q \in \gamma}sign(q)<D_q(C_x)>f_X(j(q))$.

Here, the first sum is over all triple points $p$ and the second sum is over all 
autotangencies $q$ (for the brackets compare Definition 9).

\end{definition}

The polynomials $X_K$ are evidently in general of the follwing form:

$X_K = \sum_i a_ih^i  + \sum_{j,k}a_{j, k}r^js^k$, where $i, j, k$ are natural numbers and $a_i, a_{j, k}$ are
Laurent polynomials in $A$ with integer coefficients. We suppose that in this decomposition $(j, k) \not= (0, 0)$.

\begin{lemma}
The coefficients of the polynomials $a_i$ are even for each $i$ and $a_{j, k} = a_{k,j}$ for each $j, k$.
\end{lemma}
{\em Proof:\/} Let $p$ be a triple point. Then $rot_{\pi}(p)$ is also a triple point in $rot(K)$, but of opposite sign.
Either the two triple points are of the same type or we can treat them as if they were of the same type. 
Let $C_m$ be a simplification of $p$ and let $<D_p(C_m)>$ be the corresponding Kauffman bracket. $rot_{\pi}$ maps
each Kauffman state of $D_p(C_m)$ to one of $D_{rot_{\pi}(p)}(\bar C_m)$. The two signed states give the same 
contribution to
$X_K$ besides that $j$ and $k$ are interchanged by $rot_{\pi}$. The same is true for autotangencies.
The lemma follows now easily.
$\Box$

We use the above lemma in the following definition.

\begin{definition} With the above notations we define

$X_K = 1/2\sum_i a_ih^i  + \sum_{j,k}a_{j, k}r^js^k$.
\end{definition}

\begin{theorem}
Let $\gamma$ be a generic oriented loop in $M^{braid}$.
Then $X(\gamma)$ depends only of the homology class $[\gamma] \in H_1(M^{braid})$.

Let $\gamma$ be a generic oriented loop in $M$. Then $X(\gamma) \in \mathbb{Z}/2\mathbb{Z}[A,A^{-1},h,r,s]$
 depends only of the homology class $[\gamma] \in H_1(M ;\mathbb{Z}/2\mathbb{Z})$. $X$ is trivial on sliding classes.
Consequently, $X_{\hat \beta} \in \mathbb{Z}[A,A^{-1},h,r,s]$ and $X_K \in \mathbb{Z}/2\mathbb{Z}[A,A^{-1},h,r,s]$
are well defined invariants of closed braids respectively of knots.

\end{theorem}

\subsection{The homomorphism $\Phi$}
The homomorphism $\Phi$ is constructed by using only triple crossings.

\begin{definition}
The {\em weight function} $g_i(j)$ is defined as follows: let $j \in \{1, ...,8\}$ be the type of the tripel point and let
$C_i$ with $i \in \{1, 2, 3, 4\}$ be the type of the simplification. Then

$g_i(1) = g_i(2) = g_i(3) = g_i(4) = 1$ and  $g_i(5) = g_i(6) = g_i(7) = g_i(8) = - 1$ for $i \in \{2, 3, 4\}$.

$g_1(1) =g_1(3) =g_1(4) = 1$, $g_1(5) =g_1(7) =g_1(8) = - 1$, and $g_1(2) = g_1(6) = 0$.

\end{definition}

\begin{definition}
Let $\gamma$ be a generic oriented loop in $M$. 
Then $\Phi(\gamma) \in  \mathbb{Z}[A,A^{-1}]$ is defined
as follows:

$\Phi(\gamma) = (A + A^{-1})\sum_{p \in \gamma}sign(p)[<D_p(C_1)> - <D_p(\bar C_1)>]g_1(j(p))
+ \sum_{p \in \gamma}sign(p) \sum_{i \in \{2, 3, 4\}}[<D_p(C_i)> - <D_p(\bar C_i)>]g_i(j(p))$.

\end{definition}

\begin{theorem}
Let $\gamma$ be a generic oriented loop in $M$.  
Then $\Phi(\gamma) \in  \mathbb{Z}[A,A^{-1}]$ depends only on the homology class 
$[\gamma] \in H_1(M)$ and it is zero on all sliding classes. Consequently, $\Phi$ is a well
defined homomorhism and $\Phi_K$ is a knot invariant.

\end{theorem}

\subsection{Some generalizations and refinements in the case of closed braids}
The case of closed braids in the solid torus is rather special for several reasons: there are no triple points of the types
 $2$ and $6$ and there are no autotangencies of the types $3$ and $4$. Consequently, several 2-faces of the cube are
no longer relevant. Besides the relations from double edges only the {\em long cycle in the cube}, which connects all 
the six remaining vertices, could give a relation for our variables $C_i(j)$.  
 Moreover, the Whitney index in the annulus of each crossing is zero.

As already stated in Theorem 4, the homomorphism $X$ can be lifted to a homomorphism with values in integer 
Laurent polynomials in the case of closed braids.

There are four reasons that our homomorphisms in general are only defined for knots and not for links:

a) there is no canonical choice of a meridian for a link

b) triviality on sliding loops could be only proved for knots (in this case we slide the curl over each crossing twice,
in difference to the case of links).

c) the correction factor with $n(d(p))$ is not well defined for links

d) the homological marking $[d(p)]$ is defined only if the two branches which cross at $d(p)$ belong to the same component.

(For all this compare the next section.)

Closed braids are already in the solid torus and hence, there is no problem with a). There are no sliding loops and 
$n(d(p))$ is always zero. Consequently, there are no problems with b) and c) neither.

\begin{definition}
Let $d(p)$ be a distinguished crossing of a triple point. If the two branches of the diagram, which cross at $d(p)$,
are in the same component of the oriented link, then we define $[d(p)]$ as previously. Otherwise we set $[d(p)] = 0$.
\end{definition}

We define now $S, S^+, S^-$ for links which are closed braids in the solid torus by the same expressions as previously,
but by using Definition 21 for $[d(p)]$.

The Theorems 2, 3, 4 and 5 together with the above considerations imply immediately the following proposition.

\begin{proposition}
The homomorphisms 

$S, S^+, S^- : H_1(M^{braid}) \rightarrow  \mathbb{Z}[A,A^{-1},h,t,t^{-1}]$ 

are well defined
for all closed braids. The same is true for the homomorphisms
 
$X : H_1(M^{braid}) \rightarrow \mathbb{Z}[A,A^{-1},h,r,s]$ and 
$\Phi : H_1(M^{braid}) \rightarrow \mathbb{Z}[A,A^{-1}]$.

\end{proposition}

Moreover, in the case of closed braids (and of those almost closed braids which are knots) there is a striking refinement of 
our  invariants by using the trace graph $TG(\hat \beta)$ (see \cite{FK} and \cite{F2}). It was shown in \cite{FK} that the pre-canonical 
loop of a knot is never tangential to $\Sigma^{(1)}_a$. It follows that the trace graph of a closed braid or an almost closed braid 
never changes by a Morse modification of index 1 (compare \cite{F2}). Consequently, the non contractible components 
of (a resolution of) the trace graph, which is associated to an arbitrary multiple of the canonical loop, are invariants 
of the closed braid or almost closed braid. All points in a component of the trace graph correspond to crossings of closed braids
with the same homological marking. However, for non trivial multiples of the canonical loop it really can happen that there
are {\em different} components of the trace link which correspond to the {\em same} homological marking of crossings
(se \cite{F2}). Consequently, if we replace in the definition of $S, S^+, S^-$ the homological markings $[d(p)]$ by the 
corresponding components of the trace graph which contain $d(p)$, then we obtain finer information.

The following lemma is implicitely contained in \cite{FK}.

\begin{lemma}
Two canonical loops in $M^{braid}$ are homologic if and only if they are homotopic without ever being tangential to
$\Sigma^{(1)}_a$. The same is true for their multiples. 
\end{lemma}

Let $\hat \beta$ be a closed braid, let $rot(\hat \beta)$ be the corresponding canonical loop and let $l$ be any  fixed integer.
We consider the loop $l(rot(\hat \beta))$, i.e. we go $l$ times along $rot(\hat \beta)$. Let $T$ be any fixed component of the trace graph $TG(l(rot(\hat \beta)))$.

\begin{definition}
The polynomials $S, S^+ ,S^-(l(rot(\hat \beta)))$ are defined as previously, but by summing {\em only} over those triple points 
and autotangencies which have $d(p) \in T$. We denote the corresponding polynomials by $S_{\hat \beta, l, T}$,
$S^+_{\hat \beta, l, T}$, $S^-_{\hat \beta, l, T}$.

\end{definition}

The following proposition is an immediate consequence of the Theorems 2 and 3 together with the above considerations.
\begin{proposition}
$S_{\hat \beta, l, T}, S^+_{\hat \beta, l, T}, S^-_{\hat \beta, l, T} \in  \mathbb{Z}[A,A^{-1},h,t,t^{-1}]$ are isotopy invariants of 
closed braids.

\end{proposition}

\begin{remark}
Notice, that we need no longer to distinguish wether the two branches at $d(p)$ belong to the same component of the link
or not. 

The trace graph depends in a non trivial way from $l$, compare \cite{F2}. As a consequence, the invariants from Proposition 4 behave non linear with 
respect to $l$.

\end{remark}

We left the elaboration of the refined invariants in the case of almost closed braids to the reader.

Finally, let us mention that there is another refinement of our invariants in the very special case of closed 3-braids.
Indeed, there are never quadruple points in isotopies of closed 3-braids. Consequently, we do not need to solve the 
tetrahedron equation but only the cube equations. Therefore, e.g. from the definition of $\Phi$ the sum 
$(A + A^{-1})\sum_{p \in \gamma}sign(p)[<D_p(C_1)> - <D_p(\bar C_1)>]g_1(j(p))$ and the sum
$\sum_{p \in \gamma}sign(p) \sum_{i \in \{2, 3, 4\}}[<D_p(C_i)> - <D_p(\bar C_i)>]g_i(j(p))$ are both already 
invariants. Let us call the first of these two sums $\Phi^{(1)}$ and the second $\Phi^{(2,3,4)}$.

For the same reason $X_i(\gamma) = \sum_{p \in \gamma}sign(p) \sum_i[<D_p(C_i)> - <D_p(\bar C_i)>]f_i(j(p))
+ \sum_{q \in \gamma}sign(q)<D_q(C_x)>f_X(j(q))$ is already an invariant for each $i \in \{3, 4, 5\}$.

\subsection{First examples and applications}
\begin{proposition}
Let $\hat \beta$ be a closed braid in the solid torus. Then all $S, S^+, S^-(\gamma)$ are polynomials of positive degree 
with respect to $t$.
\end{proposition}
{\em Proof:\/} Each subdiagram of a closed braid, which is obtained by smoothing a crossing with respect to the orientation, 
is again a closed braid. Consequently, it represents a positive homology class in $H_1(V)$.
$\Box$

It follows that the above proposition can be used to answer sometimes in the negative on the question, wether a given
knot in the solid torus is isotopic to any closed braid.

\begin{proposition}
Let $\hat \beta$ be a positive (negative) closed braid. Then for the reductions with $\mathbb{Z}/2\mathbb{Z}$-coefficients
 $S^-_{\hat \beta} = 0$ (respectively, $S^+_{\hat \beta} = 0$).
Moreover, $S^+_{\hat \beta} = (A^{-2} - A^2)S_{\hat \beta}\quad  (mod 2)$ (respectively, $S^-_{\hat \beta} = (A^{-2} - A^2)S_{\hat \beta}\quad  (mod 2)$).
\end{proposition} 
{\em Proof:\/} Let $\beta$ be a positive braid. We use the combinatorial canonical loop from Definition 4.
There are 
only Reidemeister II moves in
$\beta \rightarrow \Delta^2\Delta^{-2}\beta$. Pushing $\Delta^2$ through $\beta$ contains only positive triple points.
$\Delta^{-2}\Delta^2\beta \rightarrow \beta$ contains again only Reidemeister II moves. The latter are exactly the 
mirror images
of the Reidemeister II moves we started with. Consequently, they have the opposite weight (they are dual) and the 
opposite sign.
It follows that their contribution to $S^{+(-)}_{\hat \beta}$ has even coefficients. The triple points give no contribution to 
$S^-_{\hat \beta}$ because they are all positive. The proposition follows now directly from the definitions.
$\Box$

The proposition can be used to answer sometimes in the negative on the question, wether a given knot in the solid torus is 
isotopic to any positive (negative) closed braid.

Let $L$ be an oriented link in the standard solid torus $V$ in  $S^3$ and let $flip$ be the "hyperelliptic"
involution of the solid torus $V$ as shown in Fig.25. The link $L$ is called {\em invertible} in $V$ if $L$ is isotopic to
\begin{figure}
\centering 
\psfig{file=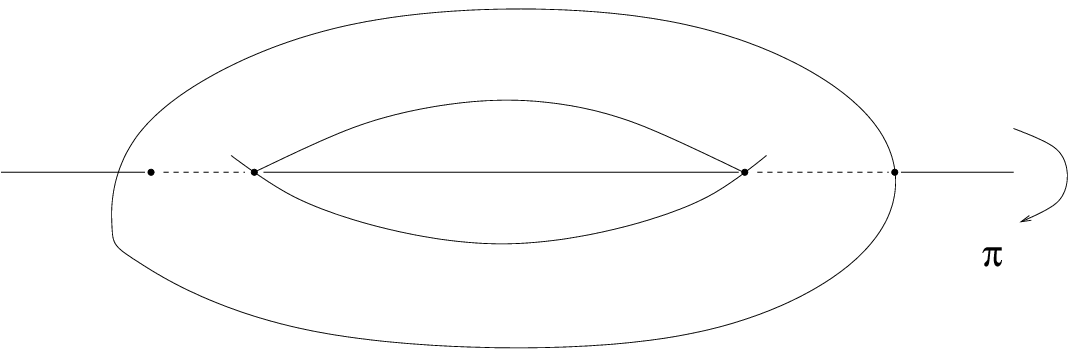}
\caption{}
\end{figure}
the inverse link $-flip(L)$ in $V$. Notice, that the link $L$ is invertible in $V$ if and only if the link 

$L \cup \{complementary\quad axes\quad of\quad V\quad in\quad S^3\}$ is invertible as a link in $S^3$. For classical links
this coincides with the usual definition of invertibility (compare
 \cite{F1} and \cite{F2}).

The behaviour of our invariants with respect to taking the mirror image or the inverse knot is a puzzling matter.
However, it becomes simpler in the case of positive closed braids. Let $\hat \beta$ be a positive closed braid.
We denote its mirror image by $\hat \beta!$ (i.e. all crossings are switched), and we denote the inverse closed braid 
by $-flip(\hat \beta)$. We will consider only Alexanders case.

\begin{proposition}
$\Phi_{\hat \beta!}(A) = - \Phi_{\hat \beta}(A^{-1})$ and $\Phi_{-flip(\hat \beta)}(A) =  - \Phi_{\hat \beta}(-A^{-1})$
.
\end{proposition}
{\em Proof:\/} Let us consider first the mirror image. Instead of $\Delta^2$ we push now $\Delta^{-2}$ through the braid $\beta!$.
Reidemeister II moves do not matter for $\Phi$ and there is a natural one to one correspondence between the triple points
for $\beta$ and $\beta!$. Let $p$ and $p!$ be any couple of corresponding triple points. All triple points for $\beta$
are positive (type $1$) and all triple points for $\beta!$ are negative (type $8$). For each marking $C_i$ of $p$ we consider
the same marking for $p!$. Then there is a natural 
one to one correspondence for the Kauffman states of $D_p(C_i(1))$ and those of $D_{p!}(C_i(8))$. Evidently,
this correspondence interchanges $A$ with $A^{-1}$ at the double points. Moreover, there are exactly three global sign 
changes:
$sign(p) = - sign(p!)$ (compare Fig. 6) and $g_i(1) = - g_i(8)$. But pushing $\Delta^{-2}$ through the braid instead of 
$\Delta^2$ corresponds to the opposite orientation of the canonical loop and consequently, changes one more time
the intersection indices of $rot(\hat \beta)$ with $\Sigma^{(1)}_t$. The first part of the proposition follows.

Let $p$ be a triple point for $\beta$ and let $-flip(p)$ be the corresponding triple point for $-flip(\beta)$. (In fact, 
$-flip(\beta)$ is just the braid $\beta$, but which has to be red backwards. If the braid is written horizontally, then 
we have to reflect it at a vertical line.) The involution $-flip$ maps Kauffman states to Kauffman states. One easily establishes the 
following (just consider vertical reflections and
orientation reversing in Fig.20):

$type(p) = type(-flip(p))$, $sign(p) = sign(-flip(p))$. If the marking of $p$ is of type $C_1$ then the corresponding 
marking of $-flip(p)$ is of type $\bar C_1$. If it is of type $C_2$ then the corresponding marking is of type $C_3$ and vice versa.
If it is of type $C_4$ then it stays invariant. The involution $-flip$ interchanges $A$ with $-A^{-1}$ at double points.

$g_2(1) = g_3(1)$ and $\bar C_1$ enters with a different sign as $C_1$ in $\Phi_K$, but they enter both with a factor
$A + A^{-1}$. We do not have to switch the sign of this factor.  

We are not yet done, because the meridian (or the axes of the complementary solid torus) has not the right orientation
after the flip.
Therefore, we have to switch the orientation of the annulus too. Notice, that this does not affect the minus signs at the
double points, neither the signs and the types of the triple points. But $C_i$ and $\bar C_i$ interchange for each $i \in \{1, 2, 3, 4\}$.

The second part of the proposition follows now easily.
$\Box$

We have calculated by hand three examples, all in the most simplest case, namely closed 3-braids. (Even in the case
of the most simple classical knots there are already many triple points from sliding curls along the knot in order to
performe Whitney tricks.) 

\begin{example}
Let $\beta = \sigma_1^3\sigma_2^5 \in B_3$. Then

$S_{\hat \beta} = ht^2(-2A^{-17} +A^{-13} -A^{-9} +2A^{-5} -A^{-1} +A^3) +ht(-2A^{-17} +A^{-13} -A^{-9} +2A^{-5} -A^{-1} +A^3)$.

(One easily sees, that the variables $h$ and $t$ are not interesting in the very special case of closed 3-braids.)

It follows, in particular, that the canonical loop $rot(\hat \beta)$ is not homologic to its inverse $-rot(\hat \beta)$.
 
Notice, that the reduction $mod 2$ of $S_{\hat \beta}$ is already non trivial. It follows from Proposition 6 that $S^+_{\hat \beta}$
is also non trivial. This implies of course, that $S^-_{\hat \beta!}$ is non trivial too.

\end{example}

\begin{example}
Let $\beta = \sigma_1^3\sigma_2^2 \in B_3$. As we have mentioned in the previous subsection that

$X_i(\gamma) = \sum_{p \in \gamma}sign(p) \sum_i[<D_p(C_i)> - <D_p(\bar C_i)>]f_i(j(p))
+ \sum_{q \in \gamma}sign(q)<D_q(C_x)>f_X(j(q))$  

for each $i \in \{3, 4, 5\}$ are already invariants. Notice, that the canonical loop $rot(\hat \beta)$ is invariant under 
the involution $rot_{\pi}$. It follows that already the h-part of each $X_i(rot(\hat \beta))$ was divisible by 2. We have only calculated  
the reduction $mod 2$ of the h-part coming from the triple points of 
$X_4(rot(\hat \beta))$. 

$X_4(rot(\hat \beta)) = h(A^{-8} +A^{-6} +A^{-5} +A^{-4} +A^{-1} +1 +A +A^3 +A^4) +\{terms\quad with\quad r\quad or\quad s\}  +
\{terms\quad coming\quad from\quad the\quad autotangencies\}$.

We will see in Subsection 4.6. that we can forget about the autotangencies in $X_4(rot(\hat \beta))\quad mod2$ 
if we add the supplementary relation $A^4 = 1$.

Consequently, the mod2 reduction of the h-part of $X_4(rot(\hat \beta))$ with $A^4 = 1$ reduces to $h(A  +A^2 +A^3)$ and is already 
non trivial.

\end{example}

\begin{example}
Let $\beta$ be the same as in Example 2. For the same reason 

$\Phi^{(1)}_{\hat \beta} = (A + A^{-1})\sum_{p \in rot(\hat \beta)}sign(p)[<D_p(C_1)> - <D_p(\bar C_1)>]g_1(j(p))$

is already an invariant. A calculation gives

$\Phi^{(1)}_{\hat \beta} = 2(A +A^{-1})(-3A^{-2} +4 -3A^2)$.

It follows from Proposition 7 that $\Phi^{(1)}_{\hat \beta!} = -  2(A +A^{-1})(-3A^{-2} +4 -3A^2)$. Hence, we have detected
with $\Phi$ that the closed braid is not amphichiral.

\end{example}

\section{Proofs}
In this section we show first how the 48 different tetrahedron equations reduce to a single equation plus the cube
equations. We give then a complete list of the solutions in Jones case and we give the unique solution of the 
tetrahedron equation in Alexanders case. Finally, we show that all our homomorphisms are trivial on sliding classes.
We finish the paper with a list of open questions.
\subsection{Reduction to the positive tetrahedron equation}
Let us remind the situation for one parameter families of diagrams. Each loop is homotopic to a loop which intersects 
$\Sigma_{t}^{(1)}$ only in strata which correspond to positive triple points. But the homotopy adds lots of intersections 
with
$\Sigma_{a}^{(1)}$ (compare \cite{F1}). the replacement of one type of triple points by another type is a local operation. 
It follows that all triple points can be replaced by only positive ones, because the graph 
$G$ is connected (compare Fig. 12).

Let us consider the generic degenerations of quadrupel points in 3-parameter families of diagrams.
Using the same methods as in the Appendix of \cite{FK} one proofs the following lemma (we left this to the reader).

\begin{lemma}
The boundary of $\Sigma_q^{(2)}$ in $\Sigma$ is the closure of only the following two types of strata of codimension 3
(in $M$):

(1) diagrams with an ordinary quintuple point

(2) diagrams with an ordinary autotangency through which two mutually transversal branches pass transversally.

At each point of (1) five strata of $\Sigma_q^{(2)}$ intersect mutually transversal. The intersection of a normal
3-disc of (2) with $\Sigma_q^{(2)}$ contains exactly two different types of strata as shown in Fig. 26.
\begin{figure}
\centering 
\psfig{file=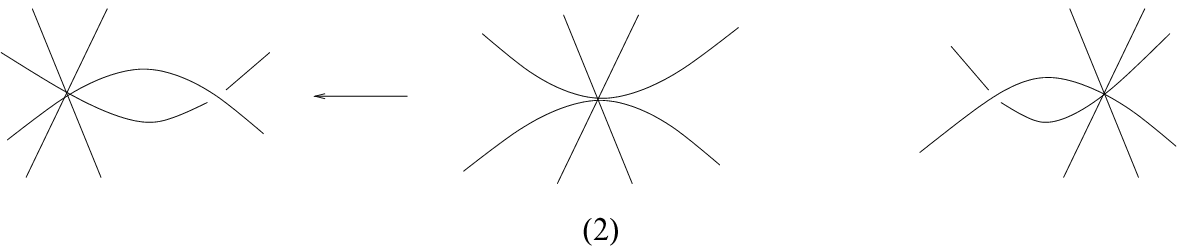}
\caption{}
\end{figure}
\end{lemma}

Notice, that this lemma is the exact analogue of the corresponding lemma for $\Sigma_{t}^{(1)}$ (compare Fig. 1 and 
Fig. 2).

 We proced now for two parameter families in a completely analogous way. 
Let $G'$ be the following graph: the vertices correspond exactly to the types of strata of $\Sigma_q^{(2)}$.
Two vertices are connected by an edge if and only if the two types of quadruple points come together in a stratum
of type (2) (see Lemma 11). 

The graph $G'$ is complicated but the following lemma is evident.

\begin{lemma}
The graph $G'$ is connected.
\end{lemma}

Using Lemmas 11 and 12 we replace by a 
homotopy (rel boundary) in each two parameter family of diagrams all intersections with $\Sigma_q^{(2)}$ 
by intersections which correspond only to positive quadruple points. Of course, this homotopy adds lots of intersections 
with $\Sigma_{a-t}^{(2)}$.

Consequently, in order to solve the tetrahedron equation (a) for arbitrary quadruple points (see the Introduction) it 
suffices to solve the 
tetrahedron equation only for the unique positive quadruple point and to solve all the cube equations. Notice, that we 
had to solve the cube equations in any case.

\begin{remark}
For each of the two quadruple points which come together in a stratum of type (2) there are exactly two triple points 
which share the same distinguished crossing (compare Remark 10). One easily sees, that these distinguished 
crossings for both quadruple points are in the same component of the trace graph. In particular, all four distinguished 
crossings have the same homological marking.
\end{remark}

\subsection{Solutions of the positive tetrahedron equation using the extended Kauffman bracket}
We consider only the 12 simplifications which are not invariant under the involution $rot_{\pi}$ (compare Fig. 15 and 16).
In Fig. 27 - 34 we show the 12 simplifications of the 8 triple points near the quadruple point. Outside of these figures all
 diagrams are identical.
\begin{figure}
\centering 
\psfig{file=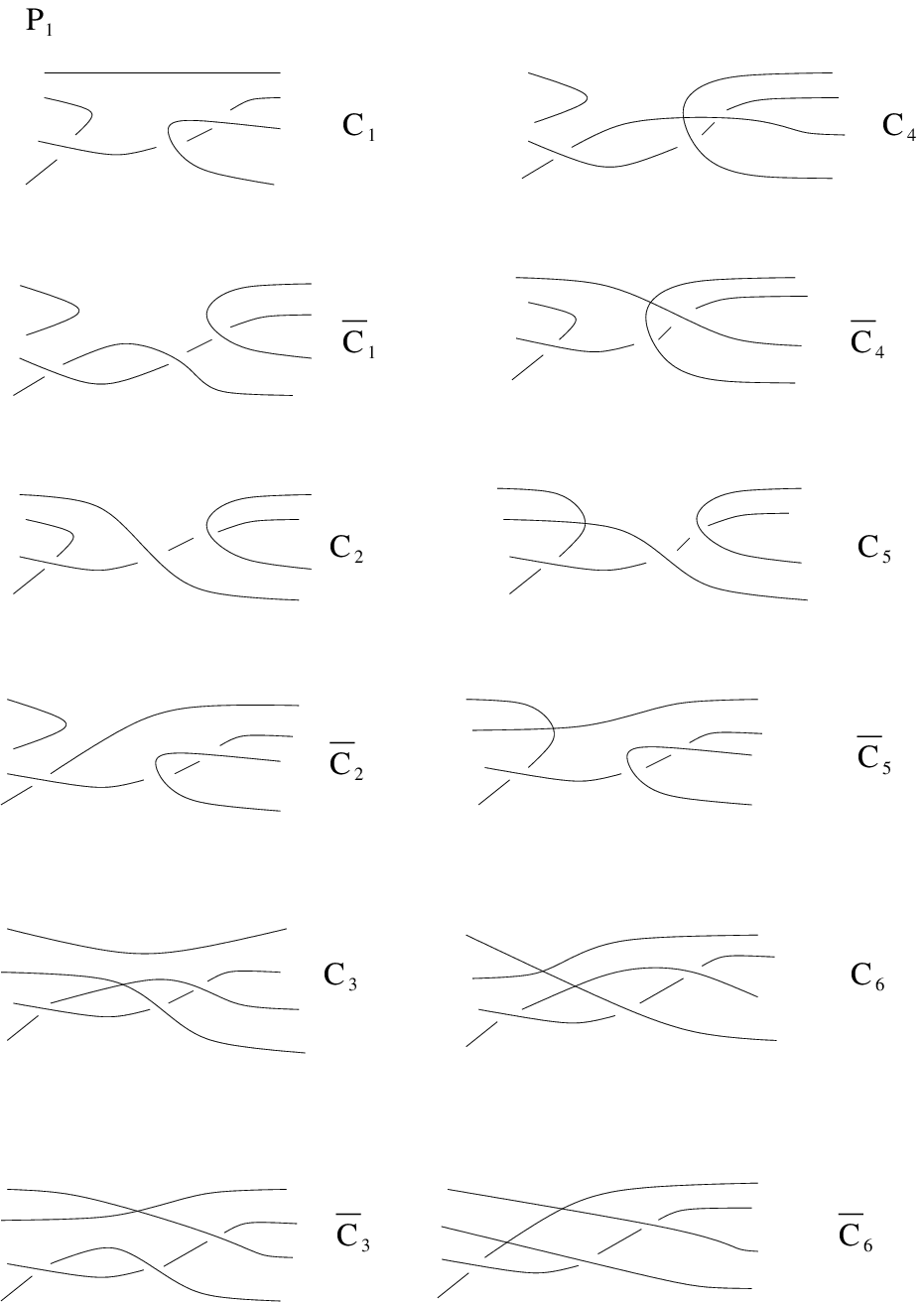}
\caption{}
\end{figure}
\begin{figure}
\centering 
\psfig{file=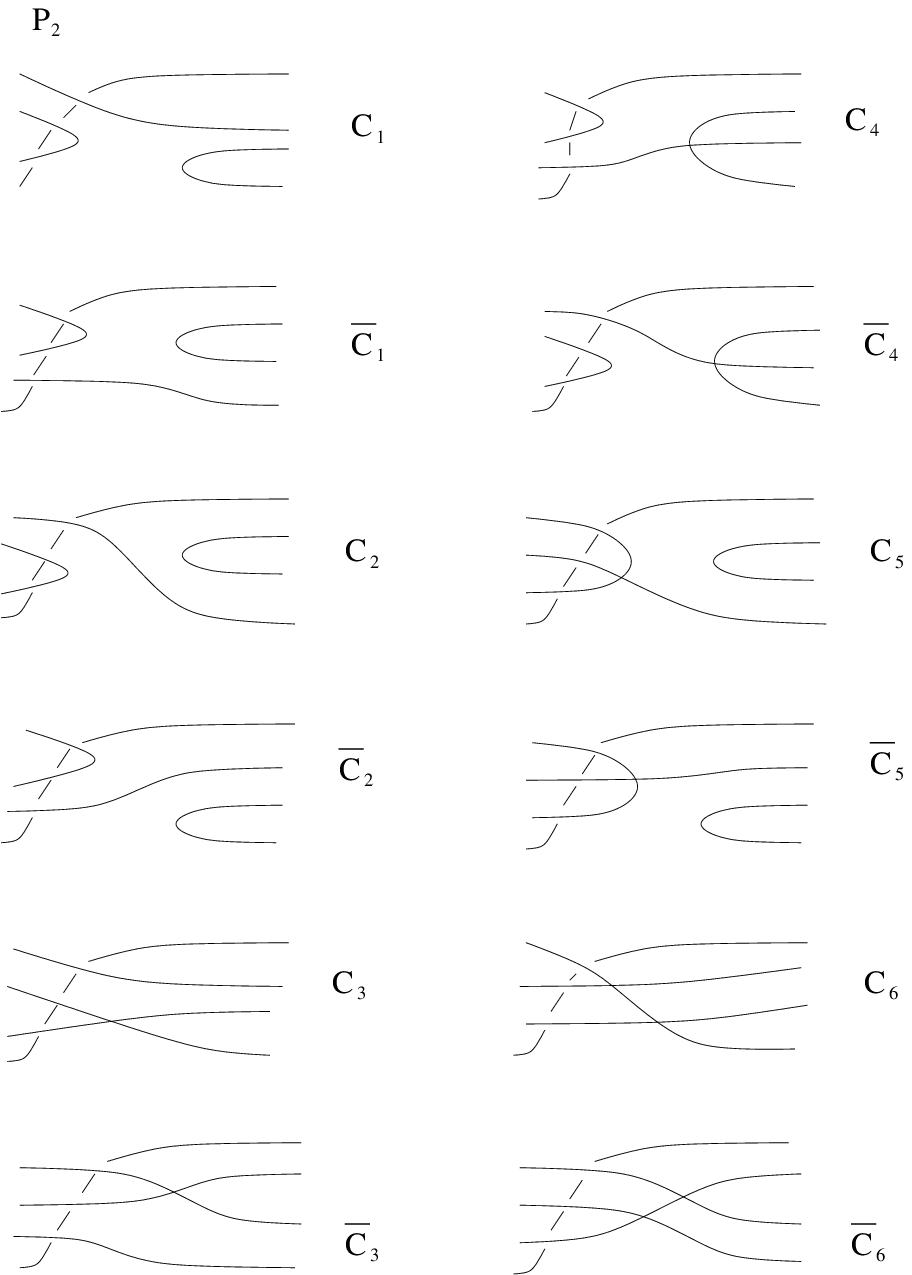}
\caption{}
\end{figure}
\begin{figure}
\centering 
\psfig{file=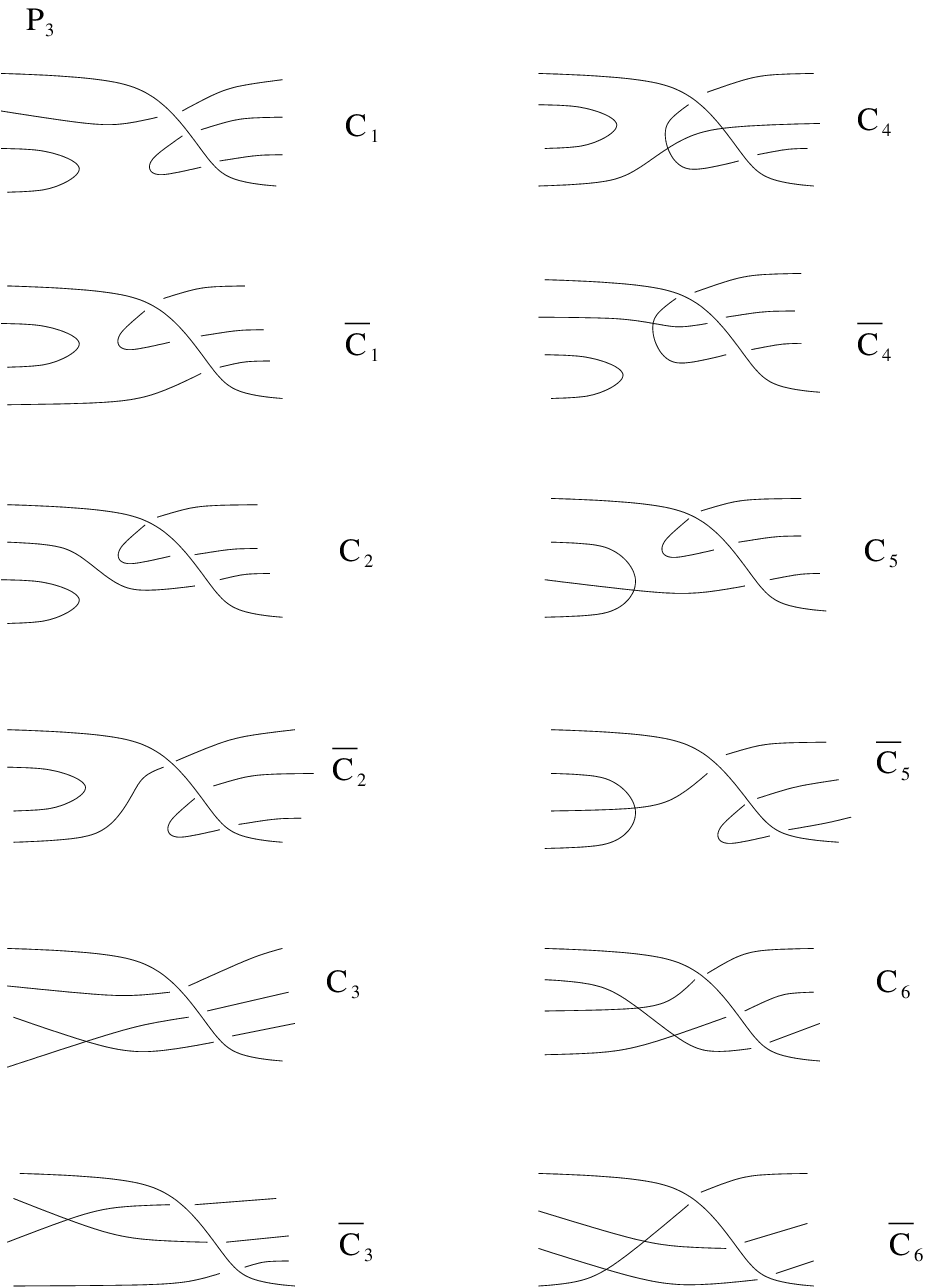}
\caption{}
\end{figure}
\begin{figure}
\centering 
\psfig{file=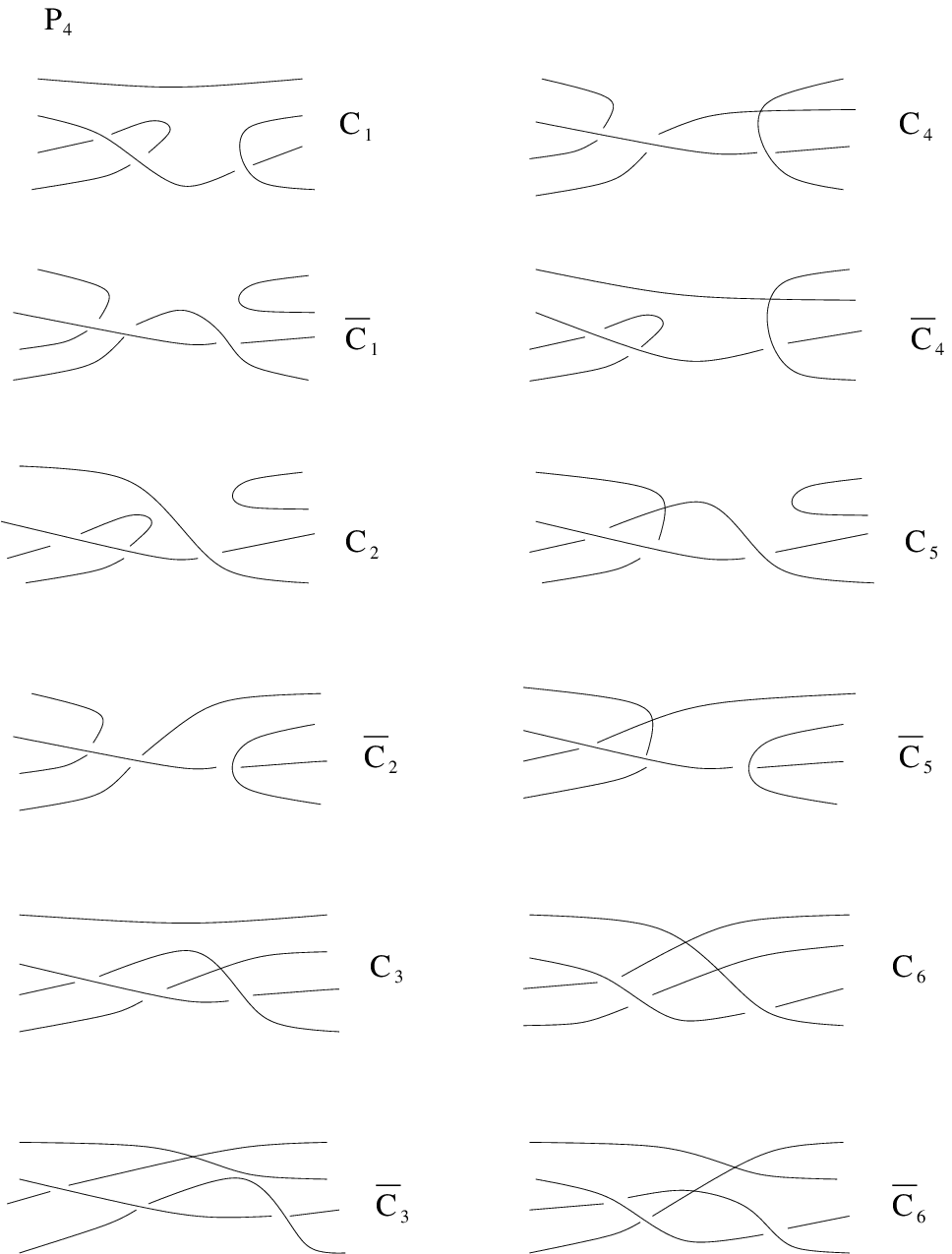}
\caption{}
\end{figure}
\begin{figure}
\centering 
\psfig{file=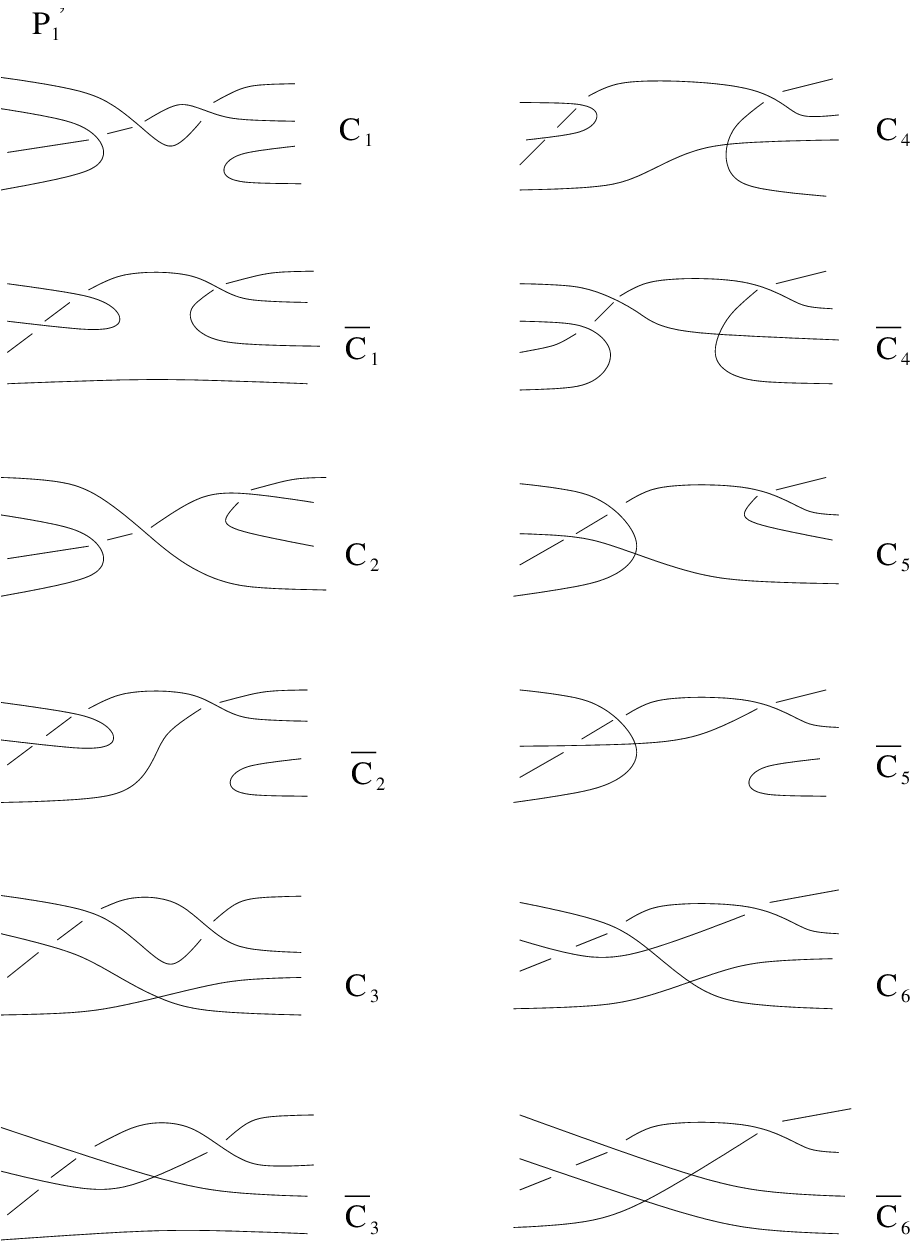}
\caption{}
\end{figure}
\begin{figure}
\centering 
\psfig{file=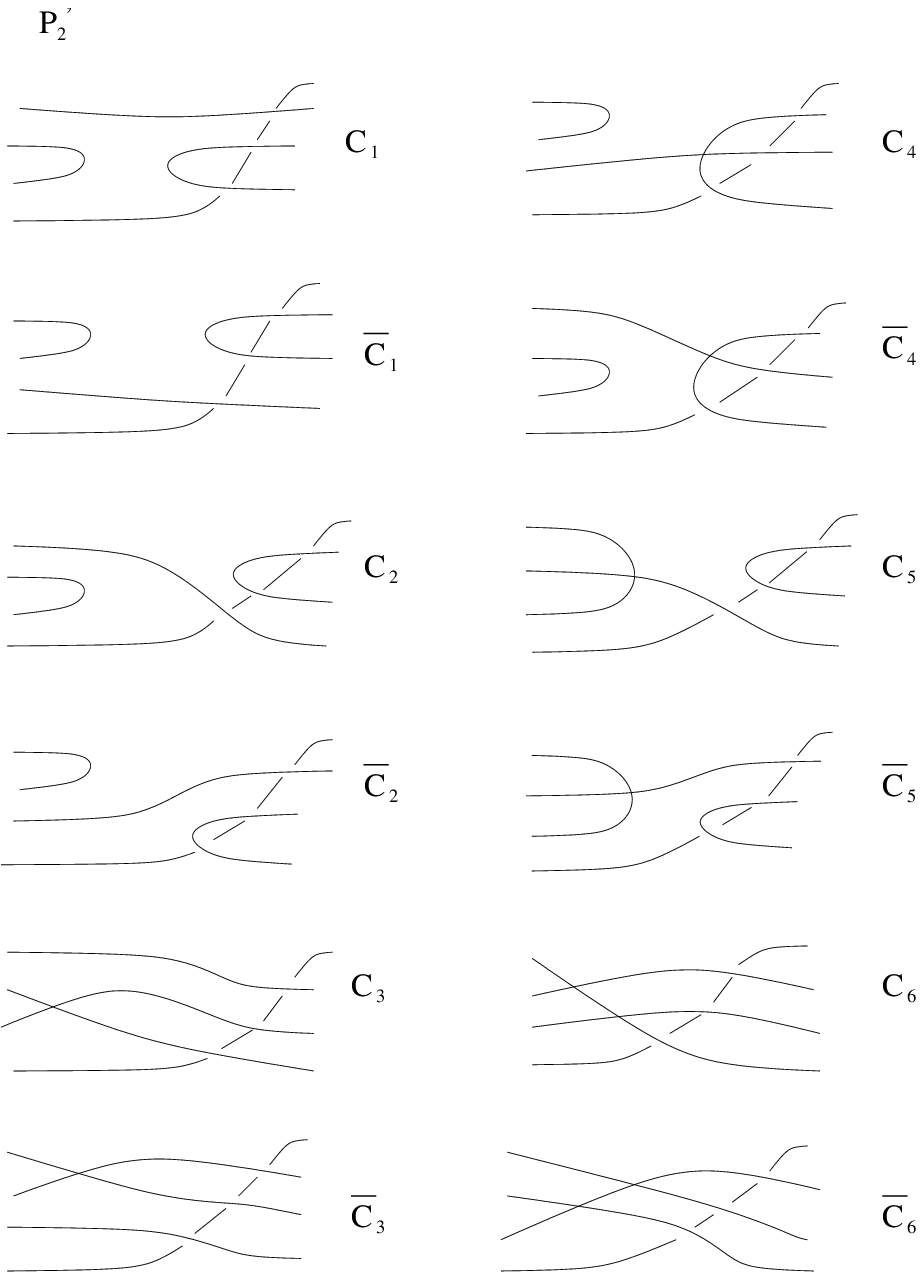}
\caption{}
\end{figure}
\begin{figure}
\centering 
\psfig{file=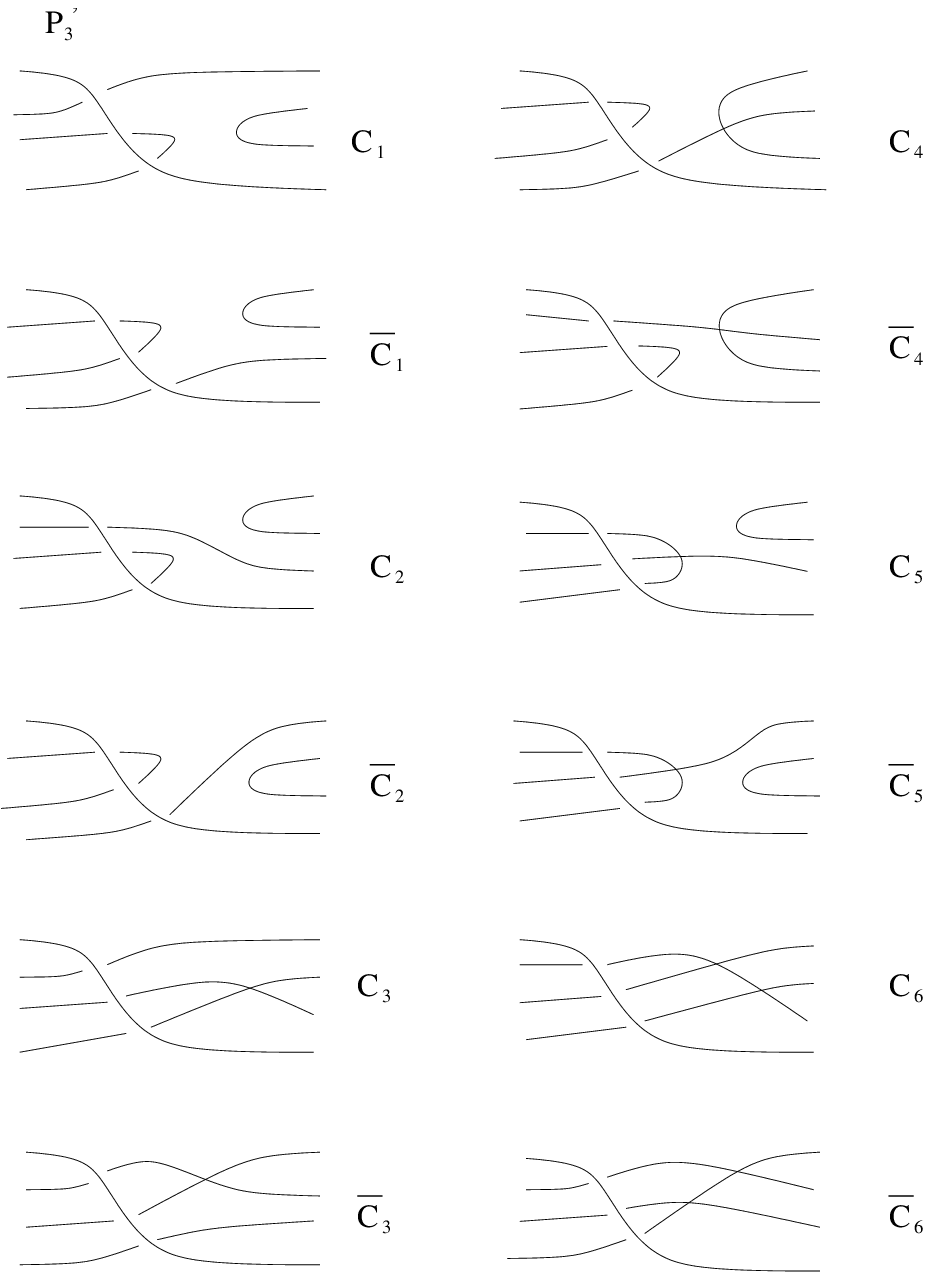}
\caption{}
\end{figure}
\begin{figure}
\centering 
\psfig{file=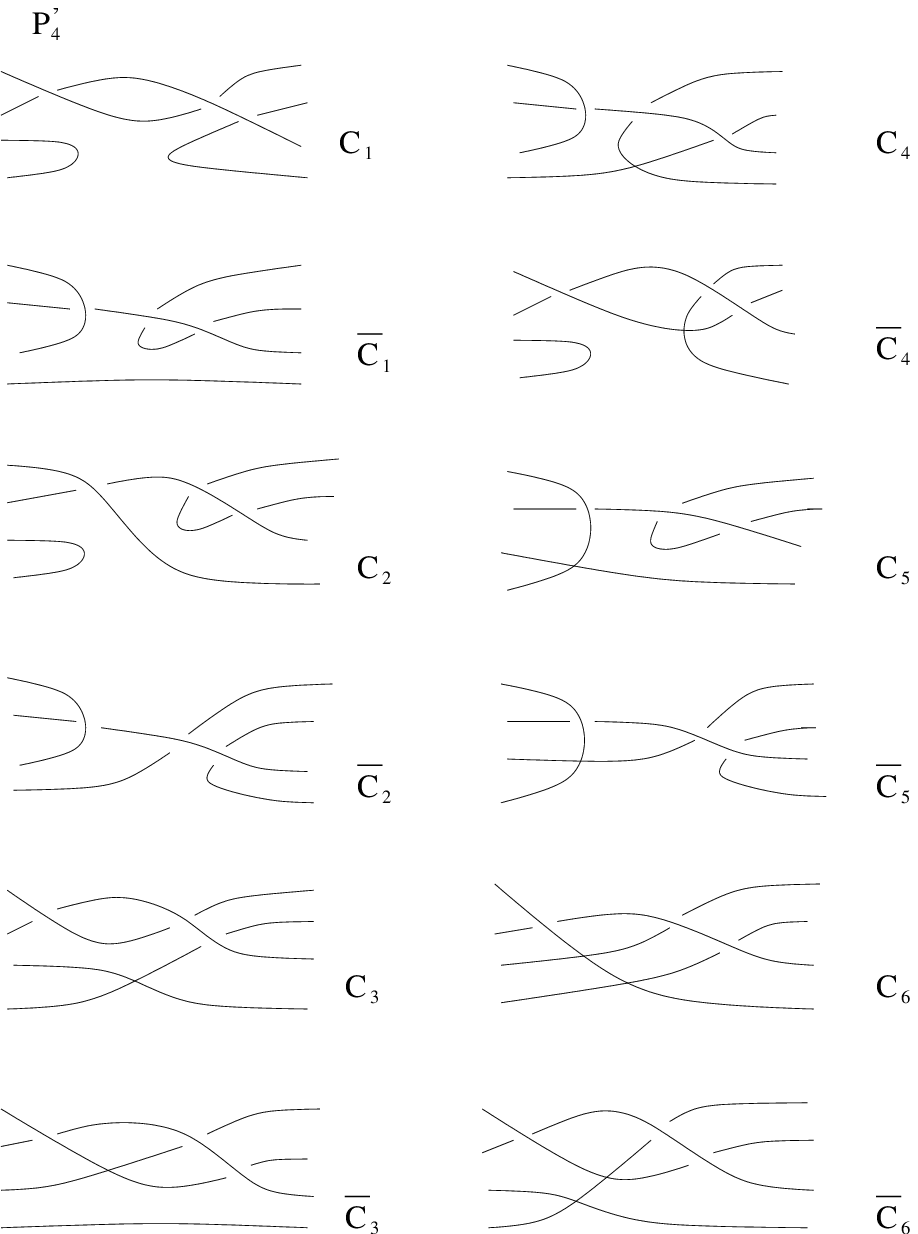}
\caption{}
\end{figure}
\begin{remark}
Notice, that for the smoothings $C_1$ and $\bar C_1$ we do not need to consider the triple points $2, 2', 3, 3'$!
Indeed, for each of these smoothings the diagram from 2 (respectively, 3) represents a link which is regularly isotopic 
to the link which is represented by the diagram from 2' (respectively, 3'). Consequently, they have identical extended 
Kauffman brackets, but which enter in the equation (c) with opposite signs.
\end{remark}

Let $T_4$ be the free $\mathbb{Z}[A,A^{-1}]$-modul associated to the quadruple point as introduced in the Introduction.
We give names to the generators of $T_4$ in Fig. 35 - 37. (The remaining generators do never occure and they are therefore
\begin{figure}
\centering 
\psfig{file=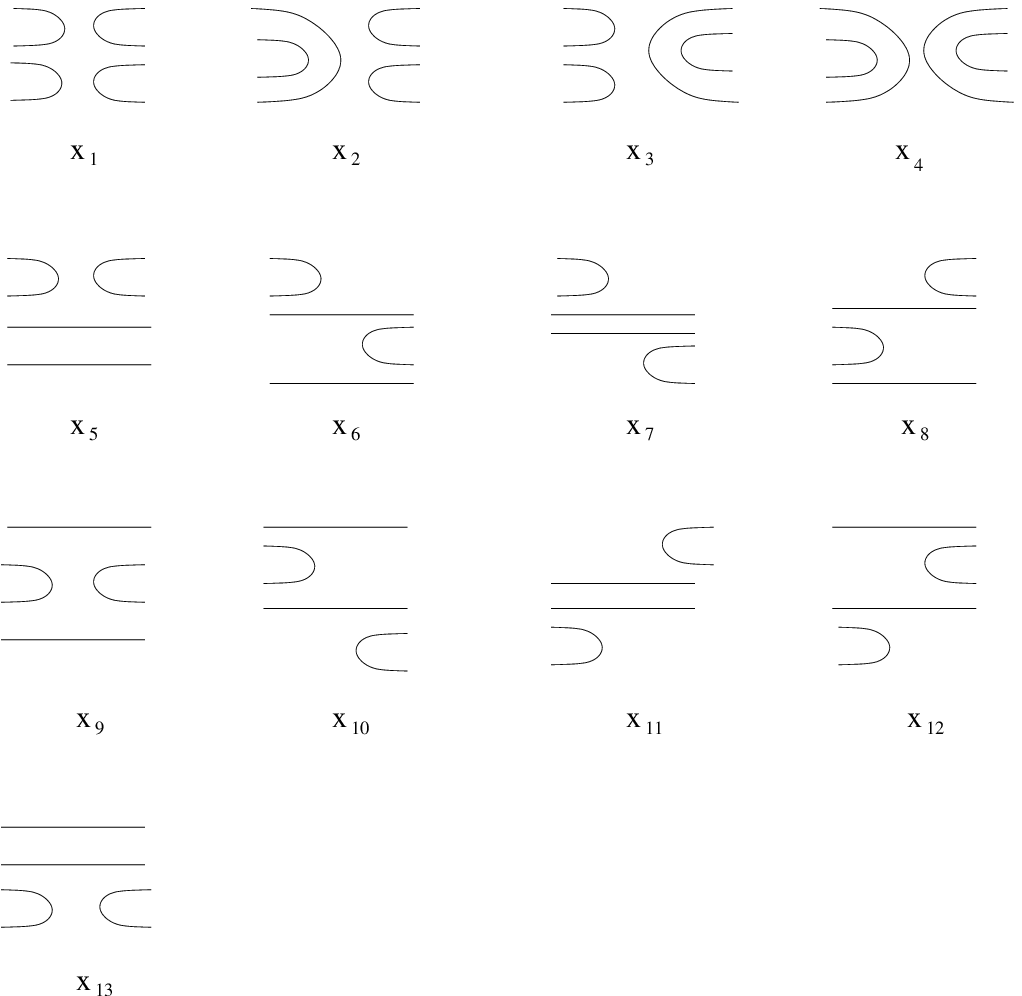}
\caption{}
\end{figure}
\begin{figure}
\centering 
\psfig{file=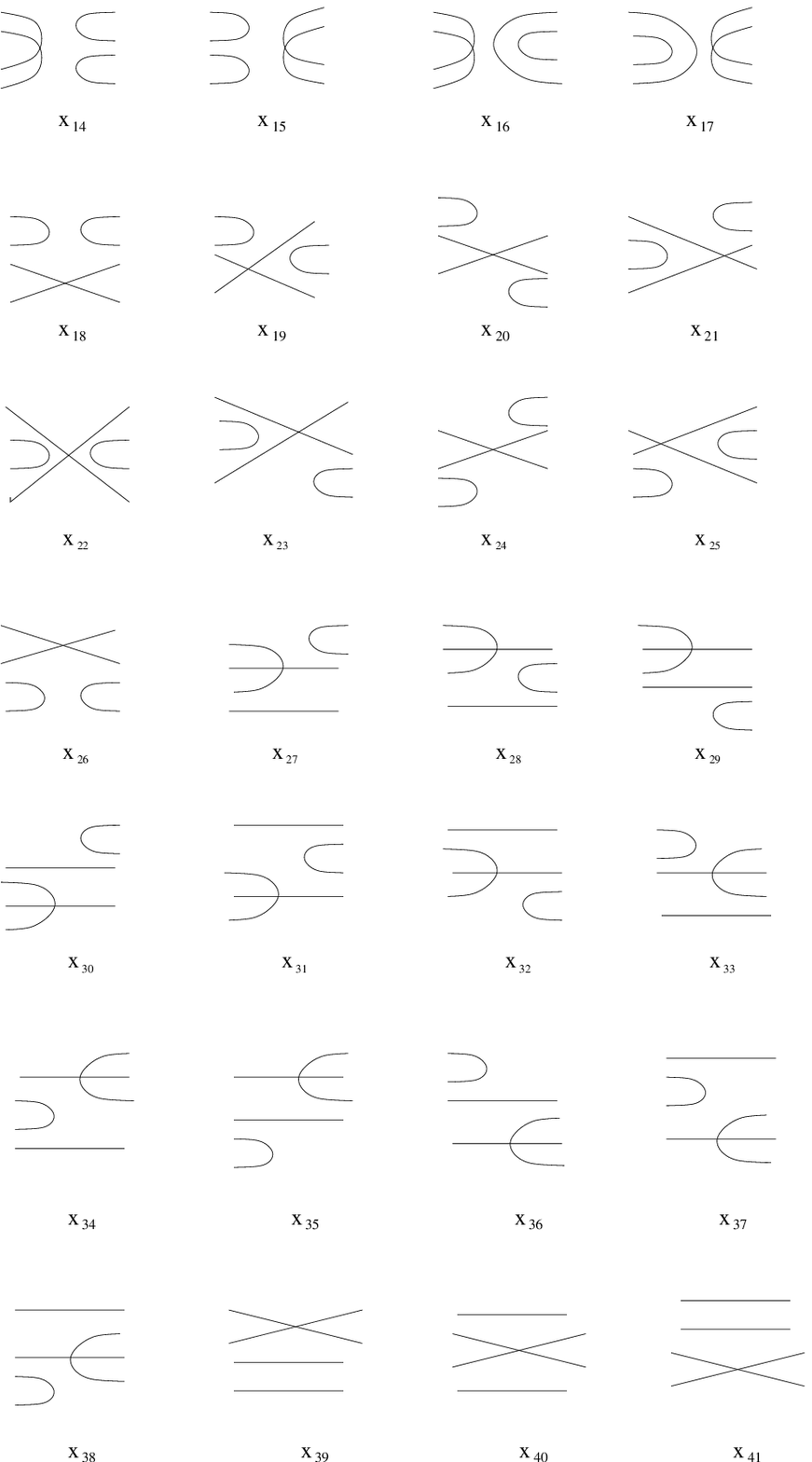}
\caption{}
\end{figure}
\begin{figure}
\centering 
\psfig{file=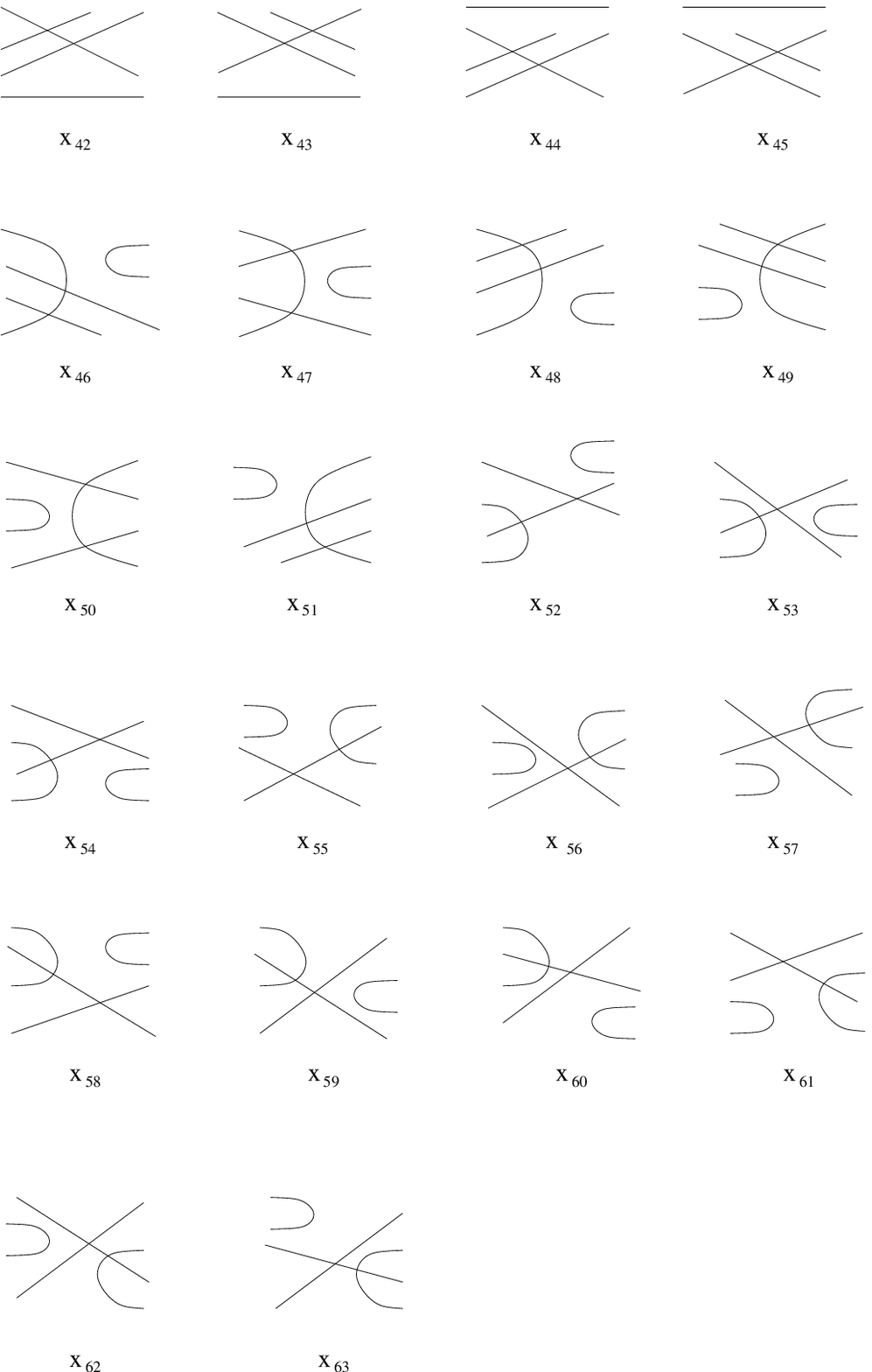}
\caption{}
\end{figure}
 ignored.)
The endpoints of the chord diagrams $x_i$ on the right correspond to the endpoints of the four arrows of the quadruple 
point.
Each of the simplifications $C_i$ or $\bar C_i$ of each of the triple points $p_j$ or $p'_j$ near the quadruple point 
defines an element $t_i(j)$ 
in $T_4$ by calculating its extended Kauffman bracket. Remember, that we use the identifications shown in Fig. 38
\begin{figure}
\centering 
\psfig{file=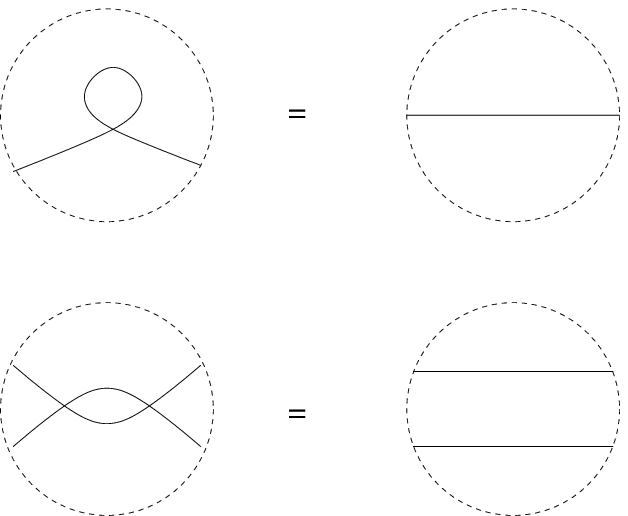}
\caption{}
\end{figure} 
(compare  Definition 9). The tetrahedron equation (c) becomes then the following equation (compare Fig. 14)

$\sum_i(-t_i(1) +t_i(1') -t_i(2) +t_i(2') +t_i(3) -t_i(3') +t_i(4) -t_i(4')) = 0$.        (e)

In this equation the $C_i$ and $\bar C_i$ are the variables. Because the $x_k$ are all independent in $T_4$, the 
equation
(e) splitts into 63 equations, one for each $x_k$.

We give below the contribution of each triple point $p_i$ and $p'_i$ in $T_4$. 
Notice, that we do not consider smoothings of $p_2$ and $p_3$ because they cancel out with the smoothings of
$p'_2$ respectively $p'_3$. 

 $p_1$:

$A^3\bar C_1x_5 +(A\bar C_1 + AC_4+A\bar C_6)x_6 +(A\bar C_1 -A^{-3}\bar C_1 +A^{-1}C_2 +A^{-1}C_4 +A^{-1}C_5)x_1 
+(A^{-1}\bar C_1 -A^{-5}\bar C_1 +A^{-3}C_2 +A^{-3}C_4 +A^{-1}C_4 +A^{-3}C_5+A^{-1}\bar C_6)x_3 +AC_1x_{10} +(A^{-1}C_1+A^{-1}C_3)x_{13} 
+A^3C_2x_8 + (AC_2 +A\bar C_4+AC_3)x_9 +AC_2x_2 +AC_2x_{11} +(A^{-1}C_2 +A^{-1}\bar C_4+A^{-3}C_3+A^{-1}C_3)x_{12} +A^{-1}C_2x_4 
+A^3C_4x_{33} +AC_4x_{20} +AC_4x_{15} +A^{-1}C_4x_{19} +A^3\bar C_4x_{34} +A\bar C_4x_{23}+A\bar C_4x_{35}
+(A^{-3}\bar C_4+A^{-1}C_5+A^{-3}C_6+A^{-1}C_6+A^{-3}\bar C_6)x_{25}+(A^{-1}\bar C_4+A^{-1}\bar C_5+A^{-1}C_6+A^{-1}\bar C_6)x_{26}+A^{-1}\bar C_4x_{22}+A\bar C_5x_{29}
+A^3C_5x_{27}+(AC_5+AC_6)x_{28}+AC_5x_{14}+ AC_5x_{24}+A^{-1}C_5x_{16}+A^3C_3x_{40}+AC_3x_{32}+AC_3x_{38}+
A^{-1}C_3x_{31}+A^3\bar C_3x_{39}+A\bar C_3x_{28}+A\bar C_3x_{26}-A^{-3}\bar C_3x_{26}+A^{-1}\bar C_3x_{25}-A^{-5}\bar C_3x_{25}+
+A^3C_6x_{42}+AC_6x_{48}+AC_6x_{61}+A^{-1}C_6x_{47}+A^3\bar C_6x_{43}+A\bar C_6x_{54}+A\bar C_6x_{49}+A^{-1}\bar C_6x_{53}$

$p'_1$:

$A\bar C_1x_6+A^{-1}\bar C_1x_5+A^3C_1x_{13}+(AC_1-A^{-3}C_1+A^{-1}C_2+A^{-1}\bar C_4+A^{-1}\bar C_5)x_1+(AC_1+A\bar C_5)x_{10}
+(A^{-1}C_1-A^{-5}C_1+A^{-3}C_2+A^{-3}\bar C_4+A^{-3}\bar C_5+A^{-1}\bar C_5)x_2+A^3C_2x_{12}+AC_2x_3+(AC_2+AC_5)x_9
+A^{-1}C_2x_4+(A^{-1}C_2+A^{-1}C_5)x_8+A\bar C_3x_9+(A^{-3}\bar C_3+A^{-1}\bar C_3)x_8+A^{-1}\bar C_3x_5+A\bar C_6x_{10}
+A^{-1}\bar C_6x_2+AC_4x_{36}+(A^{-1}C_4+A^{-1}C_5)x_{18}+A^3\bar C_4x_{38}+A\bar C_4x_{24}+A\bar C_4x_{15}
+A\bar C_4x_{37}+A^{-1}\bar C_4x_{17}+(A^{-1}\bar C_4+A^{-3}C_5)x_{21}+A^3C_5x_{31}+AC_5x_{30}+AC_5x_{19}
+A^{-1}C_5x_{22}+A^3\bar C_5x_{32}+A\bar C_5x_{14}+A\bar C_5x_{20}+A^{-1}\bar C_5x_{23}+A^3C_3x_{41}+AC_3x_{18}
-A^{-3}C_3x_{18}+AC_3x_{37}+A^{-1}C_3x_{21}-A^{-5}C_3x_{21}+A^3\bar C_3x_{40}+A\bar C_3x_{27}+A\bar C_3x_{33}
+A^{-1}\bar C_3x_{34}+A^3C_6x_{44}+AC_6x_{58}+AC_6x_{51}+AC_6x_{37}+A^{-3}C_6x_{21}+A^{-1}C_6x_{50}
+A^{-1}C_6x_{21}+A^{-1}C_6x_{18}+A^3\bar C_6x_{45}+A\bar C_6x_{46}+A\bar C_6x_{55}+A^{-3}\bar C_6x_{21}
+A^{-1}\bar C_6x_{56}+A^{-1}\bar C_6x_{18}$

$p_2$:

$AC_5x_9+(AC_5+A^{-1}C_5)x_{12}+A^{-1}C_5x_3+A\bar C_5x_{10}+(A\bar C_5+A^{-1}\bar C_5)x_{13}+A^{-1}\bar C_5x_1
+AC_6x_{12}+A^{-1}C_6x_3+A\bar C_6x_{10}+A^{-1}\bar C_6x_{13}+AC_3x_{13}+A^{-1}C_3x_1+A\bar C_3x_9
+A^{-1}\bar C_3x_{12}+AC_4x_{36}+A^{-1}C_4x_{15}+A\bar C_4x_{37}+A^{-1}\bar C_4x_{17}+A^3C_5x_{31}
+AC_5x_{19}+A^{-3}C_5x_{25}+A^{-1}C_5x_{22}+A^3\bar C_5x_{32}+A\bar C_5x_{20}+A^{-3}\bar C_5x_{26}
+A^{-1}\bar C_5x_{23}+A^3C_3x_{41}+AC_3x_{18}+AC_3x_{37}+A^{-3}C_3x_{24}+A^{-1}C_3x_{38}+A^{-1}C_3x_{21}
+A^3\bar C_3x_{40}+A\bar C_3x_{33}+A\bar C_3x_{38}+A^{-3}\bar C_3x_{35}+A^{-1}\bar C_3x_{15}+A^{-1}\bar C_3x_{34}
+A^3C_6x_{44}+A^{-3}C_6x_{49}+AC_6x_{51}+AC_6x_{37}+(A^{-1}C_6+A\bar C_6)x_{38}+A^{-1}C_6x_{50}+A^3\bar C_6x_{45}
+A^{-3}\bar C_6x_{57}+A\bar C_6x_{55}+A^{-1}\bar C_6x_{15}+A^{-1}\bar C_6x_{56}$

$p'_2$:

$(AC_4+A^{-1}C_4)x_5+AC_4x_6+A^{-1}C_4x_1+(A\bar C_4+A^{-1}\bar C_4)x_8+A\bar C_4x_9+A^{-1}\bar C_4x_2+AC_3x_9
+A^{-1}C_3x_8+A\bar C_3x_5+A^{-1}\bar C_3x_1+AC_6x_8+A^{-1}C_6x_2+A\bar C_6x_6+A^{-1}\bar C_6x_5+AC_5x_{28}
+A^{-1}C_5x_{16}+A\bar C_5x_{29}+A^{-1}\bar C_5x_{14}+A^3C_4x_{33}+AC_4x_{20}+A^{-3}C_4x_{18}+A^{-1}C_4x_{19}
+A^3\bar C_4x_{34}+A\bar C_4x_{23}+A^{-3}\bar C_4x_{21}+A^{-1}\bar C_4x_{22}+A^3C_3x_{40}+AC_3x_{27}+AC_3x_{32}
+A^{-3}C_3x_{30}+A^{-1}C_3x_{31}+A^{-1}C_3x_{14}+A^3\bar C_3x_{39}+A\bar C_3x_{28}+A\bar C_3x_{26}+A^{-3}\bar C_3x_{24}
+A^{-1}\bar C_3x_{25}+A^{-1}\bar C_3x_{27}+A^3C_6x_{42}+A^{-3}C_6x_{46}+AC_6x_{28}+AC_6x_{48}+A^{-1}C_6x_{47}
+(A^{-1}C_6+A\bar C_6)x_{27}+A^3\bar C_6x_{43}+A^{-3}\bar C_6x_{52}+A\bar C_6x_{54}+A^{-1}\bar C_6x_{53}+A^{-1}\bar C_6x_{14}$

$p_3$:

$AC_4x_9+(AC_4+A^{-1}C_4)x_{10}+A^{-1}C_4x_2+A\bar C_4x{12}+(A\bar C_4+A^{-1}\bar C_4)x_{13}+A^{-1}\bar C_4x_1
+AC_3x_{13}+A^{-1}C_3x_1+A\bar C_3x_9+A^{-1}\bar C_3x_{10}+AC_6x_{12}+A^{-1}C_6x_{13}+A\bar C_6x_{10}
+A^{-1}\bar C_6x_2+AC_5x_{30}+A^{-1}C_5x_{14}+A\bar C_5x_{31}+A^{-1}\bar C_5x_{16}+A^3C_4x_{37}+AC_4x_{21}
+A^{-3}C_4x_{23}+A^{-1}C_4x_{22}+A^3\bar C_4x_{38}+A\bar C_4x_{24}+A^{-3}\bar C_4x_{26}+A^{-1}\bar C_4x_{25}
+A^3C_3x_{41}+AC_3x_{31}+AC_3x_{18}+A^{-3}C_3x_{20}+A^{-1}C_3x_{19}+A^{-1}C_3x_{32}+A^3\bar C_3x_{40}
+A\bar C_3x_{32}+A\bar C_3x_{27}+A^{-3}\bar C_3x_{29}+A^{-1}\bar C_3x_{28}+A^{-1}\bar C_3x_{14}+A^3C_6x_{44}
+AC_6x_{32}+AC_6x_{58}+A^{-3}C_6x_{60}+A^{-1}C_6x_{59}+A^{-1}C_6x_{14}+A^3\bar C_6x_{45}+A\bar C_6x_{31}
+A\bar C_6x_{46}+A^{-3}\bar C_6x_{48}+A^{-1}\bar C_6x_{47}+A^{-1}\bar C_6x_{32}$

$p'_3$:

$(AC_5+A^{-1}C_5)x_5+AC_5x_8+A^{-1}C_5x_1+(A\bar C_5+A^{-1}\bar C_5)x_6+A\bar C_5x_9+A^{-1}\bar C_5x_3+AC_3x_9
+A^{-1}C_3x_6+A\bar C_3x_5+A^{-1}\bar C_3x_1+AC_6x_8+A^{-1}C_6x_5+A\bar C_6x_6+A^{-1}\bar C_6x_3+AC_4x_{34}
+A^{-1}C_4x_{17}+A\bar C_4x_{35}+A^{-1}\bar C_4x{15}+A^3C_5x_{27}+AC_5x_{24}+A^{-3}C_5x_{18}+A^{-1}C_5x_{21}
+A^3\bar C_5x_{28}+A\bar C_5x{25}+A^{-3}\bar C_5x_{19}+A^{-1}\bar C_5x_{22}+A^3C_3x_{40}+AC_3x_{38}
+AC_3x_{33}+A^{-3}C_3x_{36}+A^{-1}C_3x_{15}+A^{-1}C_3x_{37}+A^3\bar C_3x_{39}+A\bar C_3x_{26}+A\bar C_3x_{34}
+A^{-3}\bar C_3x_{20}+A^{-1}\bar C_3x_{33}+A^{-1}\bar C_3x_{23}+A^3C_6x_{42}+AC_6x_{61}+AC_6x_{33}
+A^{-3}C_6x_{63}+A^{-1}C_6x_{15}+A^{-1}C_6x_{62}+A^3\bar C_6x_{43}+A\bar C_6x_{49}+A\bar C_6x_{34}
+A^{-3}\bar C_6x_{51}+A^{-1}\bar C_6x_{33}+A^{-1}\bar C_6x_{50}$

$p_4$:

$A^3\bar C_1x_5+A^3\bar C_2x_6+(A\bar C_1-A^{-3}\bar C_1+A^{-1}\bar C_2)x_1+A\bar C_2x_3+A^{-1}\bar C_2x_{10}
+A^{-1}C_1x_{13}+A\bar C_1x_8+A\bar C_2x_9+(A^{-1}\bar C_1-A^{-5}\bar C_1+A^{-3}\bar C_2)x_2+AC_1x_{12}
+A^{-1}\bar C_2x_4+A\bar C_2x_7+AC_3x_9+A^{-3}C_3x_{10}+A^{-1}C_3x_{10}+A^{-1}C_3x_{13}+AC_6x_8+A^{-1}C_6x_2
+A\bar C_4x_{35}+(A^{-1}\bar C_4+A^{-1}\bar C_5)x_{26}+A^3C_4x_{33}+AC_4x_{20}+AC_4x_{15}+AC_4x_{34}
+A^{-1}C_4x_{17}+(A^{-1}C_4+A^{-3}\bar C_5)x_{23}+A^3C_5x_{27}+AC_5x_{14}+AC_5x_{24}+A^{-1}C_5x_{21}
+A^3\bar C_5x_{28}+A\bar C_5x_{29}+A\bar C_5x_{25}+A^{-1}\bar C_5x_{22}+A^3C_3x_{40}+AC_3x_{32}+AC_3x_{38}
+A^{-1}C_3x_{37}+A^3\bar C_3x_{39}+A\bar C_3x_{26}-A^{-3}\bar C_3x_{26}+A\bar C_3x_{34}+A^{-1}\bar C_3x_{23}
-A^{-5}\bar C_3x_{23}+A^3C_6x_{42}+AC_6x_{48}+AC_6x_{61}+A^{-3}C_6x_{23}+A^{-1}C_6x_{62}+A^{-1}C_6x_{26}
+A^3\bar C_6x_{43}+A\bar C_6x_{54}+A\bar C_6x_{49}+A\bar C_6x_{34}+A^{-3}\bar C_6x_{23}+A^{-1}\bar C_6x_{50}
+A^{-1}\bar C_6x_{23}+A^{-1}\bar C_6x_{26}$

$p'_4$:

$A^{-1}\bar C_1x_5+A^{-1}\bar C_2x_6+(AC_1-A^{-3}C_1+A^{-1}\bar C_2)x_1+(A^{-1}C_1-A^{-5}C_1+A^{-3}\bar C_2)x_3
+A^3\bar C_2x_{10}+A^3C_1x_{13}+A\bar C_1x_8+A\bar C_2x_9+A\bar C_2x_2+AC_1x_{12}+A^{-1}\bar C_2x_4
+A\bar C_2x_7+A^{-3}\bar C_3x_6+A\bar C_3x_9+A^{-1}\bar C_3x_5+A^{-1}\bar C_3x_6+AC_6x_{12}+A^{-1}C_6x_3
+AC_5x_{30}+(A^{-1}C_5+A^{-1}C_4)x_{18}+A^3\bar C_5x_{32}+A\bar C_5x_{31}+A\bar C_5x_{14}+A\bar C_5x_{20}
+(A^{-1}\bar C_5+A^{-3}C_4)x_{19}+A^{-1}\bar C_5x_{16}+A^3C_4x_{37}+AC_4x_{21}+AC_4x_{36}+A^{-1}C_4x_{22}
+A^3\bar C_4x_{38}+A\bar C_4x_{24}+A\bar C_4x_{15}+A^{-1}\bar C_4x_{25}+A^3C_3x_{41}+AC_3x_{31}
+AC_3x_{18}-A^{-3}C_3x_{18}+A^{-1}C_3x_{19}-A^{-5}C_3x_{19}+A^3\bar C_3x_{40}+A\bar C_3x_{27}+A\bar C_3x_{33}
+A^{-1}\bar C_3x_{28}+A^3C_6x_{44}+AC_6x_{58}+AC_6x_{51}+A^{-3}C_6x_{19}+A^{-1}C_6x_{18}+A^{-1}C_6x_{59}
+A^3\bar C_6x_{45}+A\bar C_6x_{31}+A\bar C_6x_{46}+A\bar C_6x_{55}+A^{-3}\bar C_6x_{19}+A^{-1}\bar C_6x_{18}
+A^{-1}\bar C_6x_{19}+A^{-1}\bar C_6x_{47}$

We obtain the following system (*) of 39 equations (we drop the trivial equations), where the expression after each $x_i$
is equal to $0$:

$x_1$:   $-A^{-1}\bar C_5-A^{-1}C_5+A^{-1}\bar C_4+A^{-1}C_4$

$x_2$:   $(A^{-1}-A^{-5})C_1+(A^{-1}-A^{-5})\bar C_1+(A^{-3}-A)C_2+(A^{-3}-A)\bar C_2+(A^{-3}+A^{-1})C_4+(A^{-3}+A^{-1})\bar C_4
+(A^{-3}+^{-1})C_5+(A^{-3}+A^{-1})\bar C_5+2A^{-1}C_6+2A^{-1}\bar C_6$

$x_3$:    $-[(A^{-1}-A^{-5})C_1+(A^{-1}-A^{-5})\bar C_1+(A^{-3}-A)C_2+(A^{-3}-A)\bar C_2+(A^{-3}+A^{-1})C_4+(A^{-3}+A^{-1})\bar C_4
+(A^{-3}+^{-1})C_5+(A^{-3}+A^{-1})\bar C_5+2A^{-1}C_6+2A^{-1}\bar C_6]$

$x_5$:  $(A+A^{-1})C_4-(A+A^{-1})C_5-A^{-1}C_6+A^{-1}\bar C_6$

$x_6$:  $(A^3-A^{-1})\bar C_2-A^{-1}C_3-(A^{-1}+A^{-3})\bar C_3-A^{-1}C_4-(A+A^{-1})\bar C_5-A\bar C_6$

$x_8$:  $(A^{-1}-A^3)C_2+A^{-1}C_3+(A^{-1}+A^{-3})\bar C_3+(A+A^{-1})\bar C_4+A^{-1}C_5+AC_6$

$x_{10}$: $(A^{-1}-A^3)\bar C_2+(A^{-1}+A^{-3})C_3+A^{-1}\bar C_3+(A+A^{-1})C_4+A^{-1}\bar C_5+A\bar C_6$

$x_{12}$: $(A^3-A^{-1})C_2-(A^{-1}+A^{-3})C_3-A^{-1}\bar C_3-A^{-1}\bar C_4-(A+A^{-1})C_5-AC_6$

$x_{13}$: $(A+A^{-1})\bar C_4-(A+A^{-1})\bar C_5+A^{-1}C_6-A^{-1}\bar C_6$

$x_{14}$: $A^{-1}C_3+A^{-1}\bar C_3+A^{-1}C_5+A^{-1}\bar C_5+A^{-1}C_6+A^{-1}\bar C_6$

$x_{15}$: $-A^{-1}C_3-A^{-1}\bar C_3-A^{-1}C_4-A^{-1}\bar C_4-A^{-1}C_6-A^{-1}\bar C_6$

$x_{18}$: $A^{-3}C_4-A^{-3}C_5$

$x_{19}$: $A^{-5}C_3-A^{-3}C_4-(A^{-1}+A^{-3})\bar C_5-A^{-3}C_6-(A^{-1}+A^{-3})\bar C_6$

$x_{20}$: $A^{-3}C_3-A^{-3}\bar C_3+AC_4-A\bar C_5$

$x_{21}$: $-A^{-5}C_3+(A^{-1}+A^{-3})\bar C_4+A^{-3}C_5+(A^{-1}+A^{-3})C_6+A^{-3}\bar C_6$

$x_{23}$: $-A^{-5}\bar C_3+(A^{-1}+A^{-3})C_4+A^{-3}\bar C_5+A^{-3}C_6+(A^{-1}+A^{-3})\bar C_6$

$x_{24}$: $-A^{-3}C_3+A^{-3}\bar C_3+A\bar C_4-AC_5$

$x_{25}$: $A^{-5}\bar C_3-A^{-3}\bar C_4-(A^{-1}+A^{-3})C_5-(A^{-1}+A^{-3})C_6-A^{-3}\bar C_6$

$x_{26}$: $A^{-3}\bar C_4-A^{-3}\bar C_5$

$x_{27}$: $AC_3+(A+A^{-1})\bar C_3-A^3C_5+A^{-1}C_6+A\bar C_6$

$x_{29}$: $A^{-3}\bar C_3+A\bar C_5$

$x_{30}$: $A^{-3}C_3+AC_5$

$x_{32}$: $(A+A^{-1})C_3+A\bar C_3-A^3\bar C_5+AC_6+A^{-1}\bar C_6$

$x_{33}$: $-AC_3-(A+A^{-1})\bar C_3+A^3C_4-AC_6-A^{-1}\bar C_6$

$x_{35}$: $-A^{-3}\bar C_3-A\bar C_4$

$x_{36}$: $-A^{-3}C_3-AC_4$

$x_{38}$: $-(A+A^{-1})C_3-A\bar C_3+A^3\bar C_4-A^{-1}C_6-A\bar C_6$

$x_{46}$: $A^{-3}C_6+A\bar C_6$

$x_{48}$: $AC_6+A^{-3}\bar C_6$

$x_{49}$: $-A^{-3}C_6-A\bar C_6$

$x_{51}$: $-AC_6-A^{-3}\bar C_6$

$x_{52}$: $A^{-3}\bar C_6$

$x_{54}$: $A\bar C_6$

$x_{55}$: $-A\bar C_6$

$x_{57}$: $-A^{-3}\bar C_6$

$x_{58}$: $AC_6$

$x_{60}$: $A^{-3}C_6$

$x_{61}$: $-AC_6$

$x_{63}$: $-A^{-3}C_6$

There is a good control of eventual errors, because each equation appears several times.

One easily calculates that this system of equations has exactly two non trivial solutions.

\begin{definition}
{\em Solution 1} :

$C_1 = -\bar C_1$ and all other are 0.

\end{definition}

\begin{lemma}
Let $S$ be the sum from Definition 11 and let $s$ be a meridian of $\Sigma_{q}^{(2)}$ at a stratum which corresponds
to a positive quadruple point. Then $S(s) = 0$.
\end{lemma}
{\em Proof:\/} It follows from Remarks 9 and 15 that all the triple points which contribute non trivially to $S(s)$
have the {\em same} distinguished crossing. Moreover, $n(d(p))$ and $[d(p)]$ do not depend neither of the smoothing
nor of the extended Kauffmann state. The weight function doesn't matter, because all eight triple points are of the
same type. The lemma follows now from the fact that $C_1 = -\bar C_1$ is a
solution of the system (*).
$\Box$

\begin{definition}
{\em Solution 2} :

$C_3 = -\bar C_3$, $C_4 = -\bar C_4$, $C_5 = -\bar C_5$, $C_4 = C_5$, $C_3 = -A^4C_4$ and all other are 0.

\end{definition}

\begin{remark}
The Solution 1 uses only smoothings and Solution 2 uses only simplifications with exactly one double point.
\end{remark}

Notice, that we can no longer multiply the extended Kaufmann bracket by $t^{[d(p)]}$ in the case of Solution 2, because 
all eight triple points from a quadruple point contribute now non trivially to (e).

Let $H_1(M_0)^{rot_\pi}$ be the equivariant
first homology group of $M_0$ with respect to the involution $rot_\pi$. The canonical loop induces a well defined
 homology class $[rot]$ modulo sliding classes in $H_1(M_0)^{rot_\pi}$ as well as in 
$H_1(M_0/rot_\pi,\mathbb{Z}/2\mathbb{Z})$.

 The involution $rot_\pi$ interchanges 
$C_i(p)$ with $\bar C_i(\bar p)$ and changes the sign: $sign(p) = -sign(\bar p)$. Therefore, if we would set $C_i = \bar C_i$, then 
the contributions of these two simplifications would always cancel out in the canonical loop.

\begin{remark}
Instead of the tetrahedron equation we could consider the {\em equivariant tetrahedron equation}: let $s$ be an oriented
meridian of $\Sigma_{q}^{(2)}$ at a stratum which corresponds
to a positive quadruple point. Then the equation (a) could be replaced by 

$L(s) + L(rot_{\pi}(s)) = 0$    (f).

(Here, the loop $rot_{\pi}(s)$ has of course the induced orientation.)

However, it turns out that the equivariant tetrahedron equation has no additional solutions. We left the verification to the reader.
\end{remark}

\begin{lemma}
Let $X$ be the sum introduced in Definitions 17, 18 and let $s$ be a meridian of $\Sigma_{q}^{(2)}$ at a stratum which 
corresponds
to a positive quadruple point. Then $X(s)  = 0$.
\end{lemma}
{\em Proof:\/}
Contractible circles (or d-circles) are traded to factors $-A^2 -A^{-2}$ in the bracket, but the configuration in the annulus 
of the non contractible circles (or h-circles) is an invariant (compare \cite{APS} and \cite{AF}). A double circle 
(compare Definition 9) splitts the annulus into three regions. One of these regions is contractible and hence, can not 
contain h-circles. The configuration of the h-circles in the remaining two regions is determined by their numbers in 
these regions. No h-circle can slide over the double circle, because each bracket is defined as a relative knot invariant
of the diagram with respect to the triple point. 
The lemma follows now from Solution 2 and the definition of the weight 
function $f_i(1)$ for $i \in \{3, 4, 5\}$.
$\Box$

\subsection{Solutions of the positive tetrahedron equation using the extended Kauffman state model 
for the Alexander polynomial}

We proced as in the previous subsection but with the simplifications replaced by the markings. The calculations take 105 
pages in the handwritten manuscript. We therefore decided not to include them into the paper, but we send copies of the calculations 
by request.

Let $T_4$ be the free 
$\mathbb{Z}[A,A^{-1}]$-modul generated by the 56 {\em markings} with exactly three dots  of the oriented chord diagram
 $X_4$ shown in Fig.39. We number the eight regions as shown in Fig. 39 too. We give now names to the generators 
by using the lexicographical order.
\begin{figure}
\centering 
\psfig{file=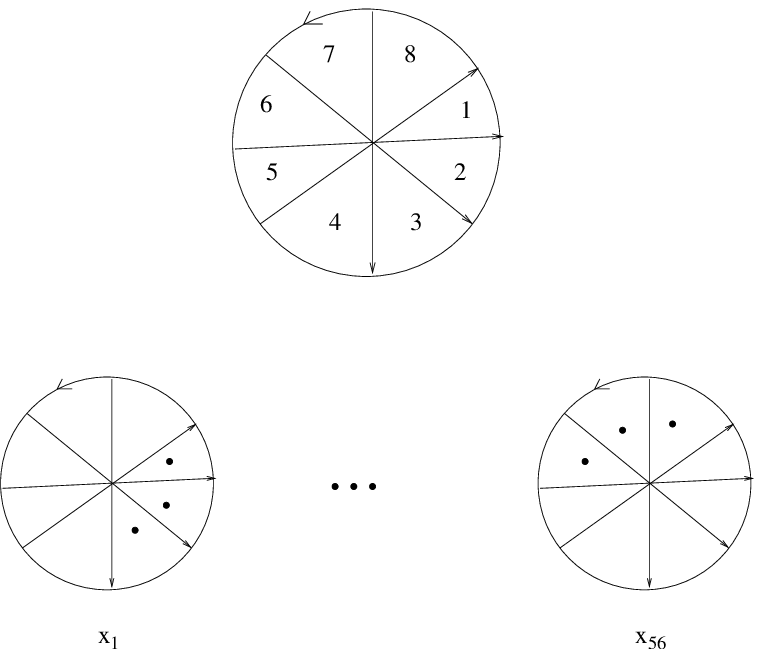}
\caption{}
\end{figure}
Each of the markings $C_i$ or $\bar C_i$ of each of the triple points $p_j$ or $p'_j$ near the quadruple point 
defines an element $t_i(j)$ 
in $T_4$ by calculating its extended Kauffman bracket. We obtain the analogue tetrahedron 
equation (e). This time, (e) splitts into a system (**) of 48 equations (we drop trivial equations):

$x_2$: $-A^2C_1-A^3\bar C_2-A\bar C_4-C_5$

$x_3$: $2AC_1+\bar C_2+A^2\bar C_2+\bar C_4+A^2\bar C_4+(A+A^{-1})C_5$

$x_4$: $-(1+A^2)C_1-(A+A^{-1})\bar C_2-(A+A^3)\bar C_4-(1+A^2)C_5$

$x_5$: $AC_1+\bar C_2+A^2\bar C_4+AC_5$

$x_7$: $-A\bar C_1-C_2-A^2C_4-A\bar C_5$

$x_8$: $(1+A^2)\bar C_1+(A+A^{-1})C_2+(A+A^3)C_4+(1+A^2)\bar C_5$

$x_9$: $-AC_1-A\bar C_1-C_2-\bar C_2-A^2C_4-A^2\bar C_4-AC_5-A\bar C_5$

$x_{10}$: $(1+A^2)C_1+(A+A^{-1})\bar C_2+(A+A^3)\bar C_4+(1+A^2)C_5$

$x_{11}$: $-AC_1-\bar C_2-A^2\bar C_4-AC_5$

$x_{12}$: $(A-A^{-1})C_1-(1+A^{-2})\bar C_2+(1+A^2)\bar C_3+(A+A^{-1})\bar C_5+(A+A^3)\bar C_6$

$x_{13}$: $-C_1+A^{-1}\bar C_2+A^{-1}C_3-A\bar C_3+A^{-1}C_4-\bar C_5-(1+A^2)\bar C_6$

$x_{14}$: $AC_1+A^{-1}\bar C_1+\bar C_2-\bar C_3+(A^2-A^{-2})\bar C_4+AC_5+A^{-1}C_6$

$x_{15}$: $-C_1+A^{-1}C_2-A\bar C_2-AC_3-A\bar C_4-(1-A^{-2})C_5-C_6$

$x_{16}$: $2A^{-1}C_1-(1+A^{-2})C_3-(1+A^{-2})C_4+(A+A^{-1})\bar C_6$

$x_{17}$: $-(1+A^{-2})\bar C_1+(A+A^{-1})\bar C_3+(A^{-1}+A^{-3})\bar C_4-(1+A^{-2})C_6$

$x_{18}$: $(A^{-1}-A)C_1-(1+A^{-2})C_2+(1+A^2)C_3+A^{-1}C_5+AC_6$

$x_{19}$: $A^{-1}\bar C_1-A^{-2}\bar C_4+A^{-1}C_6-\bar C_3$

$x_{20}$: $A^{-1}C_2-AC_3+(A^{-1}-A)\bar C_4-C_5-C_6$

$x_{23}$: $A\bar C_1+C_2+A^2C_4+A\bar C_5$

$x_{24}$: $-(1+A^2)\bar C_1-(A+A^{-1})C_2-(A+A^3)C_4-(1+A^2)\bar C_5$

$x_{25}$: $2A\bar C_1+(1+A^2)C_2+(1+A^2)C_4+(A+A^{-1})\bar C_5$

$x_{26}$: $-A^2\bar C_1-A^3C_2-AC_4-\bar C_5$

$x_{27}$: $-A^{-1}\bar C_2+A\bar C_3+(A-A^{-1})C_4+\bar C_5+\bar C_6$

$x_{28}$: $(1+A^{-2})\bar C_2-(1+A^2)\bar C_3+(A^{-2}-A^2)C_4-(A+A^{-1})\bar C_5-(A+A^{-1})\bar C_6$

$x_{29}$: $\bar C_1+AC_2-A^{-1}\bar C_2+A\bar C_3+AC_4+(1+A^{-2})\bar C_5+\bar C_6$

$x_{30}$: $-AC_1-A\bar C_1-A^2C_2-A^2\bar C_2-C_4-\bar C_4-A^{-1}C_5-A^{-1}\bar C_5$

$x_{31}$: $(1+A^{-2})C_1-(A+A^{-1})C_3-(A^{-1}+A^{-3})C_4+(1+A^{-2})\bar C_6$

$x_{32}$: $-A^{-1}C_1-A^{-1}\bar C_1+C_3+\bar C_3+A^{-2}C_4+A^{-2}\bar C_4-A^{-1}C_6-A^{-1}\bar C_6$

$x_{33}$: $C_1-A^{-1}C_2+A\bar C_2+AC_3+A\bar C_4+(1+A^{-2})C_5+C_6$

$x_{34}$: $(1+A^{-2})\bar C_1-(A+A^{-1})\bar C_3-(A^{-1}+A^{-3})\bar C_4+(1+A^{-2})C_6$

$x_{35}$: $(1+A^{-2})C_2-(1+A^2)C_3+(A^{-2}-A^2)\bar C_4-(A+A^{-1})C_5-(A+A^{-1})C_6$

$x_{36}$: $-A^{-1}C_2+AC_3+(A-A^{-1})\bar C_4+C_5+C_6$

$x_{38}$: $A^{-1}\bar C_2-A\bar C_3+(A^{-1}-A)C_4-\bar C_5-\bar C_6$

$x_{39}$: $(A^{-1}-A)\bar C_1-(1+A^{-2})\bar C_2+(1+A^2)\bar C_3+(A^{-1}+A^{-3})\bar C_5+(A+A^{-1})\bar C_6$

$x_{40}$: $-\bar C_1-AC_2+A^{-1}\bar C_2-A\bar C_3-AC_4-(1+A^{-2})\bar C_5-\bar C_6$

$x_{41}$: $A^{-1}C_1-C_3-A^{-2}C_4+A^{-1}\bar C_6$

$x_{42}$: $-(1+A^{-2})C_1+(A+A^{-1})C_3+(A^{-1}+A^{-3})C_4-(1+A^{-2})\bar C_6$

$x_{43}$: $A^{-1}C_1+A\bar C_1+C_2-C_3+(A^2-A^{-2})C_4+A\bar C_5+A^{-1}\bar C_6$

$x_{44}$: $2A^{-1}\bar C_1-(1+A^{-2})\bar C_3-(1+A^{-2})\bar C_4+(A+A^{-1})C_6$

$x_{45}$: $-\bar C_1+A^{-1}C_2-AC_3+A^{-1}\bar C_3+A^{-1}\bar C_4-C_5-(1+A^2)C_6$

$x_{46}$: $(A-A^{-1})\bar C_1-(1+A^{-2})C_2+(1+A^2)C_3+(A+A^{-1})C_5+(A+A^3)C_6$

$x_{48}$: $-A^{-1}C_1+C_3+A^{-2}C_4-A^{-1}\bar C_6$

$x_{49}$: $C_1-A^{-1}\bar C_2-A^{-1}C_3+A\bar C_3-A^{-1}C_4+\bar C_5+(1+A^2)\bar C_6$

$x_{50}$: $A^{-2}\bar C_1-A^{-3}\bar C_3-A^{-1}\bar C_4+C_6$

$x_{51}$: $-A^{-1}C_1-A^{-1}\bar C_1+A^{-2}C_3+A^{-2}\bar C_3+C_4+\bar C_4-AC_6-A\bar C_6$

$x_{52}$: $\bar C_1-A^{-1}C_2+AC_3-A^{-1}\bar C_3-A^{-1}\bar C_4+C_5+(1+A^2)C_6$

$x_{54}$: $A^{-2}C_1-A^{-3}C_3-A^{-1}C_4+\bar C_6$

$x_{55}$: $-A^{-1}\bar C_1+\bar C_3+A^{-2}\bar C_4-A^{-1}C_6$

Again, there is a good control about eventual errors, because the system is over determined and many equations 
appear several times.

One easily calculates that the system (**) has exactly one solution.

\begin{definition}
Solution:

$C_i = -\bar C_i$ for $i \in \{1, 2, 3, 4\}$, $C_2 = C_3 = C_4$, $C_1 = (A + A^{-1})C_2$ and all other are 0.
\end{definition}

\begin{lemma}
Let $\Phi$ be the sum introduced in Definition 20 and let $s$ be a meridian of $\Sigma_{q}^{(2)}$ at a stratum 
which  corresponds
to a positive quadruple point. Then $\Phi(s) = 0$.
\end{lemma}
{\em Proof:\/} The lemma follows immediately from the definition of  $\Phi$ 
and the solution of (**).
$\Box$

Again, the equivariant tetrahedron equation (f) has no new solution.

\subsection{$S$ as solution of the cube equations}

Taking the mirror image of all diagrams correspond to interchanging the cube equations by the "antipodal map"
of the cube. It is clear that the solutions of the cube equations are invariant under this operation. Consequently,
it suffices to consider the relations from six double edges and from three 2-faces which have a commun vertex.

Let $T_3$ be the free $\mathbb{Z}[A,A^{-1}]$-modul generated by all chord diagrams with exactly three chords in the 
disc with distinguished end points and which have no more than one
 double point (and such that this number is minimal).
The generators of $T_3$ correspond exactly to the simplifications $C_1$ up to $C_6$ and $\bar C_1$ up to $\bar C_6$
(see Fig. 15). We call the corresponding generators $x_i$ respectively $\bar x_i$. (The generator $x_7$ never occurs 
and can be ignored.)
 Let $s$ be a meridian of $\Sigma_{a-t}^{(2)}$. Each of the 
simplifications $C_i$ or $\bar C_i$ of each of the 
two triple points $p_j$ or $p'_j$ near the point in $\Sigma_{a-t}^{(2)}$
defines an element $t_i(j)$  in $T_3$ by calculating its extended Kauffman bracket. The edges of the cube relate the 
variables for different types of triple points.

\begin{lemma}
The relations from the edges of the cube are given in Fig. 40 (which means that $C_1(7) = -A^2C_1(1)$ ,
$C_1(2) = -A^{-2}C_1(7)$ and so on). This implies that $C_1(j) = -\bar C_1(j)$ for all $j$.  Moreover, double edges
do not impose any new relation.
\begin{figure}
\centering 
\psfig{file=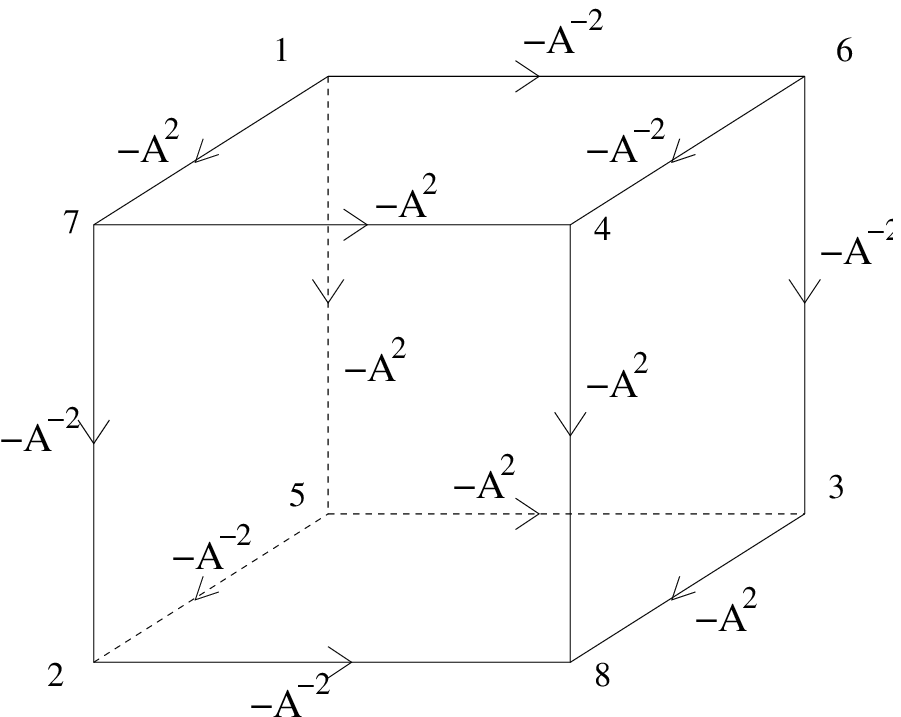}
\caption{}
\end{figure}
\end{lemma}
{\em Proof:\/} The generic diagrams near a stratum of $\Sigma_{a-t}^{(2)}$ which corresponds to e.g.
the edge 1-1-5 of the cube were shown in Fig. 10. The condition $L(s) = 0$ gives us the following system of equations:

$x_1$:  $A^{-1}C_1=-A^{-3}C_1(5)$

$\bar x_1$:  $-A^3\bar C_1=A\bar C_1(5)$

$x_2$:  $AC_1=A^{-1}\bar C_1(5)$

In Fig. 41 we show the generic diagrams near a stratum of $\Sigma_{a-t}^{(2)}$ which corresponds to the edge 5-2-1
of the cube. We obtain the equations:
\begin{figure}
\centering 
\psfig{file=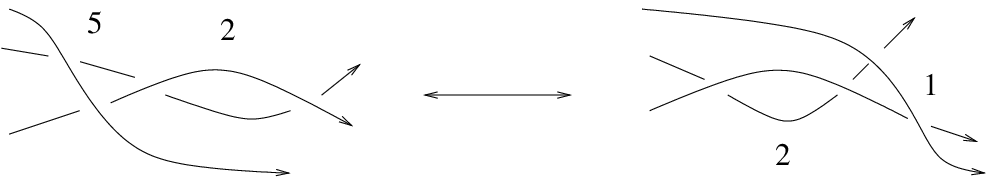}
\caption{}
\end{figure}

$x_1$:  $AC_1(5)=-A^3C_1$

$\bar x_1$:  $-A^{-3}\bar C_1(5)=A^{-1}\bar C_1$

$x_2$:  $A^{-1}C_1(5)=A\bar C_1$

The equations for the two edges are different, but they have the same solution.

The calculations for all other edges are completely analogous. We left the verification to the reader.
$\Box$

Let us consider now the cube equations which come from the 2-faces of the cube. It follows directly from Lemma 16
that the 2-faces corresponding to the relatives of the positive triple point (i.e. 1-5-2-7) and to the relatives of the 
negative triple point (i.e. 8-3-6-4) do not impose new relations. Indeed, e.g. the first 2-face gives:

$C_1(1) = -A^2C_1(5) = C_1(2) = -A^2C_1(7) = C_1(1)$.

However, one easily sees from Fig. 40 that all other 2-faces would imply the relation:  $A^8 = 1$.

We will eleminate this relation by introducing a correction term containing $n(d(p))$ in the definition of $S$.

\begin{lemma}
Let $p(i)$ and $p(j)$ be two vertices in the cube which are connected by an edge. Then $n(d(p(i))) = n(d(p(j)))$
besides for the edges 1-6 and 2-8. In the cases 1-3-6 and 1-4-6 we have $n(d(6)) = n(d(1)) + 1$.
In the cases 2-3-8 and 2-4-8 we have $n(d(2)) = n(d(8)) - 1$. 
\end{lemma}
{\em Proof:\/} This is a case by case verification. The cases e.g. 2-3-8 and 2-4-8 are shown in Fig. 42.
\begin{figure}
\centering 
\psfig{file=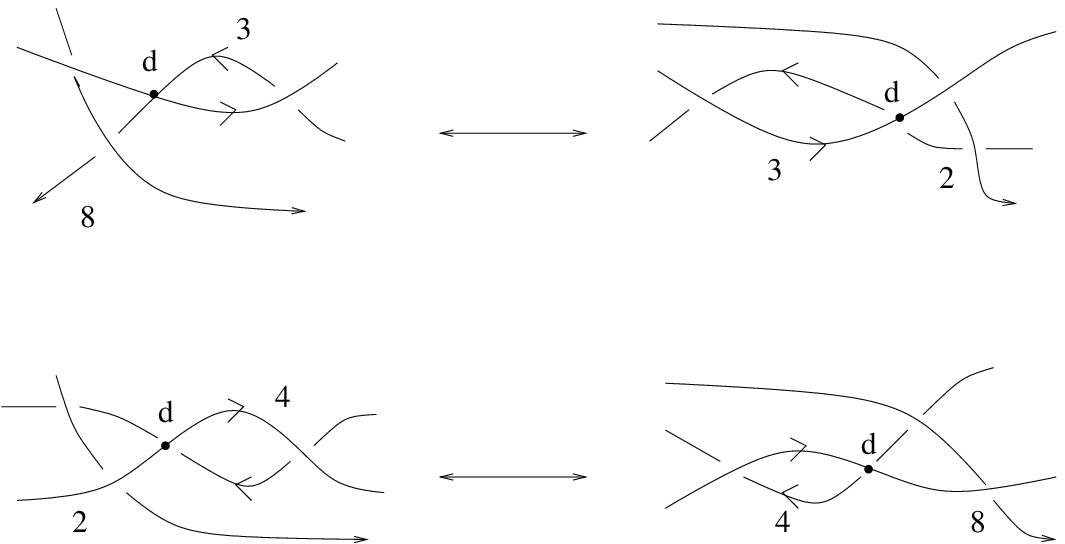}
\caption{}
\end{figure}
The lemma follows easily from the figures.
$\Box$

The above lemma has the following surprising corollary.

Let  $C = S^1 \times I$ be a cylindre in $M$ such that each generic interval $\{x\} \times Int(I)$ intersects 
$\Sigma_t^{(1)}$ transversally in exactly one point and the boundary $\delta C$ does not intersect $\Sigma_t^{(1)}$ at 
all. We call such a $C$ a {\em meridional cylindre of $\Sigma_t^{(1)}$} in $M$.

\begin{corollary}
There is no meridional cylindre $C$ of $\Sigma_t^{(1)}$ in $M$ such that
 $C$ intersects $\Sigma_{a-t}^{(2)}$ transversally in exactly four
points, such that their types  correspond to four vertices of the cube which span a 2-face and which are not all in the 
same family, e.g. 1-6, 6-4, 4-7, 7-1, or 2-8, 8-4, 4-7, 7-2 and so on.
\end{corollary}
{\em Proof:\/} To each generic interval $\{x\} \times I$ we associate $n(d(p)$, where $d(p)$ is the distinguished crossing
of the unique intersection of $\{x\} \times I$ with $\Sigma_t^{(1)}$. Going once along the circle $S^1 \times \{0\}$ in the cylindre we 
come back to the same distinguished crossing, but it follows from the above lemma that $n(d(p))$ would have changed 
by $+1$ or $-1$.
$\Box$

\begin{lemma}
$A^{-8n(d(p))}C_1(1)f_S(j(p))$ (compare Definition 10) is well defined for each triple crossing $p$ and verifies all cube 
equations.
\end{lemma}
{\em Proof:\/}
If we go along the boundary of the 2-faces 1-5-2-7 and 8-3-6-4 then $n(d(p))$ stays constant.
If we go along the boundary of the 2-face 1-6-3-5 then $n(d(p))$ increases by 1 but $C_1(1)$ has to be multiplied by 
$A^{-8}$. The same is true for 1-6-4-7. If we go along 8-2-7-4 then $n(d(p))$ decreases by 1 but $C_1(8)$ has to be 
multiplied by $A^8$. The same is true for 8-2-5-3. It follows now from Lemmas 16, 17 and Definition 10 that 
$A^{-8n(d(p))}C_1(1)f_S(j(p))$ is well defined for each 
triple point $p$ and no 2-face of the cube gives any relation for it.
$\Box$

\begin{remark}
Notice, that  $n(d(p))$ is of crucial importance in Lemma 18. This implies that Lemma 18 can not be extended to link 
diagrams, because $n(d(p))$ is no longer defined. 
\end{remark}

Theorem 2 follows now from Theorem 1, Lemmas 11, 12, 13, 16, 18 together with Remarks 9, 15.

\subsection{$S^+$ and $S^-$ as solutions of the cube equations}
We will use the autotangencies too for the cube equations. Let $s$ be a meridian of $\Sigma_{a-t}^{(2)}$. It contains exactly two 
diagrams with an autotangency. One easily sees that the simplifications $C_0$ of these two diagrams are regularly 
isotopic links if the transverse branch in the autotangency is either the highest or the lowest branch (with respect to the z-coordinate).
In this cas, the contributions of $C_0$ cancel out in $L(s)$. It follows that Lemma 16 can be still applied in each of the
two families of triple points (the positive and its relatives and the negative and its relatives). It remains to study the 
edges 1-6, 7-4, 2-8, 5-3.

\begin{lemma}
$(A^{-2} - A^2)C_0(1) = A^3C_1(1) + A^{-3}C_1(8)$, $(A^{-2} - A^2)C_0(3) = A^{-1}C_1(1) + A^1C_1(8)$, $C_1(1) = C_1(2)$,
$C_1(5) = C_1(7) = -A^2C_1(1)$, $C_1(8) = C_1(6)$, $C_1(3) = C_1(4) = -A^{-2}C_1(8)$ and $C_1(i) = -\bar C_1(i)$ for all $i$
 is the unique
solution of the cube equations with non trivial $C_0$.
\end{lemma}
(Lemma 7 determines $C_0(2)$ and $C_0(4)$.)

{\em Proof:\/} Lets consider the edge 7-1-4. We show the simplifications in Fig. 43.
\begin{figure}
\centering 
\psfig{file=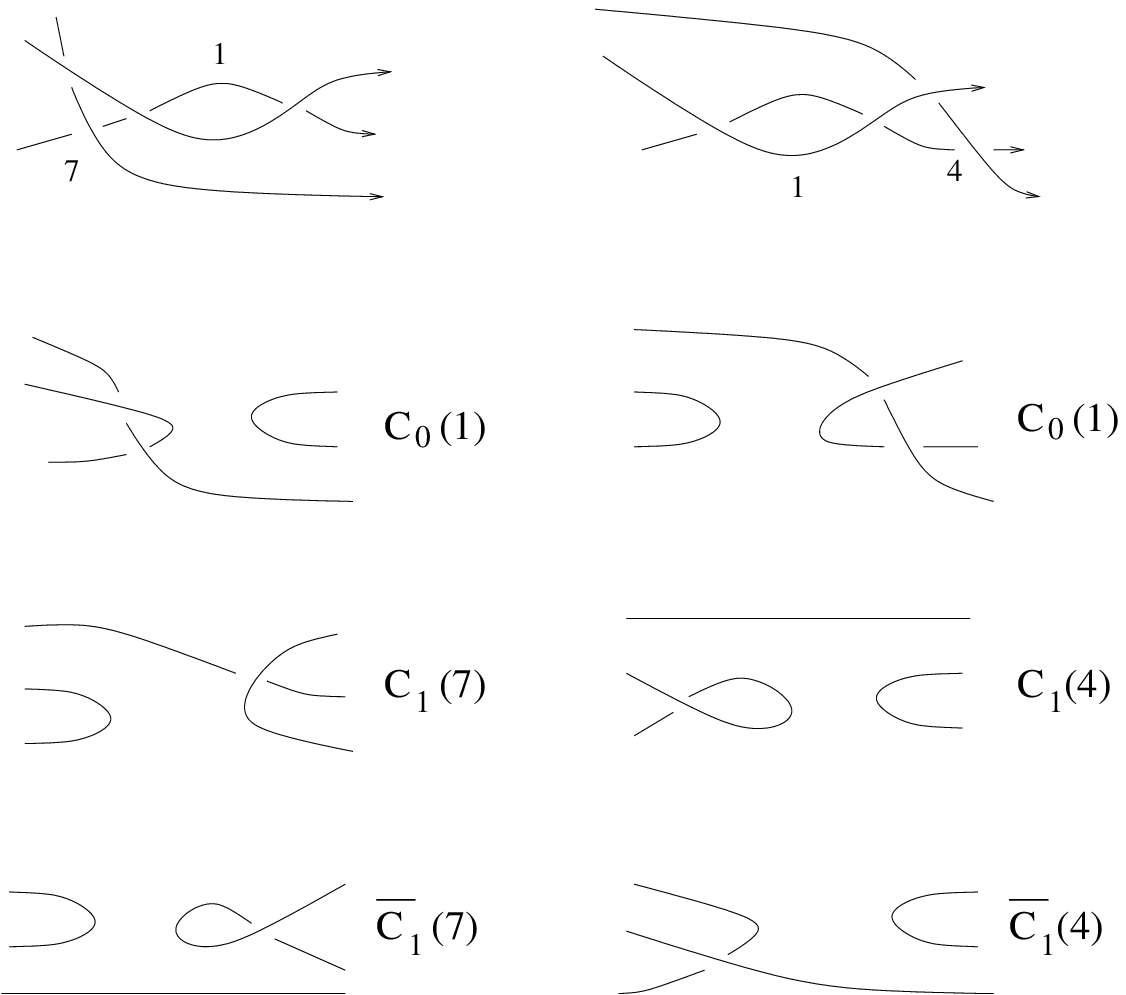}
\caption{}
\end{figure}
The calculation of the bracket gives us the following equations:

$x_1$:  $A^{-1}C_1(7)=-A^{-3}C_1(4)+(1-A^{-4})C_0(1)$

$\bar x_1$:  $(1 -A^4)C_0(1)-A^3\bar C_1(7) =A\bar C_1(4)$

$x_2$:  $AC_1(7)+A^{-2}C_0(1)=A^{-1}\bar C_1(4)+A^2C_0(1)$

Consequently, $(A^{-2} - A^2)C_0(1) = A^3C_1(1) + A^{-3}C_1(8)$ is a solution.

Lets consider the edge 2-3-8. We show the simplifications in Fig. 44.
\begin{figure}
\centering 
\psfig{file=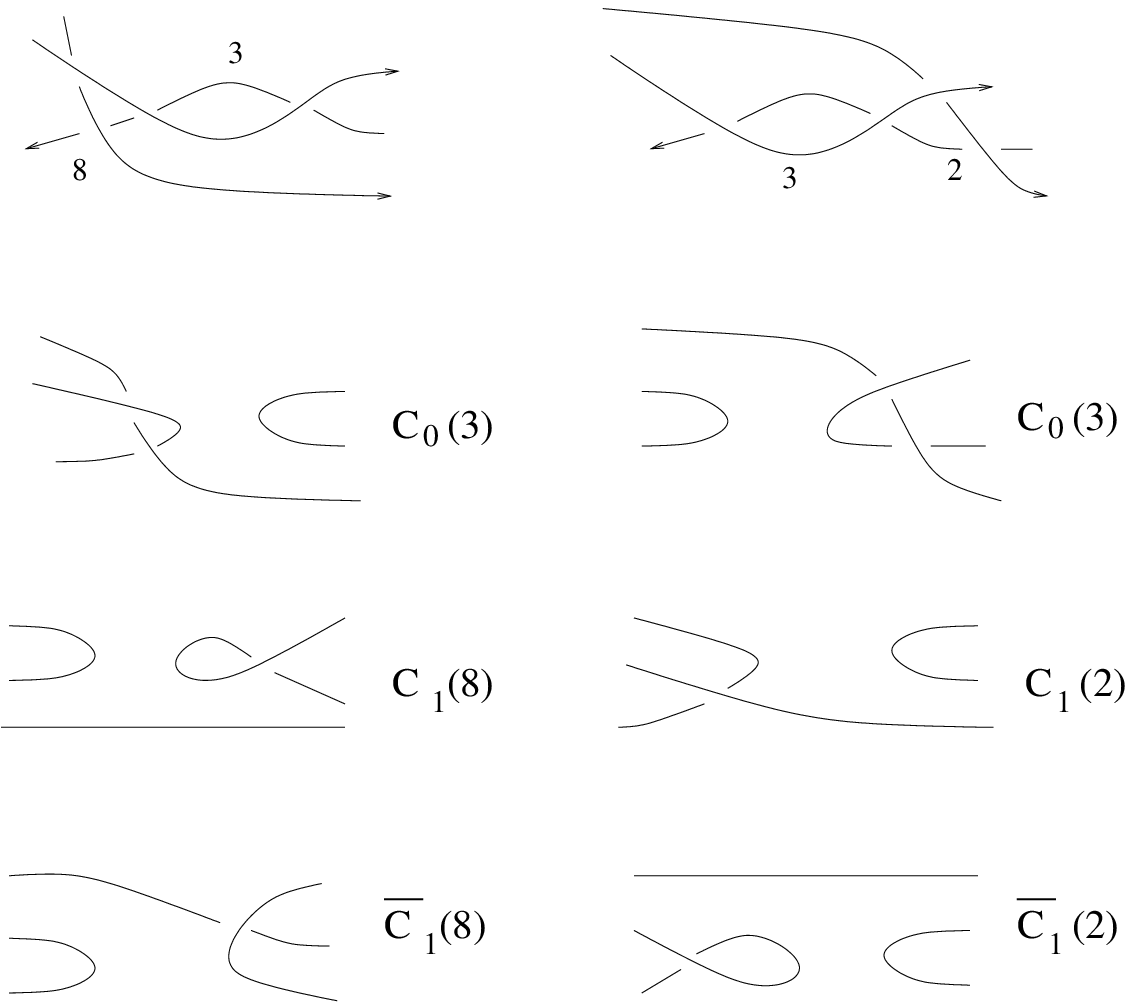}
\caption{}
\end{figure}
The calculation of the bracket gives us the following equations:

$x_1$: $A^{-1}\bar C_1(8) = -A^{-3}\bar C_1(2)+(1-A^{-4})C_0(3)$

$\bar x_1$:  $-A^3C_1(8)+(1-A^4)C_0(3)=AC_1(2)$

$x_2$:  $A\bar C_1(8)=A^{-1}C_1(2)$

Consequently, $(A^{-2} - A^2)C_0(3) = A^{-1}C_1(1) + A^1C_1(8)$ is a solution.

 All the other edges lead to the same 
solution. We left the verification to the reader. But notice, that it is essential to define the coorientation as it was done
in Fig. 8 in order that the edges 7-1-4 and 7-2-4, respectively 2-3-8 and 2-4-8,  have the same solutions.
$\Box$

Lemma 19 implies that the variables $C_1(1)$ and $C_1(8)$ become independent if we consider also the contributions
of autotangencies by $C_0$. Consequently, we can define $S^+$ and $S^-$ as done in Definition 14.
Let $s$ be a meridian of $\Sigma_q^{(2)}$ at a positive or negative quadruple point. Then $S(s) = S^+(s) = S^-(s) =0$.

Theorem 3 follows now from Theorem 1, Lemmas 11, 12, 13, 16, 19 together with Remarks 9, 15.

\subsection{$X$ as solution of the cube equations}
The homomorphism $X$ uses the Solution 2 of the positive tetrahedron equation. The calculations of
the cube equations go along the same lines as in the previous subsection but are very extensive. We just summarise 
the result of our calculations in the following lemma and we left the verification to the reader.

\begin{lemma}
$C_3(1) = -A^4C_4(1)$, $C_3(3) = C_3(4) = -C_4(1)$, $C_3(5) = C_3(7) = A^6C_4(1)$, $C_3(8) = A^2C_4(1)$,
$C_4(3) = C_4(4) = A^4C_4(1)$, $C_4(5) = C_4(7) =-A^2C_4(1)$, $C_4(8) = -A^6C_4(1)$, $C_5(j) = C_4(j)$ for each j,
$C_i(j) = -\bar C_i(j)$ for all i, j, $C_0(j) = 0$ for all j, $C_x(1) = A^5 + A$, $C_x(2) = -A^5 - A$, 
is the unique solution of the cube equations for closed braids (i.e. without triple points of type 2 and 6).

If in addition $C_3(2) = C_4(2) = C_5(2) = A^4C_4(1)$, $C_3(6) = C_4(6) = C_5(6) = A^2C_4(1)$, $C_x(3) = C_x(4) = 0$, then it becomes a 
solution for all cube equations but with $\mathbb{Z}/2\mathbb{Z}$ coefficients instead of integer coefficients.

(If we put in addition $A^4 = 1$ then evidently $C_x(j) = 0$ for all j.)

\end{lemma}

Theorem 4 follows now from Theorem 1 together with  Lemmas 9, 11, 12, 14, 20.
 
\subsection{$\Phi$ as solution of the cube equations}

We use the solution for the positive tetrahedron equation given in Definition 25.

Let $T_3$ be the free 
$\mathbb{Z}[A,A^{-1}]$-modul generated by the 15 {\em markings} with exactly two dots of the unoriented chord 
diagram $X_3$, but which has distinguished end points. Twelve of its generators correspond exactly to the markings 
defined in Fig. 20, but instead of $C_i$ or $\bar C_i$ we will call them $x_i$ respectively $\bar x_i$. The remaining 
three generators are defined in Fig. 45. 
\begin{figure}
\centering 
\psfig{file=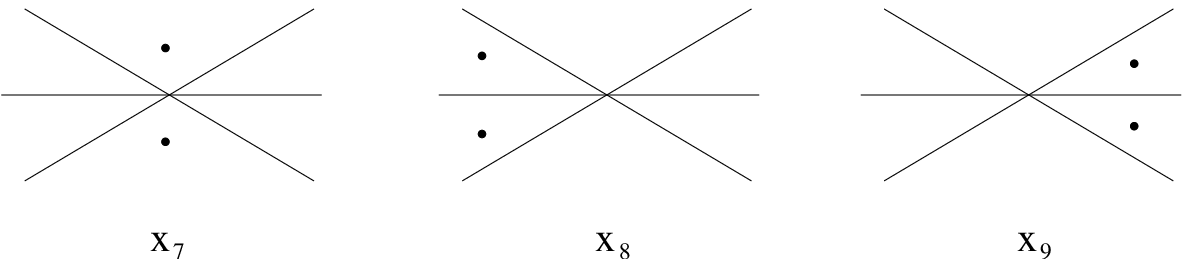}
\caption{}
\end{figure}
\begin{lemma}
$C_i(1) = C_i(2) = C_i(3) = C_i(4) = -C_i(5) = -C_i(6) = -C_i(7) = -C_i(8)$ for each $i \in \{2, 3, 4\}$,
$C_1(1) = C_1(3) = C_1(4) = -C_1(5) = -C_1(7) = -C-1(8) = (A + A^{-1})C_2(1)$, $C_1(2) = C_1(6) = 0$,
$C_i(j) = -\bar C_i(j)$ for all i,j, is a solution of the cube equations.
\end{lemma}

{\em Proof:\/} Each edge which connects two braid-like triple points, say j and j', gives the relations $C_i(j) = -C_i(j')$
for all i. Edges which connect a braid-like triple point with a star-like one give in addition the (surprising) relation
$C_1(2) = C_1(6) = 0$. Calculations are again very extensive. We consider just one example. All other cases are 
similar and are left to the reader.

The markings for the two triple points corresponding to the edge 1-1-7 are shown in Fig. 46 - 47.
\begin{figure}
\centering 
\psfig{file=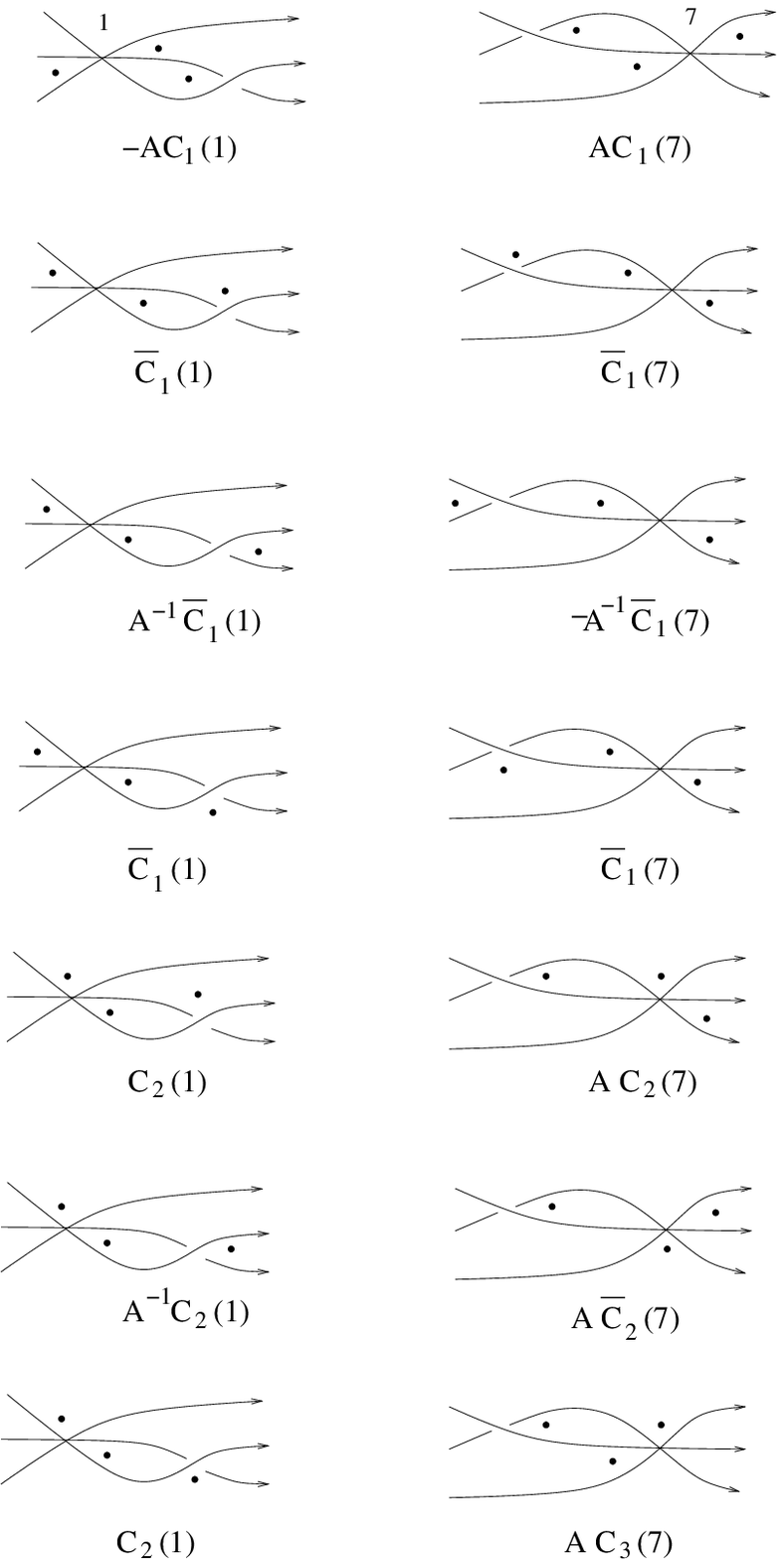}
\caption{}
\end{figure}
\begin{figure}
\centering 
\psfig{file=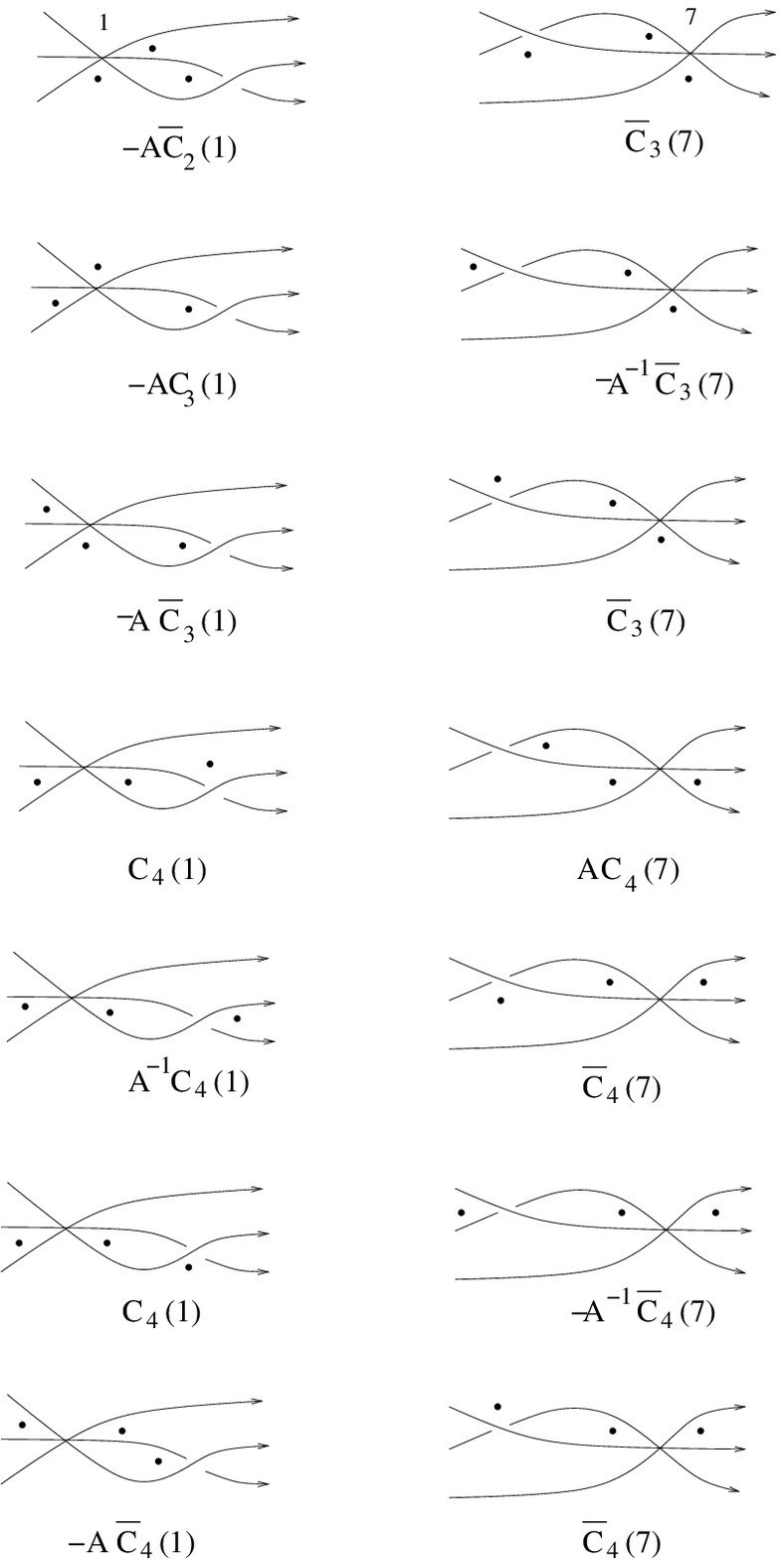}
\caption{}
\end{figure}
 We obtain in $T_3$

$t(p) =-AC_1x_1+\bar C_1\bar x_4+A^{-1}\bar C_1\bar x_1+\bar C_1\bar x_3+C_2x_5+A^{-1}C_2x_2+C_2x_7
-A\bar C_2\bar x_2-AC_3x_3-A\bar C_3\bar x_3+C_4x_1+A^{-1}C_4x_4+C_4\bar x_6-A\bar C_4\bar x_4 $

$t(p') =AC_1(7)x_1+\bar C_1(7)x_2-A^{-1}\bar C_1(7)\bar x_1+\bar C_1(7)x_4+AC_2(7)x_2+A\bar C_2(7)\bar x_2+AC_3(7)x_3
+\bar C_3(7)\bar x_6-A^{-1}\bar C_3(7)\bar x_3+\bar C_3(7)x_7+AC_4(7)x_4+\bar C_4(7)x_1-A^{-1}\bar C_4(7)\bar x_4+\bar C_4(7)x_5 $.

Comparing the coefficients of the $x_i$ gives the desired relations.
$\Box$

Theorem 5 follows now from Theorem 1 together with  Lemmas 11, 12, 15, 21.

We have not checked if there is a solution 
which uses in addition the contributions of autotangencies (compare Lemma 8).

\subsection{All the homomorphisms are trivial on  sliding classes}
Let $sl$ be a sliding loop (compare Definition 2).

\begin{lemma}
$S([sl]) = 0$ for any sliding loop $sl$.
\end{lemma} 

{\em Proof:\/} There are again many cases to distinguish, we will carry out just one of them. Let us consider a curl in 
positive direction with positive writh which passes twice through a crossing of a knot diagram, as shown in Fig. 48.
\begin{figure}
\centering 
\psfig{file=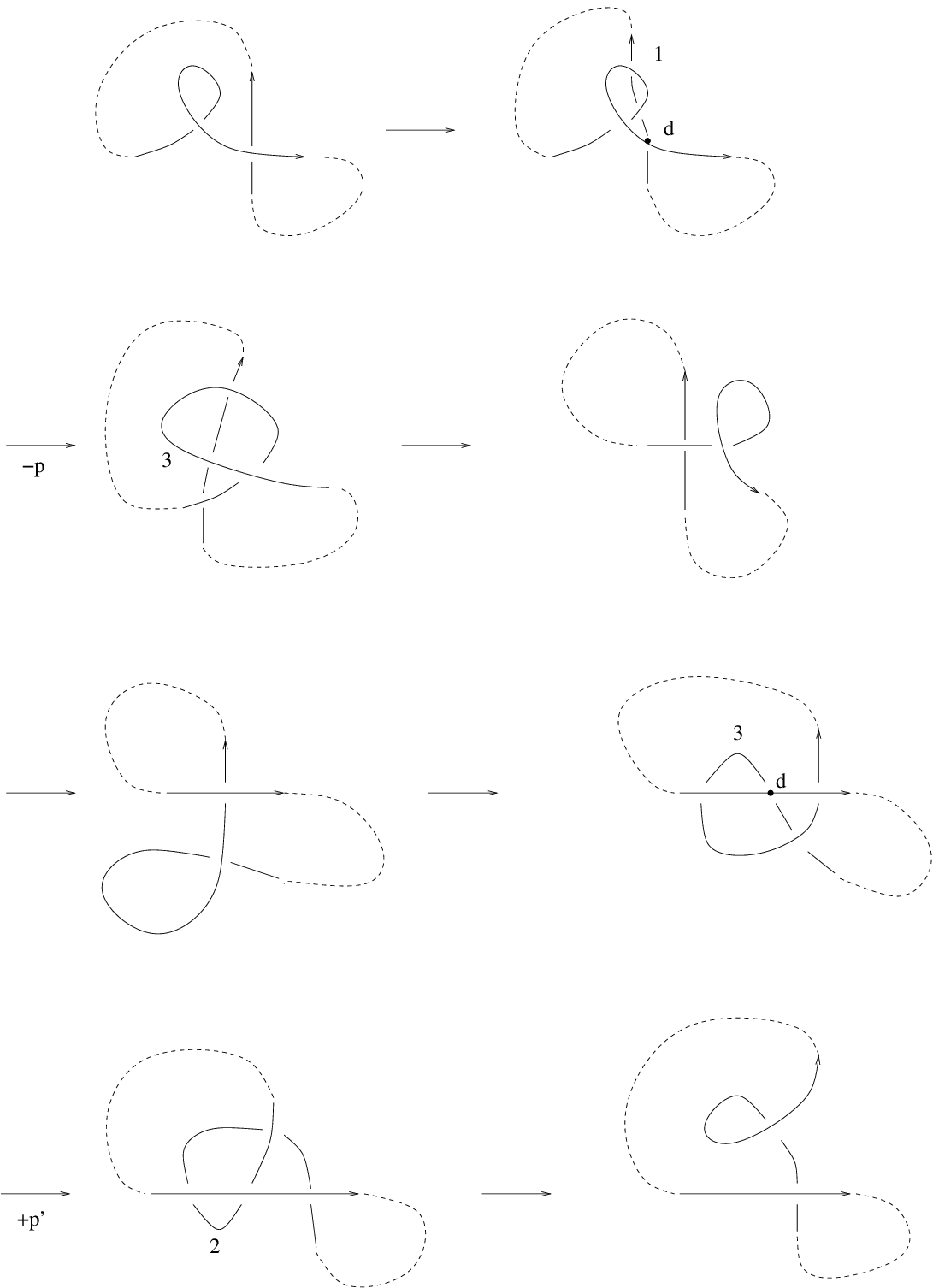}
\caption{}
\end{figure}
We show the contributions to $S([sl])$ of the two triple points in Fig. 49. Notice, that $[d(p)] = [d(p')]$ and $n(d(p)) = n(d(p'))$.
It follows that the contributions in $S$ of $p$ and $p'$ cancel out.
\begin{figure}
\centering 
\psfig{file=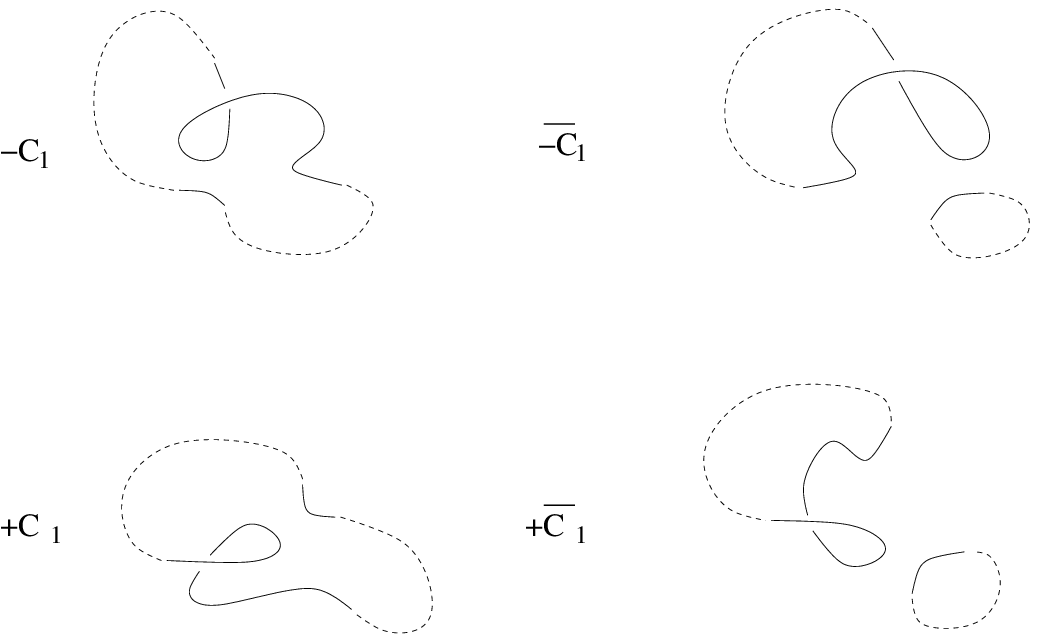}
\caption{}
\end{figure}
One easily sees, that $p$ and $p'$ are always of the same type besides in the following two cases: if one is of type 3 
then  the other is of type 4 and if one is of type 5 then the other is of type 7. But all our invariants never distinguish the 
types in these two couples.
$\Box$

\begin{lemma}
$S^+([sl]) = 0$ and $S^-([sl]) = 0$.

\end{lemma}
{\em Proof:\/} We know already from the previous lemma that the contributions of the triple points cancel out.
There are exactly four autotangencies (compare Fig. 48). 
So, we have only to verify that the contributions of the autotangencies cancel out too. We do this in Fig. 50.
(Remember that $C_0(1) = - C_0(2)$.)
\begin{figure}
\centering 
\psfig{file=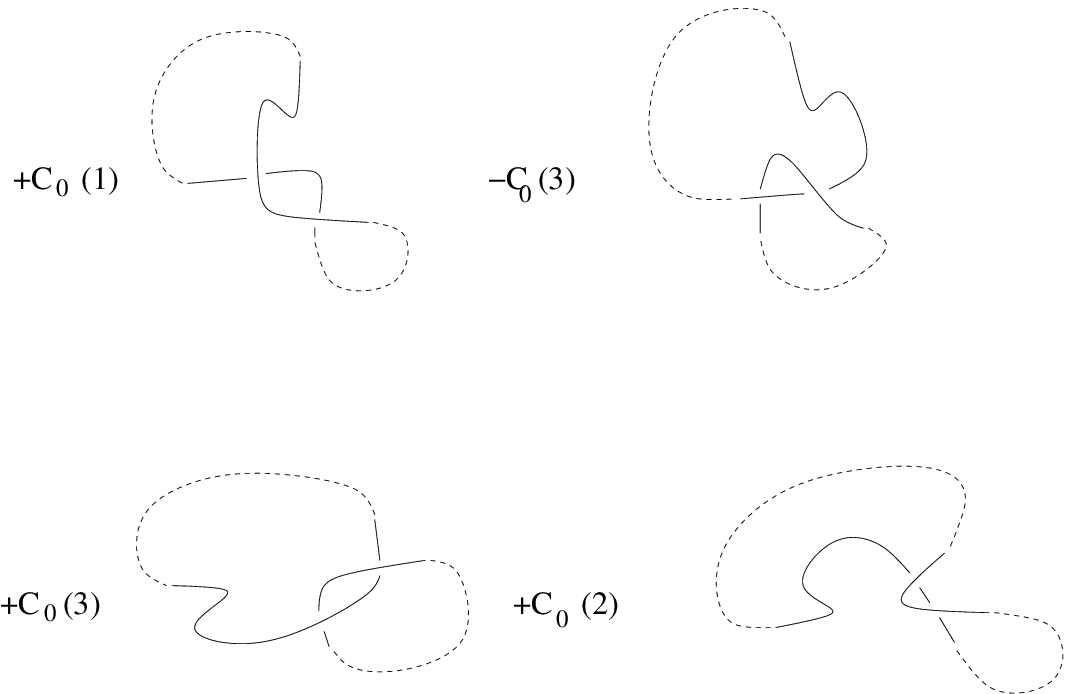}
\caption{}
\end{figure}
$\Box$

\begin{lemma}
$X([sl]) = 0$.
\end{lemma}
{\em Proof:\/} We consider the same case as previously. The contributions of the triple point $p$ and of two
autotangencies is shown in Fig. 51. 
\begin{figure}
\centering 
\psfig{file=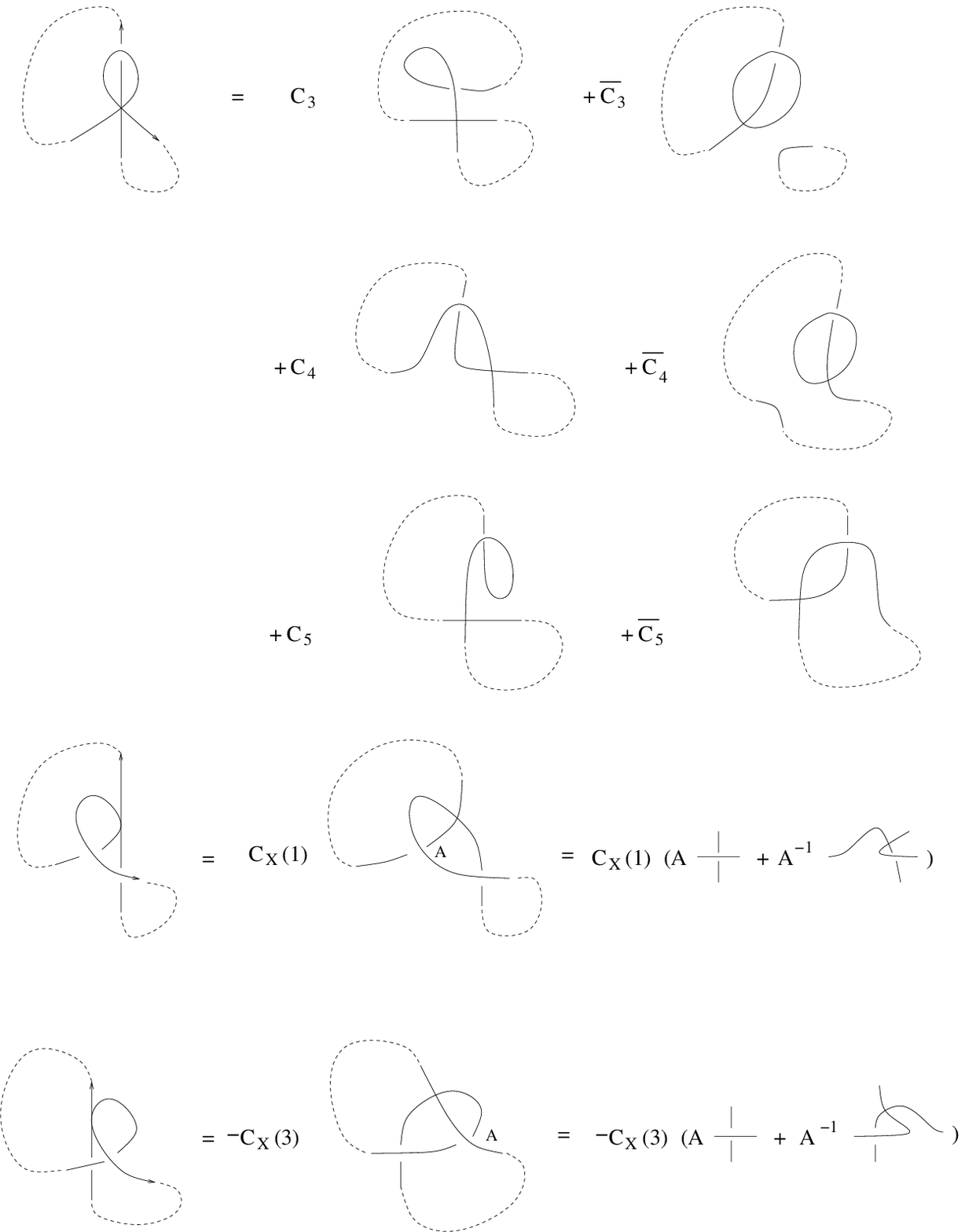}
\caption{}
\end{figure}
The contribution of $p'$ and of the other two autotangencies is shown in Fig. 52.
\begin{figure}
\centering 
\psfig{file=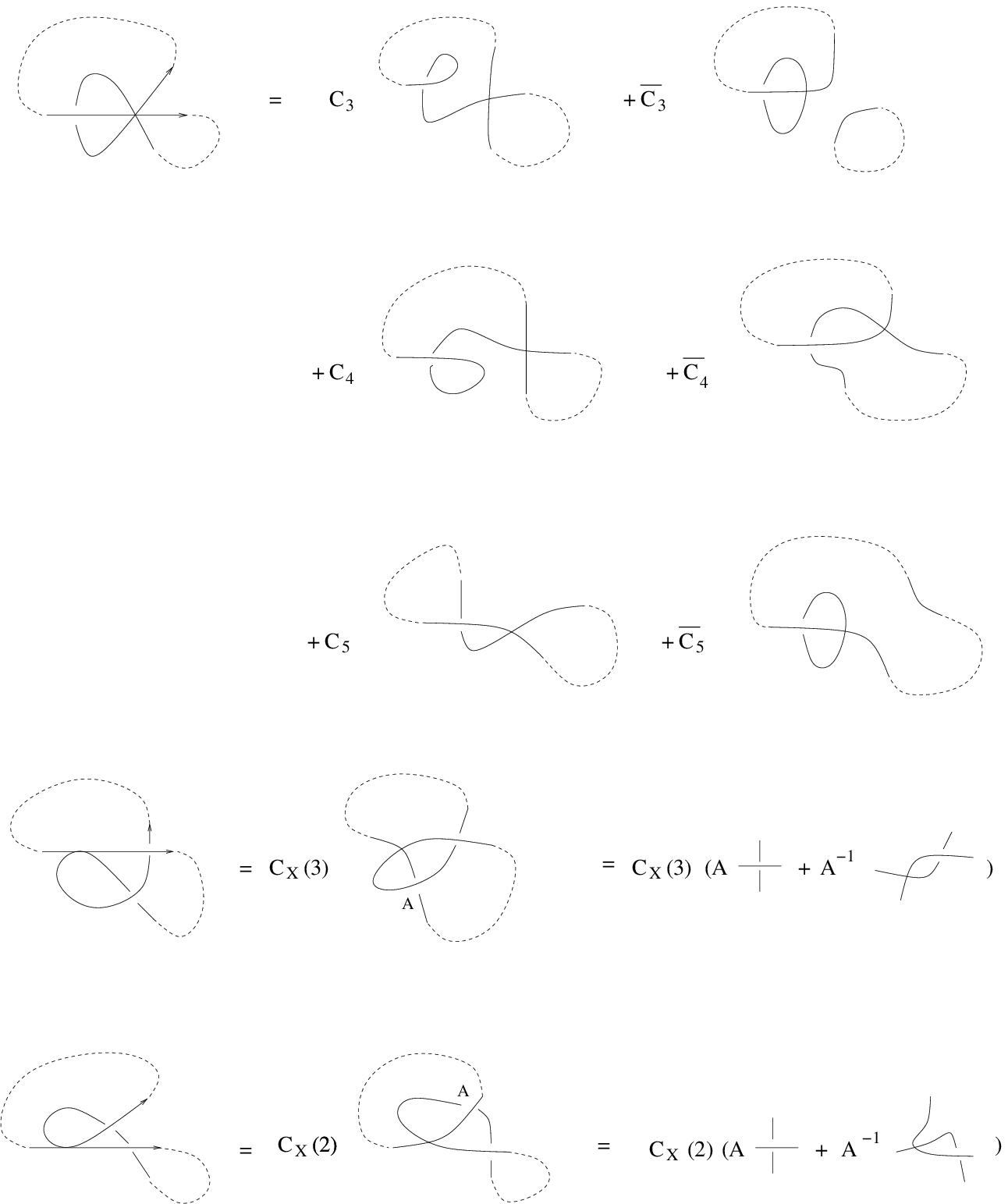}
\caption{}
\end{figure}
One easily sees that they cancel out again (remember that $p$ and $p'$ have different signs) .
$\Box$

\begin{lemma}
$\Phi([sl]) = 0$.
\end{lemma}
{\em Proof:\/} The contributions of $p$ and $p'$ are shown in Fig. 53 respectively Fig. 54.
\begin{figure}
\centering 
\psfig{file=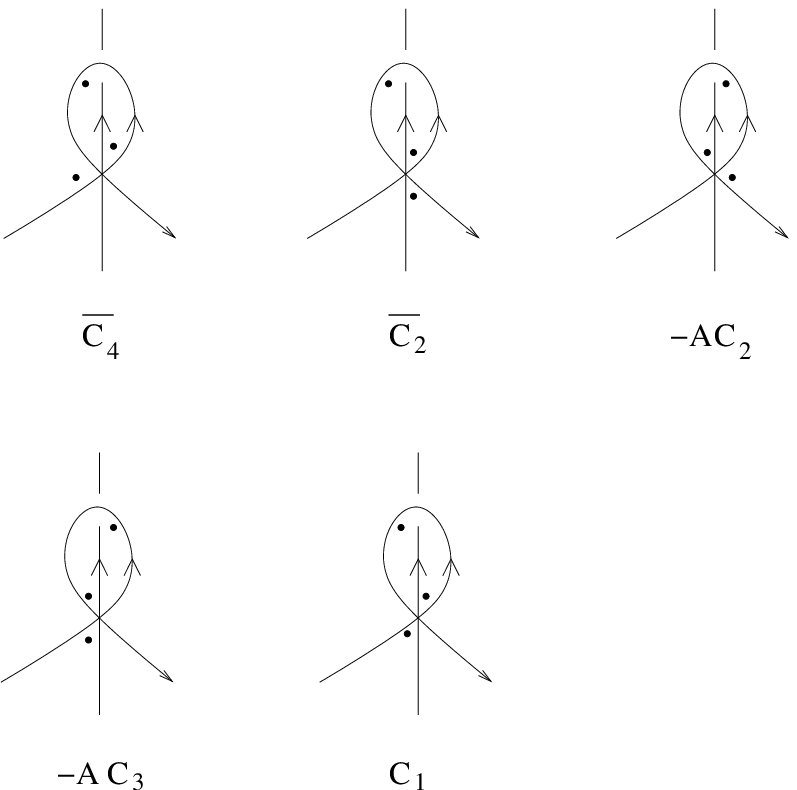}
\caption{}
\end{figure}
\begin{figure}
\centering 
\psfig{file=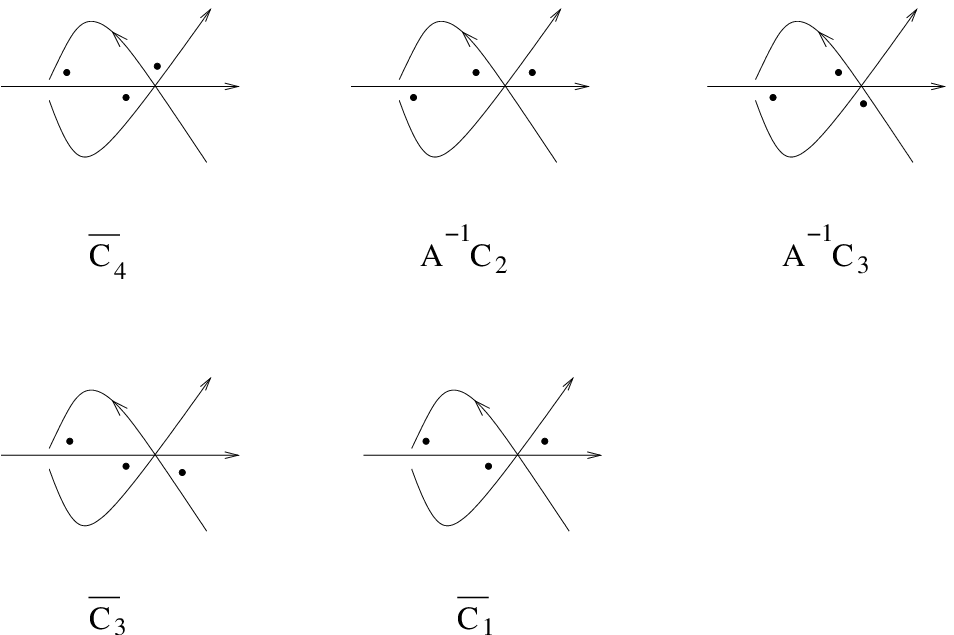}
\caption{}
\end{figure}
 Each of the contributions 
corresponds to 
a single marking for the chord diagram $X_2$ (compare Subsection 2.7.). Consequently, we obtain exactly four 
equations:

$-AC_2 -A^{-1}C_2 -\bar C_1 = 0$

$\bar C_2 -\bar C_3 = 0$

$-AC_3 +C_1 -A^{-1}C_3 = 0$

$\bar C_4 -\bar C_4 = 0$

It follows from the definition of $\Phi$ that these equations are always satisfied.
$\Box$

The fact, that for all our homomorphisms $\rho$ we have $\rho([sl]) = 0$, is sufficient to garanty that $\rho([rot(K)]$ is
a knot invariant. Indeed, an isotopy of $K$ in $V$ induces a homotopy of the pre-canonical loop $rot(K)$ in $M'$ (compare
 Sect. 2.2). We have to study how the approximation of $rot(K)$ in $M$ changes under such a homotopy of the pre-canonical
loop. To this end we have just to consider those local moves of the trace graph $TG$ from Fig. 10 in \cite{FK}, which 
involve a cusp in the knot projection. The autotangencies correspond to black dots in Fig. 10 in \cite{FK}. Consequently,
we have only to check invariance under passing a {\em ramphoid cusp} (Fig. 10v) and under the {\em extreme pair 
move} (Fig. 10x). 

The changing of the approximation of $rot(K)$ under passing a ramphoid cusp is just part of a Whitney trick.

In Fig. 55  we reproduce Fig. 10x  together with its interpretation for the knot isotopy.

\begin{figure}
\centering 
\psfig{file=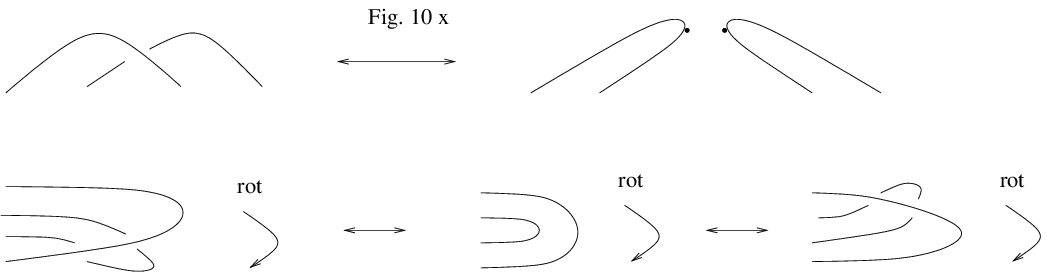}
\caption{}
\end{figure}

After passing an extreme pair (i.e. $K$ touches from the same side the same disc of the disc fibration of $V$
at two different points, compare \cite{FK}) two new autotangencies appear in the approximating loop.
But Fig. 55 shows that these two autotangencies are in fact the same but with opposit signs. Consequently, 
$\rho([rot(K)]$ stays invariant.

It remains to study the effect of Whitney tricks. Let us call a Whitney trick {\em positive}, if it adds two new crossings 
to the diagram, and {\em negative} otherwise. Evidently, Whitney tricks are local, i.e. they change $\rho([rot(K)]$ in exactly
the same way if they are performed at different places in the knot diagram. It remains to notice, that the number of positive 
Whitney tricks in the approximation of $rot(K)$ is equal to the number of singularities of $\phi$ restricted to $K$, i.e.
the number of times that the knot $K$ is tangent to the disc fibration of $V$ (compare Sect. 2.1). (Remember that each
such tangency produces two cusps in the knot projections, one for each orientation of the tangent line.)
In the end of the approximation of $rot(K)$ we performe only negative Whitney tricks in order to eleminate all the small
curls which were created before. Hence, all the small curls which were produced by positive Whitney tricks, have to disappear.
Consequently, the number of negative Whitney tricks is equal to the number of positive Whitney tricks in the approximation.
It follows that $\rho([rot(K)]$ does not depend on the approximation of $rot(K)$ modulo sliding loops.

\subsection{Open questions}
We formulate eight questions which could perhaps indicate possible directions for the further developement of the subject.

\begin{question}
Write a computer program in order to calculate the invariants in examples. Can they detect mutation?
Can $\Phi_K$ detect non invertibility?
\end{question}

\begin{question}
How do the invariants change under Reidemeister I moves? How do they change under Markov moves (see e.g. \cite{B})?
\end{question}

\begin{question}
Can Jaegers (more complicated) state model for the HOMFLY-PT polynomial (see \cite{J} and also \cite{K3}) be used in the same way
 in order to define homomorphisms? A natural try is to replace state models by skein relations. We have carried this out by  using
Kauffman's polynomial for unoriented framed links. However, it did not work ! In order to obtain a solution of the 
tetrahedron equation it seems to be essential that in the state model all permutations induced by the tangles are different.
(Switching a crossing does not change the permutation induced by the branches of the tangle.)
\end{question}

\begin{question}
Are the new invariants related to some new representation theory? Let $M_{col}^{pure}$ be the space of all colored pure positive
closed braids.
Does the operator equation (b) (compare the Introduction) has a non trivial solution in the case of colored pure positive closed
braids? (One easily sees that the matrices associated to autotangencies would cancel out in a meridian of $\Sigma^{(2)}_{a-t}$.
Consequently, the matrice associated to a triple point could not depend on the type of the triple point. This forces us
to consider only positive triple points and hence, only isotopies through positive closed braids.) 
If such a solution exists, then we could take the product of these matrices along the canonical loop. In this way we 
would produce for each component of $M_{col}^{pure}$ a matrix 
(over some
 commutative ring) which is well defined up to conjugacy. The relations of these matrices for different components of 
$M_{col}^{pure}$ would be much more complicated as those for the matrices obtained from the representation theory 
related to the 
Yang-Baxter equation.

More general, it seems to be tempting to associate to each triple point a matrix in a local way, e.g. to the triple point 
$\sigma_i\sigma_{i+1}\sigma_i = \sigma_{i+1}\sigma_i\sigma_{i+1}$ we associate a matrix $A_i^{sign}$. However,
a closer look to Fig. 14 shows that then the tetrahedron equation would become an identity. But our invariants are non 
trivial solutions of a tetrahedron equation !

In general, the algebra which could be behind our invariants stays rather mysterious (for me).
\end{question}

\begin{question}
Our homomorphisms can be considered as one dimensional cocycles on $M$ which are defined in a combinatorial way
by using the discriminant $\Sigma$. The infinite dimensional space $M$ has a natural approximation by finite 
dimensional smooth manifolds 
(possibly with singularities). Do there exist  "closed polynomial valued 1- forms" on $M$ which define the same 
cohomology classes as our 1-cocycles?
\end{question}

\begin{question}
As well known, knot polynomials (coming from quantum knot invariants) can be decomposed into series of finite type
knot invariants (see \cite{BN} and references therein). Do our homomorphisms have an analogue decomposition, i.e. 
into series of integer valued (respectively, 
$\mathbb{Z}/2\mathbb{Z}$ valued) one 
cocycles, such that the value of each of them is calculable in polynomial time with respect to the number of 
intersections of the 
loop $\gamma$ with the discriminant $\Sigma^{(1)}$ and of the number of crossings of the diagrams at these 
intersections
(this is a perhaps more general definition of {\em finite type one cocycles}
as that of Vassiliev, see \cite{V} and references therein)?
\end{question}

\begin{question}
Let us consider the homomorphisms $S, S^+$ or $S^-$.
Each smoothing of a triple point or an autotangency leads to an ordinary (non oriented) link diagram in the annulus.
It is not hard to see (but we will not carry this out), that we can define an orientation in a canonical way on diagrams 
which are obtained by {\em smoothings}
of triple points or autotangencies of oriented knot diagrams $D$.
Consequently, its Khovanov homology is well defined (see \cite{Kh} and \cite{APS}) and it is a 
relative knot invariant 
of the oriented diagram $D$ with respect to the triple point or the autotangency (compare the Introduction). Let us try to proceed as in the 
construction of cellular homology from singular homology. We take as chain groups the direct sums of all the Khovanov 
homology groups of all smoothings of all triple points and autotangencies. There are natural gradings on these groups.
Do there exist  differentials such that the corresponding homology groups of the complex depend only on the homology
class of the loop $\gamma$ in $M$ and such that the graded Euler characteristics are equal to $S([\gamma])$,
$S^+([\gamma])$ or $S^-([\gamma])$?

In the case of $X$ and $\Phi$ the situation is even more complicated. The simplifications or markings of triple points 
and autotangencies do not correspond to ordinary link diagrams. Here the first question is the following: can the
definition of Khovanov homology, respectively knot Floer homology (see  \cite{MOS}), be extended in order to give 
relative knot invariants of the oriented diagram $D$ with respect to the triple point or the autotangency? 
\end{question}

\begin{question}
Can our approach be generalized to produce new knot polynomials from higher dimensional families of knots?
If a knot diagram in the solid torus is not a satellite of another non trivial knot diagram (other than the core of
the solid torus) , in particular it has no curls, then very often the corresponding component of $M$ retracts by deformation onto 
a 2-torus (compare the Introduction). But this 2-torus does not intersect $\Sigma^{(2)}$ at all ! 

However, we can replace
the diagram by a fixed satellite. The product of the canonical loops is now a 2-torus which intersects $\Sigma^{(2)}_{q}$.
A meridional 2-sphere $S$ of a stratum of $\Sigma^{(3)}$, which corresponds to an ordinary quintupel point, cuts 
$\Sigma^{(2)}_{q}$ transversally in ten points. We would have to define coorientations for the strata of $\Sigma^{(2)}_{q}$
and to define extensions of the Kauffman state sums for diagrams with an ordinary quadruple point in the projection.
The analogue of the tetrahedron equation would then be that the signed contributions of the ten quadruple points in $S$
cancel out. But the calculations are so complex that one would need to carry them out with a computer.

Another possibility would be to consider the intersections of the 2-torus with the strata which are transverse 
intersections in $\Sigma^{(1)}_t \cap \Sigma^{(1)}_t$, i.e. diagrams which have two triple crossings at the same time.
But notice, that  there is no natural coorientation on such strata. However, there is a canonical coorientation on
the strata in $\Sigma^{(1)}_t \cap \Sigma^{(1)}_a$.

Finally, we could go on increasing the dimension by taking satellits of satellits.
\end{question}

Laboratoire de Math\'ematiques

Emile Picard

Universit\'e Paul Sabatier

118 ,route de Narbonne 

31062 Toulouse Cedex 09, France

fiedler@picard.ups-tlse.fr

\end{document}